\documentclass[11pt]{article}
\usepackage[T1]{fontenc}
\usepackage[utf8]{inputenc}
\usepackage{lmodern}
\usepackage[a4paper,margin=29mm,headheight=14pt]{geometry}
\usepackage{amsmath,amssymb,amsthm,mathtools}
\usepackage{microtype,enumitem}
\usepackage{needspace}
\usepackage{placeins}
\usepackage{array,booktabs,tabularx}
\usepackage[dvipsnames]{xcolor}
\usepackage{tikz}
\usetikzlibrary{arrows.meta,calc,positioning}
\definecolor{ringblue}{RGB}{36,83,122}
\definecolor{ringred}{RGB}{158,62,48}
\usepackage[colorlinks=true,linkcolor=MidnightBlue,citecolor=MidnightBlue,urlcolor=MidnightBlue]{hyperref}

\setlist[enumerate,1]{label=\textup{(\roman*)},ref=\roman*,leftmargin=2em,itemsep=3pt,topsep=5pt}
\numberwithin{equation}{section}
\newtheorem{theorem}{Theorem}[section]
\newtheorem{proposition}[theorem]{Proposition}
\newtheorem{lemma}[theorem]{Lemma}

\theoremstyle{definition}

\theoremstyle{plain}
\newtheorem{setup}{Setup}[section]
\newtheorem{remark}[theorem]{Remark}
\newtheoremstyle{requirement}{\topsep}{\topsep}{\itshape}{}
 {\bfseries}{.}{.5em}
 {\thmname{#1}\thmnumber{ #2}\thmnote{ (#3)}}
\theoremstyle{requirement}
\newtheorem{condition}{Condition}
\newtheorem{parametercondition}{Condition}

\theoremstyle{plain}
\DeclareRobustCommand{\conditionref}[1]{Condition~\ref{#1} (\nameref{#1})}
\newcommand{\R}{\mathbb R}
\newcommand{\Hh}{\mathbb H}
\newcommand{\BS}{\operatorname{BS}}
\newcommand{\supp}{\operatorname{supp}}
\newcommand{\Div}{\operatorname{div}}
\newcommand{\curl}{\operatorname{curl}}
\newcommand{\diag}{\operatorname{diag}}
\newcommand{\dd}{\,\mathrm d}
\newcommand{\dnu}{\,\mathrm d\nu}
\newcommand{\eps}{\varepsilon}
\newcommand{\cE}{\mathcal E}

\newcommand{\vectomega}{\boldsymbol\omega}

\allowdisplaybreaks[1]
\hypersetup{
 pdftitle={Sharp Finite-Time Blow-Up for Unforced Axisymmetric Euler Without Swirl via Multiscale Data},
 pdfsubject={Research draft of a multiscale construction},
 pdfauthor={Jeaheang Bang and Alexey Cheskidov},
 pdfkeywords={Euler equations, axisymmetric flow, no swirl, Holder regularity, finite-time singularity}
}

\title{\bfseries Sharp Finite-Time Blow-Up for Unforced\\Axisymmetric Euler Without Swirl\\via Multiscale Data}
\author{Jeaheang Bang\thanks{Institute for Theoretical Sciences, Westlake Institute for Advanced Study, Westlake University, 600 Dunyu Road, Hangzhou, Zhejiang 310024, People's Republic of China. Email: \texttt{jhbang@westlake.edu.cn}.}
\and Alexey Cheskidov\thanks{Institute for Theoretical Sciences, Westlake University, 600 Dunyu Road, Hangzhou, Zhejiang 310024, People's Republic of China. Email: \texttt{cheskidov@westlake.edu.cn}.}}
\date{}

\begin{document}
\maketitle

\begin{abstract}
For the unforced Euler equations, we construct a compactly
supported initial vorticity within the class of axisymmetric flows
without swirl, whose solution blows up in finite time. The initial
vorticity belongs simultaneously to every H\"older space
$C^\alpha$, $\alpha<1/3$. This is sharp given that global
regularity is known when $\alpha=1/3$ in the case of axisymmetric
flows without swirl.

Previous constructions by Shkoller \cite{Shk26}, Chen \cite{Chen26I,Chen26II}, and
Shao--Wei--Zhang--Zhang \cite{SWZZ26} establish blow-up for each prescribed exponent
$\alpha$ in the same sharp range. Our construction achieves
simultaneous membership in all these spaces. In addition, our construction is short.

Our proof uses multiscale initial data, in the spirit of
Kim--Jeong~\cite{KJ22}, C\'ordoba--Mart\'inez-Zoroa--Zheng~\cite{CMZZ25}, and
Bang--Cheskidov~\cite{BC26}. The blow-up mechanism is a forward
cascade through successively smaller scales. A main difficulty in approaching the critical exponent $1/3$ is caused by large self-deformation. Our main idea to overcome this difficulty is to obtain a lower bound for strain at the stationary non-smooth point by using angular redistribution.

A related outstanding open problem is to establish an unforced Euler
blow-up from smooth data in the case of axisymmetric flows with swirl.
A recent multiscale construction by OpenAI for the unforced Euler equations is,
however, non-axisymmetric. Although we work in the no-swirl class,
the axisymmetric multiscale mechanism developed here may offer a
possible route toward the still-open problem of axisymmetric smooth
blow-up.
\end{abstract}

\tableofcontents

\section{Introduction}\label{sec:intro}

We consider the incompressible Euler equations
\begin{equation}\label{eq:euler}
 \partial_tu+u\cdot\nabla u+\nabla p=0,\qquad
 \Div u=0,\qquad u|_{t=0}=u_0
\end{equation}
on $\R^3$, where $u(t,x)\in\R^3$ is the velocity and
$p(t,x)\in\R$ is the pressure. The evolution of the vorticity
$\vectomega=\curl u$ is governed by
\[
 \partial_t\vectomega+u\cdot\nabla\vectomega
 =\vectomega\cdot\nabla u.
\]
The interaction between transport, vortex stretching, and the nonlocal
recovery of velocity from vorticity is central to the regularity theory.
In this paper we study finite-time singularity formation in the
three-dimensional axisymmetric class without swirl.

\subsection{Classical local theory and continuation}

The local theory in H\"older spaces goes back to 
\cite{Lic25}, \cite{Gun27}. In the finite-energy setting
used here, a divergence-free datum
$u_0\in C^{1,\alpha}(\R^3)\cap L^2(\R^3)$,
$0<\alpha<1$, generates a unique local classical solution.
More generally, $C^{k,\alpha}$ regularity persists locally for
integers $k\geq1$. The later Sobolev theory includes the work of
Ebin and Marsden \cite{EM70} and Kato \cite{Kat72}.
We refer to \cite{MB01,BCD11} for modern accounts of these theories.

The continuation criterion of Beale, Kato, and Majda \cite{BKM84}
identifies vorticity amplification as the obstruction to extending a
classical solution. In the H\"older setting above, a finite maximal
lifespan $T_*$ necessarily satisfies
\[
 \int_0^{T_*}\|\vectomega(t)\|_{L^\infty(\R^3)}\dd t=\infty.
\]
See
\cite[Theorem~4.2.3]{Che95}.
In two dimensions,
transport without vortex stretching conserves the supremum norm of vorticity  and yields global
regularity. Axisymmetry without swirl retains vortex stretching but provides
an additional transported quantity.

\subsection{Axisymmetry and the global regularity threshold}

Let $\mathbf{e}_1 = (1,0,0)$, $\mathbf{e}_2 = (0,1,0)$, and
$\mathbf{e}_3 = (0,0,1)$. In cylindrical coordinates,
\[
 x = (r\cos\theta,r\sin\theta,z), \qquad \mathbf{e}_r = (\cos\theta,\sin\theta,0), \qquad \mathbf{e}_\theta = (-\sin\theta,\cos\theta,0).
\]
An axisymmetric velocity without swirl is of the form
\[
 u(t,x) = u^r(t,r,z)\mathbf{e}_r+u^z(t,r,z)\mathbf{e}_3.
\]
Its vector vorticity $\vectomega = \curl u$ is purely azimuthal. With the sign convention used in the construction, the
scalar and relative-vorticity variables are
\[
 \vectomega=\omega\mathbf e_\theta,\qquad q=-\frac{\omega}{r}.
\]
Here $\omega$ follows the scalar-vorticity convention of \cite{BC26},
while $q$ is the negative of their relative vorticity. 
The sign is chosen so that $q\geq0$ in the upper half-space for our data.
The Euler equations imply
\begin{equation}\label{eq:relative-transport}
 \partial_tq+u^r\partial_rq+u^z\partial_zq = 0.
\end{equation}
Although $q$ is transported, $\omega = -rq$ changes
along trajectories as their distance from the axis changes. Moreover,
H\"older continuity of $\vectomega$ does not imply boundedness of
$q$. This distinction allows singularity formation in the class
considered here.

Global regularity for smooth, sufficiently decaying axisymmetric
Euler flows without swirl goes back to \cite{UY68}, with subsequent refinements by
\cite{SR94}, \cite{SY94}.

Danchin \cite{Dan07} extended the global theory to bounded vorticity
with suitable decay and relative vorticity in $L^{3,1}(\R^3)$.
For finite-energy axisymmetric no-swirl data with compactly supported
vorticity, this yields global classical solutions when
$u_0\in C^{1,\alpha}(\R^3)$, $1/3<\alpha<1$.
The result of Shao, Wei, and Z.~Zhang \cite{SWZ26} also covers the
endpoint $C^{1/3}$ vorticity class under the same finite-energy
and compact-support assumptions.

Thus the possible blow-up range for the compact-vorticity H\"older
data considered here lies strictly below $1/3$.

\subsection{Singularity formation below the threshold}

Elgindi \cite{Elg21} constructed self-similar axisymmetric no-swirl
blow-up solutions with $C^{1,\alpha}$ velocity for sufficiently small
$\alpha>0$.
Elgindi, Ghoul, and Masmoudi \cite{EGM21} developed stability of the
self-similar construction and finite-energy localization.

C\'ordoba, Martinez-Zoroa, and Zheng \cite{CMZZ25} in 2025 introduced a
multiscale construction based on infinitely many interacting vorticity regions
separated by vorticity-free regions, without using a self-similar construction. They obtained finite-energy blow-up
solutions with velocity in
$C^{1,\alpha}(\mathbb{R}^3)\cap C^\infty(\R^3\setminus\{0\})$ for a small
positive exponent $\alpha$. The present paper extends their multiscale construction to the sharp range $\alpha \in (0,1/3)$.

Recent works reach every fixed $0<\alpha<1/3$. Chen
\cite{Chen26I,Chen26II} constructs self-similar profiles from a
one-dimensional limiting model and proves asymptotically self-similar
blow-up with compactly supported $C^\alpha$ vorticity and finite
energy. Shao, Wei, P.~Zhang, and Z.~Zhang \cite{SWZZ26} construct
exact self-similar profiles in dimension $d\geq3$ for
$0<\alpha<1-2/d$; finite-codimensional stability yields finite-energy,
asymptotically self-similar blow-up with compact initial vorticity and
$C^{1,\alpha}$ velocity smooth away from the origin. Shkoller
\cite{Shk26} provides a Lagrangian construction based on strain rates and pressure for each
$0<\alpha<1/3$. This raises the question addressed here: can one
initial datum belong to every supercritical H\"older class and
still develop a singularity?

Table~\ref{tab:sharpcomparison} compares these sharp H\"older blow-up
results. 

\subsection{The multiscale method and smooth-data blow-up}\label{sec:multiscalehistory}

Multiscale constructions use multi-scale data for which large scales 
 stretch and amplify smaller scales successively. Kim and Jeong \cite{KJ22} used
axisymmetric multiscale vorticity rings to prove critical
Sobolev and Lorentz ill-posedness. C\'ordoba, Mart\'inez-Zoroa,
and Zheng \cite{CMZZ25} constructed H\"older-class blow-up from
infinitely many interacting vorticity regions, working backward-in-time
from terminal configurations in a small H\"older exponent regime;
see Remark~\ref{rem:forwardconstruction}.

In \cite{BC26}, we developed a forward-in-time ring cascade that accounts
for geometric self-slowdown of the strain rate $\partial_r u^r(t,0,0)$, proving vorticity norm inflation
and instantaneous blow-up for bounded vorticity with relative
vorticity in $L^{3,s}$, for every $s>1$, extending a Lorentz result by Kim--Jeong \cite{KJ22} above from large $s$ to every $s>1$. The present paper
establishes finite-time blow-up in H\"older spaces forward-in-time and extends the multiscale construction in \cite{CMZZ25} to the sharp H\"older range; see Remark~\ref{rem:selfdeformation} for comparison.

Outside the axisymmetric case without swirl, smooth data may develop a finite-time blow up. One may consider either axisymmetric flows with swirl or non-axisymmetric flows.
C\'ordoba and Mart\'inez-Zoroa \cite{CMZ23} constructed
non-axisymmetric forced Euler blow-up with force uniformly bounded
in $C^{1,1/2-\varepsilon}\cap L^2$.
Alp\"oge and Buckmaster's work-in-progress note \cite{AB26} obtained smooth-data
axisymmetric blow-up with a force smooth in space and
time through blow-up.
OpenAI \cite{OAIE26} announced unforced Euler blow-up from smooth,
compactly supported initial velocity. Their non-axisymmetric
construction successively amplifies localized oscillations and
corrects them to exact Euler solutions, whose gradients amplify
the next oscillations.
Related forced multiscale constructions include IPM \cite{CMZIPM24},
Boussinesq \cite{CLM25}, and Navier--Stokes \cite{CMZZ26NS,OAIN26},
with fractional dissipation in \cite{CMZZ26NS}.

The outstanding goal is smooth-data blow-up for \emph{unforced
axisymmetric Euler on $\R^3$}. The smooth-data results of
Chen--Hou \cite{CH25} involve a solid boundary, while the
whole-space axisymmetric result \cite{AB26} uses forcing and
OpenAI's unforced construction is non-axisymmetric.  Although our datum is not smooth at the origin and our
solution has no swirl, the unforced forward-in-scale cascade developed here,
including its control of large self-deformation (explained below in Remark~\ref{rem:selfdeformation} and in Table~\ref{tab:multiscalecomparison}), may provide a
route toward that goal. Table~\ref{tab:multiscalecomparison} compares
all these multiscale constructions. For more multiscale constructions on the Euler equations, see \cite{Choi26, HQSZ26, CZ26}.

\subsection{Main results}

Our initial datum belongs to every supercritical H\"older class and consists of sparse anisotropic vortex rings. 
The first theorem states the main blow-up result. 
The second one specifies our multiscale datum more explicitly. 

\begin{theorem}[Finite time blow-up]\label{thm:main}
There exists a divergence-free axisymmetric velocity without swirl
\[
 u_0\in L^2(\R^3)\cap\bigcap_{0<\alpha<1/3}C^{1,\alpha}(\R^3) \cap C^\infty(\R^3\setminus\{0\}),
\]
with compactly supported vorticity $\vectomega_0 = \curl u_0$ such that
the classical solution $u(t)$ starting from $u_0$ satisfies
\[
 \int_0^{T_*}\|\vectomega(t)\|_{L^\infty(\R^3)}\dd t = \infty, \qquad \vectomega = \curl u,
\]
for some finite $T_*>0$.
\end{theorem}

\begin{theorem}[Explicit multiscale family]\label{thm:constructive}
Fix a profile
$f_0\in C_c^\infty((0,\infty)\times(0,\infty))$ with
\[
 \int_0^\infty\!\int_0^\infty
 \rho f_0(\rho,\zeta)\dd\rho\dd\zeta=1, \qquad f_0  \geq 0.
\]
Let $\widetilde f_0$ be its odd extension in the second variable,
with $\widetilde f_0(\rho,0)=0$. There exists a constant
$Y_0>0$ depending only on $f_0$ such that, for every $Y\geq Y_0$,
the initial data defined below satisfies all the conclusions of
Theorem~\ref{thm:main}:
\begin{equation}\label{eq:constructive-data}
 \begin{gathered}
  \vectomega_0=\sum_{j=1}^\infty\vectomega_{j,0},\qquad
 u_0=\BS[\vectomega_0]
     =\sum_{j=1}^\infty\BS[\vectomega_{j,0}],\\
 \vectomega_{j,0}(x)=q_{j,0}(r,z)(x_2,-x_1,0), \qquad q_{j,0}(r,z)=Q_{*,j}\widetilde f_0(r/R_{j,0},z/H_{j,0}),
 \end{gathered}
\end{equation}
where
{
\begin{equation}\label{eq:constructive-scales}
 \begin{gathered}
 R_{j,0}=e^{-L_j},\qquad H_{j,0}=e^{-5\widehat L_j/6},
 \qquad Q_{*,j}=e^{(1-\alpha_j)L_j},\\
 \alpha_j=\frac{j}{3(j+1)},\qquad \beta_j=1-3\alpha_j=\frac{1}{j+1},
 \\
 L_j=\frac{Y\,100^{j-1}}{\displaystyle\prod_{i=1}^j\beta_i}
 =Y 100^{j-1} (j+1)!
 ,
 \qquad \widehat L_j=\beta_jL_j=Y100^{j-1}j!.
 \end{gathered}
\end{equation}
}
\end{theorem}

Here $\BS[\vectomega_0]$ is the Biot-Savart law defined in \eqref{eq:bs}.

\begin{remark}
We emphasize that our single initial datum belongs to all supercritical H\"older spaces simultaneously, that is,
\[
u_0 \in \bigcap _{0<\alpha<1/3} C^{1,\alpha} (\mathbb{R}^3).
\] On the other hand, in the previous sharp blow-up constructions \cite{Shk26, Chen26I, Chen26II, SWZZ26}, the initial data belongs to a H\"older space $C^{1,\alpha}$  for a prescribed exponent $\alpha\in(0,1/3)$.
\end{remark}

\begin{remark}[Backward vs forward constructions]\label{rem:forwardconstruction}
C\'ordoba, Mart\'inez-Zoroa, and Zheng \cite[Section~3]{CMZZ25}
prescribe multiscale terminal configurations, i.e., at the singularity at the time of blow-up, and determine the initial
data by backward evolution. We instead prescribe the initial datum
explicitly in Theorem~\ref{thm:constructive} and follow its rings
forward in time. The terminal configuration is determined by the evolution.
See Table~\ref{tab:multiscalecomparison} for comparisons of multiscale constructions.
\end{remark}

\begin{remark}[Large self-deformation]\label{rem:selfdeformation}
Self-deformation is the deformation caused by a vortex-ring pair's
own induced velocity. Earlier multiscale arguments \cite{KJ22,CMZZ25,BC26}
suppress this contribution, while allowing large stretching and flattening driven by outer rings.

In particular, the construction in \cite{CMZZ25} keeps self-deformation perturbative, which allows the dynamics to be approximated by a system of ODEs driven by outer strain rates. The parameter choices that ensure this control also restrict the initial regularity to a small positive H\"older exponent.

Our construction approaches the sharp H\"older threshold $1/3$ without imposing smallness of the self-deformation. Establishing a blow up through the exact PDE without relying on the approximating ODEs is therefore a central difficulty in our proof. See Section~1.6 on our main ideas to overcome this difficulty.

 See Remark~\ref{rem:CMZZansatz} for another aspect on comparison with \cite{CMZZ25}, and see Table~\ref{tab:multiscalecomparison} for comparisons of multiscale constructions.
\end{remark}

\begin{remark}[Freedom of the initial profile]\label{rem:profilefreedom}
Theorem~\ref{thm:constructive} allows any initial profile $f_0$
satisfying its stated assumptions. The same profile is used at every
scale, and $Y_0$ may depend on it. Flexibility in the initial profiles
is also noted in \cite[Remark~6]{CMZZ25}, subject to support and norm
bounds. This freedom within the prescribed multiscale family does not
assert stability under general perturbations of the initial datum.
\end{remark}

\begin{remark}[Smoothness away from the origin]
\label{rem:blowuplocation}
The solution of Theorem~\ref{thm:constructive} satisfies
\[
 u\in C^\infty\bigl((0,T_*)\times(\R^3\setminus\{0\})\bigr).
\]
Vorticity blows up at the origin: for every $\delta>0$,
\[
 \sup_{0\leq t<T_*}
 \|\vectomega(t)\|_{L^\infty(\R^3\setminus B_\delta(0))}<\infty,
 \qquad
 \int_0^{T_*}\|\vectomega(t)\|_{L^\infty(B_\delta(0))}\dd t=\infty.
\]
The proof is given in Subsection~\ref{subsec:blowuplocation-proof}.
\end{remark}

{

\begin{remark}[Global conservative weak continuation]
\label{rem:weakcontinuation}
In a subsequent paper we will show that the solution of
Theorem~\ref{thm:constructive} admits a global axisymmetric weak
Euler continuation $U$ without swirl, with $U(t)=u(t)$ for $0\leq t<T_*$, and
\[
 \begin{gathered}
 U\in C_{\mathrm{loc}}\bigl([0,\infty);
 H^1(\R^3)\cap L^\infty(\R^3)\bigr),\\
 U\in C^1_{\mathrm{loc}}\bigl([0,\infty);L^2(\R^3)\bigr).
 \end{gathered}
\]
Moreover,
\[
 \|U(t)\|_{L^2}=\|u_0\|_{L^2} \qquad t\geq 0.
\]
In particular, we have
\[
 \lim_{t\to T_*}\|u(t)-U(T_*)\|_{H^1}=0,\qquad
 \|U(T_*)\|_{L^2}=\|u_0\|_{L^2}.
\]
Such a continuation can be obtained as a subsequential
vanishing-viscosity limit of axisymmetric Navier--Stokes
solutions without swirl.
\end{remark}

}

\begin{figure}[!t]
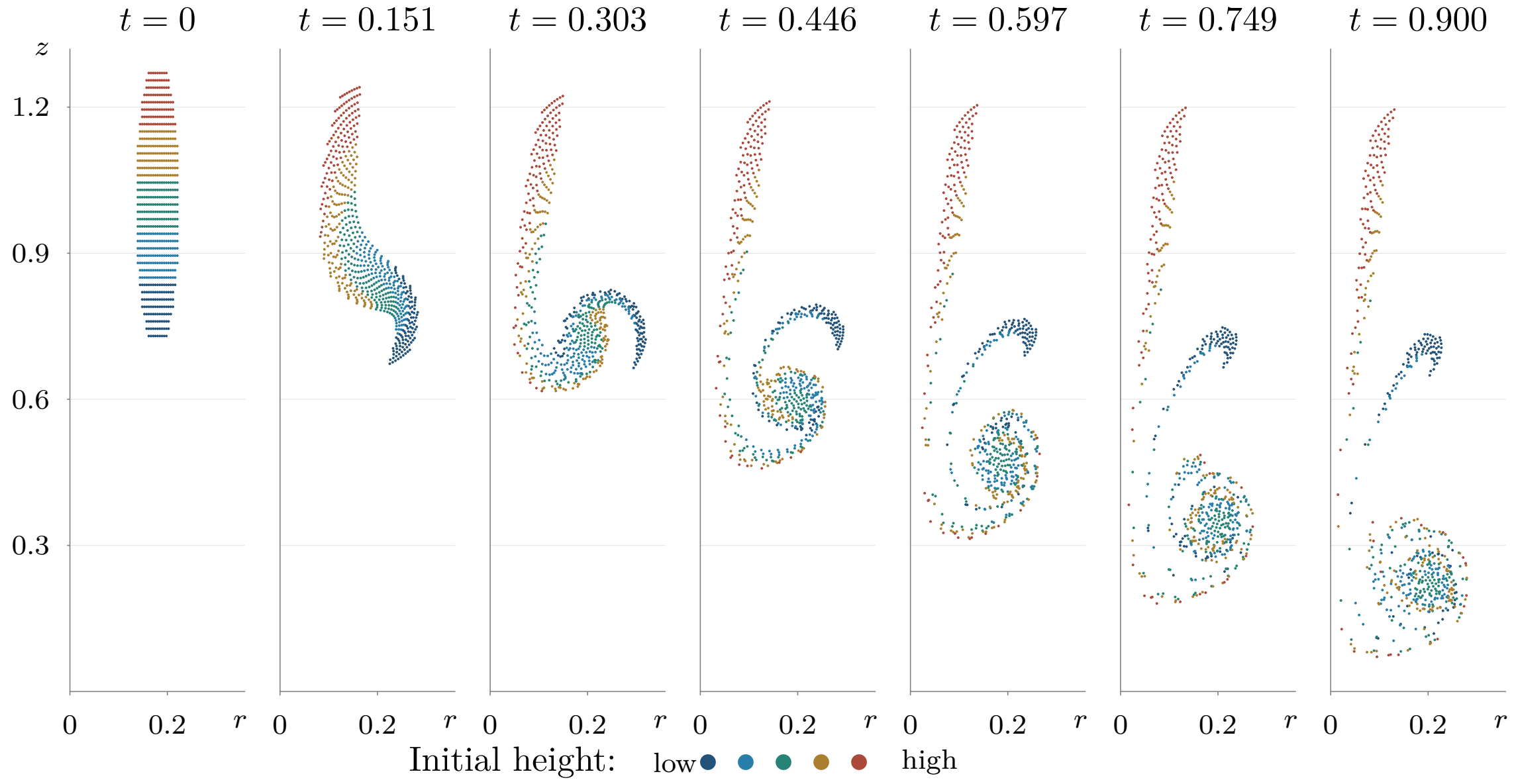

\centering
\begingroup
\definecolor{paircolor0}{HTML}{24537A}
\definecolor{paircolor1}{HTML}{287EA8}
\definecolor{paircolor2}{HTML}{278375}
\definecolor{paircolor3}{HTML}{AA802E}
\definecolor{paircolor4}{HTML}{AB493B}
\resizebox{\linewidth}{!}{%
%
}
\endgroup
\caption{Numerical self-deformation of one vortex-ring pair.
The panels show positions of numerical vortex markers in the first
meridional quadrant, with equal radial and vertical scales and the same
axes in every panel. Colors label initial heights and remain attached
to the markers; they do not represent vorticity amplitude.
The upper ring is sampled by 603 markers at
$(r,z)=(0.18+0.045\hat r,1+0.30\hat z)$, using the points of a
$19\times37$ uniform grid in $[-0.9,0.9]^2$ that satisfy
$\hat r^2+\hat z^2<1$. Each marker carries a fixed positive weight
approximating the integral of $-\omega=rq$ over a small meridional
cross-section. This integral is conserved when the cross-section is
transported by the Euler flow. The weights determine each marker's
contribution to the Biot--Savart velocity that transports the markers.
They are proportional to $(1-\hat r^2-\hat z^2)^2$ and normalized to
sum to one. The oppositely signed reflected ring below $z=0$ is
included in the dynamics but omitted from view.
We evaluate the axisymmetric Biot--Savart velocity with
the kernel factor $|x-y|^{-3}$ replaced by
$(|x-y|^2+\delta^2)^{-3/2}$, where $\delta=0.025$, using
complete elliptic integrals for the angular integration and explicit
midpoint time stepping with step $0.0005$.
All quantities are dimensionless. These are numerical trajectories
in a fixed regularized model, not schematic deformations or a
reconstruction of the support boundary. Convergence as the marker
count increases and $\delta\to0$ is not asserted.}
\label{fig:selfdeformationnumerics}
\end{figure}

Figure~\ref{fig:selfdeformationnumerics} illustrates the geometric
difficulty through a numerical evolution of an isolated vortex-ring
pair (with no imposed external field). The initially elongated
distribution bends and spreads radially, while markers from different
initial heights rearrange within it. Even in this simple configuration,
the evolving geometry is substantially more complicated than a rescaling
of the initial profile. Our proof controls moments, {strain at the origin}, and support
bounds without requiring a detailed description of this deformation.

\begin{table}[!t]
\centering
\caption{Selected multiscale Euler constructions and the smooth-data
axisymmetric target. }
\label{tab:multiscalecomparison}
\begingroup
\small
\setlength{\tabcolsep}{4pt}
\renewcommand{\arraystretch}{1.18}
\begin{tabularx}{\linewidth}{@{}
 >{\raggedright\arraybackslash}p{.15\linewidth}
 >{\raggedright\arraybackslash}p{.18\linewidth}
 >{\raggedright\arraybackslash}p{.26\linewidth}
 >{\raggedright\arraybackslash}p{.21\linewidth}
 >{\raggedright\arraybackslash}X@{}}
\toprule
Work & Time organization & Self-deformation / self-interaction
 & Symmetry & External force \\
\midrule
\multicolumn{5}{@{}l}{\textit{Critical-space local ill-posedness}}\\[3pt]
Kim--Jeong \cite{KJ22} & Forward & Small 
 & Axisymmetric, no swirl & None \\
\addlinespace[4pt]
Bang--Cheskidov \cite{BC26} & Forward & Small
 & Axisymmetric, no swirl & None \\
\midrule
\multicolumn{5}{@{}l}{\textit{Axisymmetric no-swirl blow-up in H\"older classes}}\\[3pt]
CMZZ \cite{CMZZ25} & Backward & Small 
 & Axisymmetric, no swirl & None \\
\addlinespace[4pt]
This paper & Forward & Large self-deformation allowed
 & Axisymmetric, no swirl & None \\
\midrule
\multicolumn{5}{@{}l}{\textit{Higher-regularity and smooth-data blow-up}}\\[3pt]
CM \cite{CMZ23} & Backward
 & Small; residual canceled by forcing
 
 & Non-axisymmetric & Present \\
\addlinespace[4pt]
Alp\"oge--Buckmaster \cite{AB26} & Forward
 & Higher-order corrections; residual forced
 & Axisymmetric, with swirl & Present \\
\addlinespace[4pt]
OpenAI \cite{OAIE26} & Forward
 & Small
 & Non-axisymmetric & None \\
\midrule
Smooth-data open problem & Not applicable & Not applicable
 & Axisymmetric, with swirl & None \\
\bottomrule
\end{tabularx}

\smallskip
\begin{minipage}{\linewidth}
\footnotesize
CM denotes C\'ordoba--Mart\'inez-Zoroa; CMZZ also includes Zheng.

\end{minipage}
\endgroup
\end{table}

Table~\ref{tab:sharpcomparison} compares the sharp H\"older blow-up
results. The present construction combines one datum in every
$C^{1,\alpha}$, $0<\alpha<1/3$, with smoothness away from the
origin and an explicit ring series built from any profile allowed by
Theorem~\ref{thm:constructive}. Smoothness away from the origin is
also a feature of \cite{SWZZ26}. Flexibility in the initial profile
is shared with the construction in \cite{CMZZ25}; see
Remark~\ref{rem:profilefreedom}.

\begin{table}[!t]
\centering
\caption{Sharp H\"older blow-up constructions for unforced
axisymmetric Euler without swirl on $\R^3$. All four yield
finite-energy solutions. Throughout, $0<\alpha<1/3$.}
\label{tab:sharpcomparison}
\vspace{4pt}
\begingroup
\small
\setlength{\tabcolsep}{4pt}
\renewcommand{\arraystretch}{1.18}
\begin{tabularx}{\linewidth}{@{}
>{\raggedright\arraybackslash}p{.17\linewidth}
*{4}{>{\raggedright\arraybackslash}X}@{}}
\toprule
Feature & Shkoller \cite{Shk26} & Chen \cite{Chen26I,Chen26II}
& Shao--Wei--Zhang--Zhang \cite{SWZZ26} & This paper \\
\midrule
H\"older regularity of $u_0$
& $C^{1,\alpha}$ for each prescribed $\alpha$
& $C^{1,\alpha}$ for each prescribed $\alpha$
& $C^{1,\alpha}$ for each prescribed $\alpha$
& One datum in every $C^{1,\alpha}$ simultaneously \\
\addlinespace[5pt]
Initial smoothness
& Axis cusp
& Axis cusp
& $C^\infty$ off the origin
& $C^\infty$ off the origin \\
\addlinespace[5pt]
Initial vorticity support
& Noncompact, with decay
& Compact
& Compact
& Compact \\
\midrule
Construction
& Lagrangian analysis of strain rates and pressure
& Self-similar profiles lifted from a 1D model
& Self-similar elliptic--transport profiles
& Forward multiscale ring cascade \\
\addlinespace[5pt]
Computer-assisted step
& None
& 1D profile in Part~I
& None
& None \\
\midrule
Blow-up description
& Type-I; Lagrangian collapse
& Type-I; asymptotically self-similar
& Type-I; asymptotically self-similar
& Finite lifespan; no prescribed rate \\
\addlinespace[5pt]
Stability / profile freedom
& Admissible angular perturbations
& Finite-codimensional profile stability
& Finite-codimensional profile stability
& Any admissible ring profile \\
\bottomrule
\end{tabularx}

\medskip
\begin{minipage}{\linewidth}
\footnotesize
Type-I here means
$\|\vectomega(t)\|_{L^\infty(\R^3)}\asymp(T-t)^{-1}$.
The Chen and Shao--Wei--Zhang--Zhang columns refer to their
finite-energy localized solutions. The latter result also holds in
every dimension $d\geq3$, with threshold $1-2/d$.

\smallskip
For the data and stability statements, see
\cite[Theorems~1.1--1.2 and Definition~5.1]{Shk26},
\cite[Theorem~1.2 and Remarks~1.3--1.4]{Chen26II}, and
\cite[Theorem~1.3 and Remark~1.2]{SWZZ26}.
Here profile freedom refers to the family in
Theorem~\ref{thm:constructive}; See Remark~\ref{rem:profilefreedom}. General
stability under initial-data perturbations is not asserted.
\end{minipage}
\endgroup
\end{table}

We call the vorticity described by $q_j$ a \emph{vortex-ring pair}:
its restriction to $z>0$ is the upper vortex ring, and its odd
reflection is the oppositely signed lower ring. We use this terminology
also for their transported vorticities. A vortex ring may have a
disconnected meridional cross-section; the term imposes neither a
fixed cross-sectional shape nor a thin-core assumption. When the
index is omitted, $q$ denotes one such pair.

To describe the geometry of our initial datum, which is also illustrated in
Figure~\ref{fig:strainmechanism}, enclose $\supp f_0$ in a rectangle
$[\rho_-,\rho_+]\times[\zeta_-,\zeta_+]$ with
$\rho_-,\zeta_->0$. The $j$-th upper vortex ring is supported in
\[
 \rho_-R_{j,0}\leq r\leq\rho_+R_{j,0},\qquad
 \zeta_-H_{j,0}\leq z\leq\zeta_+H_{j,0},
\]
and its lower ring is its odd reflection.
{For $Y\geq Y_0$, these support regions are pairwise disjoint
and accumulate only at the origin. The exponents $\alpha_j$ increase to the threshold
$1/3$, and the regularity gaps $\beta_j=1-3\alpha_j$ tend to zero.
The scales in Theorem~\ref{thm:constructive} satisfy
\[
 \frac{R_{j,0}}{H_{j,0}}
 =e^{-(1-5\beta_j/6)L_j}\to0,\qquad
 Q_{*,j}R_{j,0}=R_{j,0}^{\alpha_j},\qquad
 R_{j,0}=H_{j,0}^{6/(5\beta_j)}.
\]
Thus the distances to the symmetry axis shrink faster than the heights.
The datum consists of increasingly slender, vertically elongated
vortex-ring pairs with vorticity
amplitudes comparable to $R_{j,0}^{\alpha_j}$.  The motivations behind the choices of these sequences in \eqref{eq:constructive-scales} are explained in Section~\ref{sec:cascade}. For comparison with \cite{KJ22,BC26}, see Remark~\ref{rem:relative_scale}.

The proofs of
Theorems~\ref{thm:main} and~\ref{thm:constructive} are written in
Section~\ref{sec:limit}.}

{
\subsection{Idea of the construction}\label{sec:constructionidea}

\begin{figure}[!t]
\centering
\begingroup
\color{black}
\newcommand{\ep}{\varepsilon}
\newcommand{\largeint}[2]{\int_{#1}^{#2}}
\newcommand{\stepineq}[1]{\mathrel{\overset{\text{\raisebox{2pt}{%
 \fontsize{10}{12}\selectfont (#1)}}}{\gtrsim}}}
\resizebox{\linewidth}{!}{%
\begin{tikzpicture}[x=1cm,y=1cm,>=Stealth,text=black,
 font=\fontsize{10}{12}\selectfont,
 panel/.style={draw=black,line width=.65pt,rounded corners=2pt,fill=white},
 heading/.style={anchor=west,text=black,font=\fontsize{12}{14}\selectfont\bfseries},
 note/.style={anchor=west,font=\fontsize{10}{12}\selectfont,text=black},
 math/.style={anchor=west,font=\fontsize{12}{14}\selectfont},
 flow/.style={->,draw=black,line width=.9pt}]
\path[use as bounding box] (0,3.35) rectangle (16.5,21.0);

\coordinate (O) at (.80,17.12);
\draw[->,line width=.7pt] (O)--(.80,20.64) node[above,font=\fontsize{12}{14}\selectfont] {$z$};
\draw[->,line width=.7pt] (O)--(4.50,17.12) node[right,font=\fontsize{12}{14}\selectfont] {$r$};
\fill (O) circle (1.2pt);
\draw[draw=black,fill=white,line width=.85pt]
 (3.30,20.12) ellipse[x radius=.19,y radius=.48];
\draw[draw=black,fill=white,line width=.85pt]
 (2.03,19.18) ellipse[x radius=.10,y radius=.40];
\draw[draw=black,fill=white,line width=.75pt]
 (1.30,18.15) ellipse[x radius=.047,y radius=.24];
\draw[draw=black,fill=white,line width=.65pt]
 (1.00,17.53) ellipse[x radius=.022,y radius=.12];
\foreach \xx/\yy/\rr in {.89/17.31/.012,.845/17.205/.007}
  \fill[black] (\xx,\yy) circle[radius=\rr cm];
\node[anchor=west,font=\fontsize{10}{12}\selectfont,text=black] at (2.23,19.26) {$j$};
\draw[flow] (2.30,18.72) .. controls (2.48,18.20) and (2.94,18.04) .. (3.47,17.96);
\draw[flow] (1.52,17.98) .. controls (1.67,17.62) and (2.03,17.48) .. (2.46,17.42);
\node[anchor=west,text=black,font=\fontsize{10}{12}\selectfont] at (2.64,18.52)
 {Outer flow};
\node[align=center,font=\fontsize{10}{12}\selectfont] at (2.60,16.55)
 {Initial configuration\\[-1pt]{\fontsize{10}{12}\selectfont Upper rings in the $(r,z)$ plane}};

\node[align=left,text width=2.10cm,font=\fontsize{10}{12}\selectfont]
 at (6.12,19.70)
 {\makebox[0pt][r]{\textbullet\hspace{.4em}}$(0,0)$ is the\\stationary\\non-smooth\\point};
\node[align=left,text width=2.10cm,font=\fontsize{10}{12}\selectfont]
 at (6.12,17.85)
 {\makebox[0pt][r]{\textbullet\hspace{.4em}}Inner\\contribution\\$(i>j)$\\negligible};

\draw[panel] (7.30,16.10) rectangle (16.15,20.80);
\node[heading] at (7.70,20.15) {(1)\enspace Angular-flux estimate};
\node[note,anchor=north west,align=left] at (7.70,19.80)
 {first line: transport of relative vorticity.\\
 Inequality: a spatial interpolation-type estimate.};
\node[anchor=center,font=\fontsize{12}{14}\selectfont] at (11.725,17.65)
 {$\displaystyle\begin{aligned}
 |M_{\Phi,j}'(t)|\sim\frac1{m_j}&\left|\int_{\mathbb H}\nabla\Phi(\mathfrak{c}_X)\cdot u_j\dd\nu_j\right|\\[4pt]
 &\lesssim\frac{Q_j^{1/7}E_j(t)^{3/7}}{m_j}\,\partial_r u_j^r(t,0,0).
 \end{aligned}$};

\draw[draw=black,line width=1.6pt,fill=white,rounded corners=2pt]
 (4.70,15.85)--(16.15,15.85)--(16.15,11.35)--(7.30,11.35)
 --(7.30,13.00)--(4.70,13.00)--cycle;
\node[heading] at (5.18,15.23) {\underline{Large strain estimate}};
\node[note,anchor=north west,align=left] at (5.18,14.85)
 {Large self-deformation is allowed. This figure is at a heuristic level.};
\node[anchor=north west,inner xsep=0pt,font=\fontsize{12}{14}\selectfont]
 at (5.18,14.40)
 {$\displaystyle\begin{aligned}
 \largeint{0}{t_j}\partial_r u_j^r(t,0,0)\dd t\hspace{8pt}
 &\stepineq{1}\largeint{0}{t_j}
   \frac{m_j}{Q_j^{1/7}E_j(t)^{3/7}}\,|M_{\Phi,j}'(t)|\dd t\\[4pt]
 &\stepineq{2}\bar\ep_j^{-4/7}\largeint{0}{t_j}|M_{\Phi,j}'(t)|\dd t
   \mkern8mu\stepineq{3}\bar\ep_j^{-4/7}.
 \end{aligned}$};

\node[align=center,font=\fontsize{10}{12}\selectfont] at (2.40,15.05)
 {Smooth cutoff $\Phi$\\[5pt]
 $\Phi{=}0\ (c\leq\tfrac12);\ \Phi{=}1\ (c\geq\tfrac34)$.};
\node[align=center,font=\fontsize{12}{14}\selectfont] at (2.40,13.55)
 {$\displaystyle E_j(t):=\frac12\|u_j(t)\|_{L^2(\R^3)}^2.$};

\draw[panel] (.35,8.42) rectangle (6.95,12.75);
\node[heading] at (.80,12.25) {(2)\enspace Energy estimate};
\node[note] at (.80,11.85) {and the reference aspect condition};
\node[note,anchor=north west,align=left] at (.80,11.55)
 {(Energy growth is controlled by\\
 outer strain.)};
\node[math,anchor=center,font=\fontsize{12}{14}\selectfont] at (3.65,9.88)
 {$\displaystyle\frac{Q_j^{1/7}E_j(t)^{3/7}}{m_j}
   \lesssim\ep_j(t)^{4/7}\leq\bar\ep_j^{4/7}$};
\node[note,anchor=center] at (3.65,8.98) {$\ep_j(t)=R_j(t)/H_j(t),\qquad 0\leq t\leq t_j.$};

\node[anchor=west,font=\fontsize{12}{14}\selectfont] at (7.12,9.785)
 {$\displaystyle\begin{alignedat}{2}
 M_{\Phi,j}(t)&:=\frac1{m_j}\int_{\mathbb H}\Phi(\mathfrak{c}_X)\dd\nu_j,
 &\qquad\dd\nu_j&:=q_j\,r\dd r\dd z.\\[6pt]
 \mathfrak{c}_X&:=\frac{z}{\sqrt{r^2+z^2}},
 &m_j&:=\int_{\mathbb H}\dd\nu_j.
 \end{alignedat}$};

\draw[panel] (.35,3.55) rectangle (16.15,8.15);
\node[heading] at (.80,7.55)
 {(3)\enspace Height dissipation implies \underline{angular redistribution}};
\node[note] at (.80,6.98)
 {Hyperbolic flows connect the height moment $M_{z,j}$ to the cutoff angular moment $M_{\Phi,j}$.};
\node[anchor=center,font=\fontsize{12}{14}\selectfont] at (8.25,5.93)
 {$\displaystyle (M_{z,j}(t)^2)'\lesssim-\lambda_j(t)M_{\Phi,j}(t)^3,
 \quad\lambda_j(t):=\frac{Q_{*,j}R_j(t)^2}{H_j(t)},
 \quad M_{z,j}(t):=\frac1{m_jH_j(t)}\int_{\mathbb H}z\dd\nu_j.$};
\node[note] at (.80,5.08)
 {If classical evolution persists sufficiently long,};
\node[anchor=center,font=\fontsize{12}{14}\selectfont] at (8.25,4.37)
 {$\displaystyle M_{\Phi,j}(0)=1,\quad \text{but} \quad M_{\Phi,j}(t_j)=\frac12
 \quad\text{\fontsize{10}{12}\selectfont for some $t_j$ such that}\quad
 \tau_j(t_j):=\int_0^{t_j}\lambda_j(t)\dd t\lesssim1.$};

\end{tikzpicture}
}
\endgroup
\caption{This figure illustrates how angular redistribution yields a lower
bound for strain. The estimates allow large self-deformation.
The displayed argument is heuristic: the angular-flux formula omits
the contributions of the affine field and the error field, and the
height-dissipation estimate omits the error field. The outer
affine motion is accounted for in the reference scales $R_j(t), H_j(t)$;
the energy $E_j(t)$ still evolves under the full
flow. Further explanations are
given in Section~\ref{sec:constructionidea}. The initial configuration is
schematic and not to scale.}
\label{fig:strainmechanism}
\end{figure}
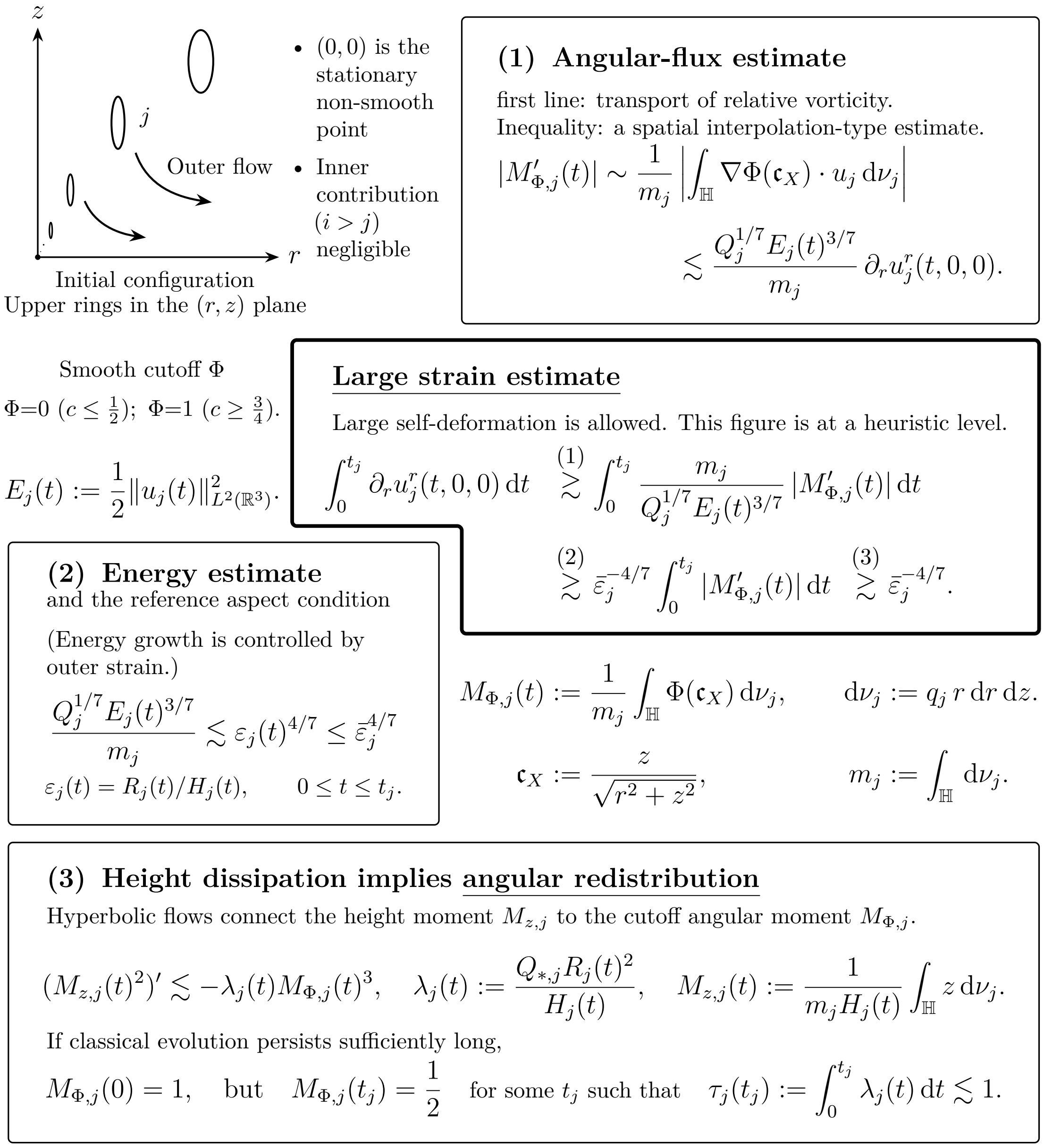

\begin{table}[!t]
\centering
\small
\setlength{\tabcolsep}{5pt}
\begin{tabularx}{\linewidth}{@{}>{\raggedright\arraybackslash}p{.43\linewidth}>{\raggedright\arraybackslash}X@{}}
\hline
Notation & Meaning in the multiscale argument\\
\hline
$u_j,\omega_j,q_j$ & Induced velocity, scalar vorticity, and $-\omega_j/r$\\
$b_j=\partial_r u^r_j(t,0,0)$ & Strain rate of pair $j$ at the origin\\
$B_j=\int_0^t b_j$ & Strain of pair $j$\\
$a_j=\sum_{i\leq j}b_i$ & Accumulated strain rate from the first $j$ pairs\\
$A_j=\sum_{i\leq j}B_i$ & Accumulated strain from the first $j$ pairs\\
$Q_j,Q_{*,j},P_j$ & Relative vorticity amplitude, rescaled relative vorticity amplitude, strain-rate capacity\\
$\mathfrak{e}_j,\ell_j,\mathfrak{E}_j$ & Error field, Lipschitz bound, integrated error\\
$R_{j,0},H_{j,0};\ R_j(t),H_j(t)$ & Initial and evolving reference lengths\\
$\eps_j=R_j/H_j$ & Reference aspect ratio\\
{$\overline A_j,\ \bar\eps_j$} & {{Upper bounds on accumulated strain and reference aspect obtained from bounds on $\tau_j$'s}}\\
{$\lambda_j, \tau_j,\tau_*$} & {Rate of angular-redistribution time, angular-redistribution time, and prescribed threshold}
\\
{$M_{z,j},\ M_{\mathfrak{c},j},\ M_{\Phi,j}$} & {Normalized height, cosine angular, and cutoff angular moments of pair $j$}\\
{$t_0,\ T,\ T_{\mathrm{cmp}}$} & {Classical interval endpoint, time bound $T\geq t_0$, and final comparison time $T_{\mathrm{cmp}}=\tau_*/\lambda_{1,0}$}\\
$T_*$ & Maximal classical lifespan\\
\hline
\end{tabularx}
\end{table}

Earlier multiscale constructions keep self-deformation perturbative, so that its dynamics can be approximated by a system of ODEs. Our construction approaches the H\"older threshold $1/3$ by allowing large self-deformation. The resulting geometry can be complicated, as Figure~\ref{fig:selfdeformationnumerics} illustrates, and we do not attempt to describe it precisely. Instead, we identify geometric quantities that can be investigated directly through the exact PDE instead of the approximating ODE.

Our vortex-ring pairs accumulate only at the origin, which is fixed by symmetry and is the only point where the initial datum is not smooth. The strain $\int \partial_r u^r_j(t,0,0)dt$ generated there by each pair $j$ determines its leading-order influence on smaller pairs. Our main idea is to estimate this strain from below through \emph{angular redistribution}. Figure~\ref{fig:strainmechanism} illustrates how angular-flux, energy, and height-dissipation estimates establish this connection without relying on the ODE approximation. In this subsection~\ref{sec:constructionidea}, we first explain these estimates in a suitable setup where ring pairs are well separated and errors of ring interaction are under control. After that, we describe how the multiscale construction verifies this setup and produces blow-up.

To explain the suitable setup, each vortex-ring pair is transported by
the full Euler flow. Assuming sufficient ring separation, the velocity $u$
near pair $j$ decomposes as
\begin{align}\label{eq:main_dec}
 u=u_{<j}+u_j+u_{>j}
 =
 \partial_r u^r_{<j}(t,0,0)Dx+u_j+\mathfrak{e}_j
\end{align}
where
\[
D=\diag(1,1,-2),\qquad
 \partial_r u^r_{<j}(t,0,0)=\sum_{i<j} \partial_r u^r_i(t,0,0).
\]
Here $u_{<j},u_j,u_{>j}$ are the velocities induced by outer pairs $i<j$, pair $j$ and inner pairs $i>j$. The first-order Taylor approximation of $u_{<j}$ at the origin is the outer affine field $\partial_r u^r_{<j}(t,0,0)Dx$, which stretches pair $j$ radially and compresses it vertically. Near $\supp q_j$, the global error field $\mathfrak{e}_j$ agrees with the local interaction error, the sum of
the outer Taylor remainder, that is $u_{<j}-\partial_r u^r_{<j}(t,0,0)Dx$, and the inner-pair contribution $u_{>j}$. The error field $\mathfrak{e}_j$ will be under control, hence we will mostly suppress it in this subsection~\ref{sec:constructionidea}.

Now, to explain our main idea of using angular redistribution, the first ingredient is the angular-flux estimate in Figure~\ref{fig:strainmechanism}(1). Define $q_j(t,r,z):=-\omega_j(t,r,z)/r$ and then $d\nu_j, m_j$ as in Figure~\ref{fig:strainmechanism}.
Consider the normalized cutoff angular moment $M_{\Phi,j}$ with the cutoff $\Phi(\mathfrak{c}_X)$, which is defined in Figure~\ref{fig:strainmechanism}.
The letter $\mathfrak{c}$ in $\mathfrak{c}_X$ represents the ``cosine'' of an angle.
The cutoff $\Phi$ selects mass near the positive $z$-axis, with
a smooth transition across an intermediate meridional sector. Thus
$M_{\Phi,j}$ measures angular redistribution without tracking the detailed
shape of the vorticity support. Transport of relative vorticity expresses
$M_{\Phi,j}'$ as an angular flux:
\[
M'_{\Phi,j}(t)
=
\frac{1}{m_j}
\int_{\Hh} \nabla \Phi(\mathfrak{c}_X) \cdot u \dd\nu_j.
\] See Remark~\ref{rem:mot_phi} for motivation of the study of $M_{\Phi,j}$ instead of another angular moment.

For a heuristic argument, we suppress, in the decomposition of $u$, the affine field $\partial_r u^r_{<j}(t,0,0)Dx$
and the error field $\mathfrak{e}_j$, as in Figure~\ref{fig:strainmechanism}. Then the spatial
estimate reads 
\[
 |M_{\Phi,j}'(t)|
 \lesssim
 \frac{Q_j^{1/7}E_j(t)^{3/7}}{m_j}\,\partial_r u^r_j(t,0,0),
\]
where $Q_j=\|q_j\|_{L^\infty(\R^3)}$ and $E_j(t)$ is the energy
of the induced velocity $u_j$. The quantities $m_j$ and $Q_j$
are conserved under the full transport. The estimate is derived by using spatial
localization, interpolation, and kernel cancellation;
{see Proposition~\ref{prop:flux} and its development in
Section~\ref{sec:angular}.}

The second ingredient, Figure~\ref{fig:strainmechanism}(2), controls
the energy factor $E_j(t)$ in this inequality. Although the full Euler flow
conserves total energy, the energy $E_j$ associated with an individual
pair can change through its interaction with the other pairs
{(Lemma~\ref{lem:energy})}.
We therefore retain the full flow without the suppression in the energy estimate.
Introduce reference lengths describing the outer affine motion,
{as in Subsection~\ref{subsec:energy-reference}},
\[
 R_j(t)=R_{j,0}e^{A_{j-1}(t)},\qquad
 H_j(t)=H_{j,0}e^{-2A_{j-1}(t)},\qquad
 A_{j-1}(t)=\int_0^t \partial_r u^r_{<j}(s,0,0)\dd s.
\]
These are reference scales, rather than actual dimensions of the evolving support. 
Pure self-interaction conserves $E_j$, but intuitively outer contributions yield energy cascade to smaller scales. Hence, it is not surprising that energy growth can be controlled by outer strain and our energy estimate {(Lemma~\ref{lem:normalizedenergy})} in fact yields
\[
 \frac{Q_j^{1/7}E_j(t)^{3/7}}{m_j}
 \lesssim
 \frac{Q_j^{1/7}\big(Q_j^2 R_j(t)^6 H_j(t)\big)^{3/7}}{Q_j R_j(t)^2 H_j(t)}
 =C\eps_j(t)^{4/7},
 \qquad
 \eps_j(t)=\frac{R_j(t)}{H_j(t)}.
\]
Consequently, under the reference aspect condition
$\eps_j(t)\leq\bar\eps_j$, a substantial change in $M_{\Phi,j}$
implies large strain when $\bar\eps_j$ is small.
See Lemma~\ref{lem:angularaspect} for its rigorous derivation.

The third ingredient is the height-dissipation estimate in
Figure~\ref{fig:strainmechanism}(3). Define the normalized height moment $M_{z,j}$ as in Figure~\ref{fig:strainmechanism}.
Building on the height identities
of Egamberganov--Yao~\cite{EY25},
{we obtain the height-dissipation estimate
shown in Figure~\ref{fig:strainmechanism}(3), with the
error field $\mathfrak{e}_j$ suppressed, using the exact transport and the
Biot--Savart law. See \eqref{eq:heightdissipation} for its rigorous
form.
Using the angular-redistribution time $\tau_j(t)$ defined in
Figure~\ref{fig:strainmechanism}(3), the displayed estimate becomes}
\[
\frac{d}{d\tau_j}\big ( M_{z,j}(t) \big) ^2\lesssim -M_{\Phi,j}(t)^3.
\]
Our initial geometry ensures $M_{\Phi,j}(0)=1$.  If $M_{\Phi,j}(t)$ remains bounded away from zero for a long angular-redistribution time $\tau_j$, then $M_{z,j}(t)^2$ must become negative at some time, which is a contradiction. Therefore,
 $M_{\Phi,j}$
must reach $1/2$ at some physical time $t_j$ satisfying
$\tau_j(t_j)\lesssim1$ if classical evolution
persists sufficiently long.

Combining the three ingredients yields the heuristic strain estimate
displayed in the central box in Figure~\ref{fig:strainmechanism}:
\[
 \begin{aligned}
\int_0^{t_j} \partial_r u^r _j(t,0,0)\dd t
 &\gtrsim
 \bar\eps_j^{-4/7}
 \int_0^{t_j}|M_{\Phi,j}'(t)|\dd t
 \gtrsim\bar\eps_j^{-4/7}.
 \end{aligned}
\]
The last step uses the change from $M_{\Phi,j}(0)=1$ to $M_{\Phi,j}(t_j)=1/2$.
Thus angular redistribution implies large strain at the origin.
The rigorous counterpart, Theorem~\ref{thm:pulse}, controls the contributions of the affine field and the error field
omitted from the figure and uses a 
different time $t_{*,j}$ other than $t_j$. This time $t_{*,j}$ is related to the prescribed threshold $\tau_*$ that appears below. The hypotheses of the theorem concern the reference aspect $\varepsilon_j(t)$,
the integrated error bound, and the duration of classical evolution.

So far we explained how to control our solution from below, but we also need to control our solution from above to preserve its desired structures. One important input for that is the strain-rate capacity bound (Lemmas~\ref{lem:capacity},~\ref{lem:cascadestrain})
\begin{equation}
\label{eq:capacity_bd}
 0\leq \partial_r u^r_j(t,0,0)\leq Q_j^{2/3}m_j^{1/3}\sim P_j,\qquad
 P_j=Q_{*,j}(R_{j,0}^2H_{j,0})^{1/3}.
\end{equation}

The H\"older regularity requirement is met by setting the amplitude rule
$Q_{*,j}=R_{j,0}^{\alpha_j-1}$, which implies
\[
 P_j^3=R_{j,0}^{3\alpha_j-1}H_{j,0}.
\]
If $P_j$ does not grow as $j\to\infty$, then in a natural scale choice, $P_j$ decays exponentially fast, and in that case, the capacity bound leads to global regularity of the solution. 

The threshold $1/3$ is decisive here.  For $\alpha_j<1/3$, sufficiently small radial scales $R_{j,0}$ allow $P_j$ to grow while the heights $H_{j,0}$ shrink. At $\alpha_j=1/3$,
this identity becomes $P_j^3=H_{j,0}$, so it cannot produce growing capacities with shrinking heights.

Finally, the large strain estimate closes the blow-up argument.
While every $\tau_j$ stays below a prescribed threshold $\tau_*$, which is related to angular redistribution,
the capacity estimate provides bounds on $A_j,\varepsilon_j$: $A_j\leq\overline A_j$ and
$\varepsilon_j\leq\bar\varepsilon_j$. We choose the scales so that
$\bar\varepsilon_j\leq\varepsilon_*$ and
$\overline A_j< c_* \bar\varepsilon_j^{-4/7}$.
Suppose  some $\tau_k$ is the first to reach $\tau_*$. Then the large strain
estimate implies $B_k\geq c_* \bar \varepsilon_{k}^{-4/7}>\overline A_k$, contradicting
$B_k\leq A_k\leq\overline A_k$. Therefore, no $\tau_k$ can reach $\tau_*$ while the solution is classical.

However, for the outermost pair, $\tau_1$  must reach $\tau_*$ at a finite time. For this pair, $a_0=0$, so its reference scales are constant and $\tau_1(t)=\lambda_{1,0}t.$ Classical continuation to the finite comparison time $T_{\mathrm{cmp}}=\tau_*/\lambda_{1,0}$ would therefore imply $\tau_1(T_{\mathrm{cmp}})=\tau_*$. This yields finite-time blow-up. 

Here, we simplified several aspects of our actual proofs to explain their essential part.
}

\subsection{Notation}

All unlabelled spatial norms are over $\R^3$. For $0<\alpha<1$,
\[
 \|f\|_{C^\alpha} = \|f\|_{L^\infty(\R^3)}+ \sup_{x\ne y}\frac{|f(x)-f(y)|}{|x-y|^\alpha}.
\]
The inhomogeneous $C^{1,\alpha}$ norm includes $\|f\|_{L^\infty(\R^3)}$,
$\|\nabla f\|_{L^\infty(\R^3)}$, and $[\nabla f]_{C^\alpha}$.
Matrix norms are operator norms, and $B_R$ denotes the ball of radius
$R$ centered at the origin. We write $A \lesssim B$ when
$A \leq CB$ for a constant $C>0$ independent of time and the scales under
consideration, and $A \sim B$ when both inequalities hold. Constants
may depend on the fixed profile $f_0$. {In the estimates
with unspecified multiscale parameters, constants may also depend on
the fixed height parameter $\kappa$; this dependence is suppressed after
its numerical value is chosen. The quantity $A_{\mathrm{Lip}}(t_0)$ is
different: it depends on the given classical solution and interval,
and the initial parameter choices are independent of it.} The angular and radial cutoffs $\Phi, \chi$
in \eqref{eq:angularcutoff} and \eqref{eq:radialcutoffs} are fixed
universally; constants depending on these cutoffs are suppressed. Additional dependence
is indicated by a subscript. Parameter choices follow
Setup~\ref{setup:parameters} (an initial family). 

For a $C^1$ axisymmetric no-swirl velocity, radial differentiation
on the axis means the continuous extension
\[
 \partial_r u^r(t,0,z)=\lim_{r\to0^+}\frac{u^r(t,r,z)}r.
\]
In particular, the strain rate at the origin is written
$\partial_r u^r(t,0,0)$. Its time integral is called strain; ``accumulated'' denotes
summation over vortex-ring pairs, for both strain rates and strains.

{For one vortex-ring pair with induced velocity $u$, we write
$b(t)=\partial_r u^r(t,0,0)$ for the strain rate it generates at the
origin and $B(t)=\int_0^t b(s)\dd s$ for its strain.
The coefficient $a(t)$ is the strain rate of the
affine field $a(t)Dx$, and $A(t)=\int_0^t a(s)\dd s$ is the
corresponding affine strain. We denote $A_{\rm tot}=A+B$; the contribution
of the error field $\mathfrak{e}$ is kept separate.}

{In the multiscale application to pair $j$, these quantities are
$(a,b)=(a_{j-1},b_j)$ and $(A,B,A_{\rm tot})=(A_{j-1},B_j,A_j)$.
Thus $b_j,B_j$ refer to the strain rate and strain generated at the
origin by pair $j$, while $a_j=\sum_{i=1}^j b_i$ and
$A_j=\sum_{i=1}^j B_i$ are the accumulated strain rate and accumulated
strain from the first $j$ pairs. In particular, $A_j'=a_j$, with
$a_0=A_0=0$. Relative to pair $j$, the outer quantities are
$a_{j-1}$ and $A_{j-1}$; relative to pair $j+1$, they are $a_j$ and
$A_j$. The notation $b_j,B_j$ follows \cite{BC26}, but here the
quantities refer to the exact transported vorticities of the
vortex-ring pairs.}

\begingroup

\section{Transported relative vorticity and kinetic energy}\label{sec:ringpairs}

In Sections~\ref{sec:ringpairs}--\ref{sec:strainbounds}, we prove
the single-pair estimates outlined in Subsection~\ref{sec:constructionidea}
and Figure~\ref{fig:strainmechanism}.
The evolving pair is transported by $U=u+aDx+\mathfrak{e}$, where $u$ is its
own induced velocity. The \emph{prescribed} external field is the sum of the affine field $aDx$ and the error field $\mathfrak{e}$.
We later apply this framework to each pair of the multiscale solution.

\subsection{The transported equation}

Let $\Hh = \{(r,z):r>0,\ z>0\}$ be the upper meridional quadrant.
The matrix $D = \diag(1,1,-2)$ describes unit radial expansion with
the axial compression required by incompressibility. We distinguish
a vortex-ring pair at a fixed time from a vortex-ring pair satisfying a transport
equation. The setups below record the shared assumptions. Later each result will also
state additional hypotheses if needed.

\begin{setup}[A fixed-time pair]\label{setup:ringpair}
Let $q\in C^{1,\alpha}(\R^3)$, for some $0<\alpha<1$, be
axisymmetric and compactly supported away from the axis and the plane
$z=0$. Assume that $q$ is odd in $z$ and nonnegative on $\Hh$.
Define $\vectomega(x)=q(r,z)\,(x_2,-x_1,0)=\omega\mathbf e_\theta$.
The self-induced velocity $u$ is axisymmetric without swirl and
divergence-free, with $u^r$ even and $u^z$ odd in $z$. It is defined by
\begin{equation}\label{eq:bs}
 u(x) = \BS[\vectomega](x) = \frac1{4\pi}\int_{\R^3} \frac{\vectomega(y)\times(x-y)}{|x-y|^3}\dd y.
\end{equation}
Thus $\vectomega$ is even,
$\vectomega(-x)=\vectomega(x)$, and $u(-x)=-u(x)$. We use
\[
 \begin{gathered}
 \dd\nu=rq\dd r\dd z\quad\text{on }\Hh,\qquad
 m=\int_{\Hh}\dnu,\qquad Q=\|q\|_{L^\infty(\R^3)},\\
 E=\tfrac12\|u\|_{L^2(\R^3)}^2,\qquad b=\partial_r u^r(0,0).
 \end{gathered}
\]
This setup imposes no evolution equation.
\end{setup}

Throughout, continuity on a time interval is understood relative to
that interval; at included endpoints, continuity and time derivatives
are one-sided. Compact time subintervals are denoted by
$I_0=[s_0,s_1]\subseteq I$ and may include endpoints of $I$.
Local integrability and local uniform convergence in time are understood
on all such compact subintervals.

\begin{setup}[An evolving pair]\label{setup:evolution}
Let $I=[0,T)$ with $0<T\leq\infty$, or $I=[0,T]$ with
$0<T<\infty$. On $I$, let the vortex-ring pair satisfy
Setup~\ref{setup:ringpair} (a fixed-time pair) at each time, and let $q$ be a classical
solution of the transport equation
\begin{equation}\label{eq:tageq}
 q_t+(u+a(t)Dx+\mathfrak{e})\cdot\nabla q = 0, \qquad u = \BS[\vectomega], \qquad a \geq 0.
\end{equation}
The support remains separated from the axis and the symmetry plane
on each compact subinterval of $I$. Here $a\in C(I)$ is a prescribed
nonnegative strain rate of the affine field $aDx$, and
$\mathfrak{e}\in C(I;C^1_{\mathrm{loc}}(\R^3))$ is an error field with the
same axisymmetry, absence of swirl, and parity as $u$, satisfying
\begin{equation}\label{eq:errorhyp}
 \Div \mathfrak{e} = 0, \qquad \|\nabla \mathfrak{e}(t)\|_{L^\infty(\R^3)} \leq \ell(t), \qquad \mathfrak{E}(t) = \int_0^t\ell(s)\dd s \leq 1.
\end{equation}
Here $0\leq\ell\in L^1_{\mathrm{loc}}(I)$, the gradient is taken in Cartesian coordinates,
and the bounds hold on $I$. We denote the full transporting velocity by $U=u+aDx+\mathfrak{e}$ whereas $u$ denotes only
the self-induced velocity of this vortex-ring pair.
\end{setup}

The notation $\ell$ stands for the letter L in the Lipschitz norm. The other notation $\mathfrak{E}$ is a Fraktur style of the letter $E$ as it is related to the integral of $\mathfrak{e}$.

Axis regularity and plane parity imply
\begin{equation}\label{eq:errorcoordinates}
 |\mathfrak{e}^r(r,z)| \leq \ell r, \qquad |\mathfrak{e}^z(r,z)| \leq \ell|z|.
\end{equation}

We explain why the error field $\mathfrak{e}$ is defined globally and how Setup~\ref{setup:evolution} (an evolving pair) applies to a multiscale solution. In the multiscale decomposition \eqref{eq:main_dec}, the interaction
error
\[
 \bigl(u_{<j}(t,r,z)-\partial_r u^r_{<j}(t,0,0)Dx\bigr)+u_{>j}(t,r,z)
\]
is controlled only near $\supp q_j$. However, one can modify the velocity away
from this support because it does not change the transport of $q_j$.

We will construct a global error field
$\mathfrak{e}_j$ (Lemma~\ref{lem:combinederror}) agreeing with this interaction error near $\supp q_j$, preserving
the required symmetries and incompressibility, and satisfying
a controlled global gradient bound. Thus, with
$a_{j-1}=\partial_r u^r_{<j}(t,0,0)$,
\[
 \partial_tq_j+(u_j+a_{j-1}Dx+\mathfrak{e}_j)\cdot\nabla q_j=0
\]
holds globally. This motivates the global formulation of $\mathfrak{e}$ in
Setup~\ref{setup:evolution} (an evolving pair) and explains how we apply this setup to a multiscale solution.

Coming back to Setup~\ref{setup:evolution} (an evolving pair),  $\dd\nu = -\omega\dd r\dd z=q \, r\dd r\dd z$ is transported, and thus
\begin{equation}\label{eq:mass}
 Q = \|q_0\|_{L^\infty(\R^3)} = \|q(t)\|_{L^\infty(\R^3)}, \qquad m = \int_{\Hh}\dnu
\end{equation}
are conserved. 

\begin{lemma}[Strain rate at the origin]\label{lem:originstrain}
Under Setup~\ref{setup:ringpair} (a fixed-time pair),
\begin{equation}\label{eq:source}
b = \partial_r u^r(0,0) = \int_{\Hh}k(r,z)\dnu,
 \qquad
k(r,z) = \frac32\frac{r^2z}{(r^2+z^2)^{5/2}} \geq 0.
\end{equation}
Moreover, $b\geq0$ and $\nabla u(0)=bD$.
\end{lemma}

\begin{proof}
Differentiating \eqref{eq:bs} and integrating in the azimuthal
variable, we obtain \eqref{eq:source}. Due to axisymmetry and absence of swirl,
$\nabla u(0)$ is diagonal with equal first two entries, and its trace
vanishes. Hence $\nabla u(0)=bD$.
\end{proof}

\subsection{Energy and reference scales}\label{subsec:energy-reference}

In this subsection, we will establish the energy estimate explained in Figure~\ref{fig:strainmechanism}(2).
The energy of one vortex-ring pair need not be conserved under
the full transporting velocity $U$. The following identity measures
its change through the affine field $aDx$ and the error field $\mathfrak{e}$. The affine field permits energy growth at rate $4a$.

\begin{lemma}[Energy of a transported vortex-ring pair]\label{lem:energy}
Under Setup~\ref{setup:evolution} (an evolving pair), let
$E(t)=\|u(t)\|_{L^2(\R^3)}^2/2$. Then
\begin{equation}\label{eq:energyidentity}
 E' = -\int_{\R^3}u_iu_j\partial_j(aDx+\mathfrak{e})_i\dd x,
\end{equation}
and
\begin{equation}\label{eq:energybound}
 E' \leq 4aE-3a\|u^r\|_{L^2(\R^3)}^2+2\ell E.
\end{equation}
\end{lemma}

\begin{proof}
Let $U=u+aDx+\mathfrak{e}$. 

Fix a compact subinterval $I_0\subseteq I$. By the transport formula,
$\vectomega$ is continuous on $I_0\times\R^3$, and there is a
compact set in $\mathbb{R}^3$ that contains 
$\supp\vectomega(t,\cdot)$ for every $t\in I_0$.

The Biot--Savart estimate
$\|u\|_{L^2(\R^3)}\lesssim\|\vectomega\|_{L^{6/5}(\R^3)}$
implies $u\in C(I_0;L^2(\R^3))$, and
$U\times\vectomega\in C(I_0;L^2(\R^3))$.
Equation \eqref{eq:tageq} implies
$\partial_t\vectomega=\curl(U\times\vectomega)$, hence
$\partial_tu=\mathbb P(U\times\vectomega)$ distributionally,
where $\mathbb P$ is the $L^2$-orthogonal Leray projection.
Thus $u\in C^1(I_0;L^2(\R^3))$. Testing the equation of $u$ with $u$ now justifies
\[
 E' = \int_{\R^3}u\cdot[(aDx+\mathfrak{e})\times\vectomega]\dd x.
\]
For a divergence-free field $w$, integrating by parts, we obtain
\[
 \int_{\R^3}u\cdot(w\times\curl u)\dd x = -\int_{\R^3}u_iu_j\partial_jw_i\dd x.
\]
Spatial cutoffs justify the identity for a field $w$ of at most linear growth. Hence \eqref{eq:energyidentity} follows.

Now to show \eqref{eq:energybound}, note that the affine contribution in \eqref{eq:energyidentity} by $aDx$ is
\[
 -a\int_{\R^3}(|u^r|^2-2|u^z|^2)\dd x = 4aE-3a\|u^r\|_{L^2(\R^3)}^2,
\]
and the contribution of the error field is at most $2\ell E$, which completes the proof.
\end{proof}

\begin{setup}[An evolving, scaled pair]\label{setup:profile}
In Setup~\ref{setup:evolution} (an evolving pair), let $(\rho,\zeta)$ denote dimensionless
radial and vertical
variables. Choose $f_0\in C_c^\infty(\Hh)$, $f_0 \geq 0$, normalized by
$\int_{\Hh}\rho f_0(\rho,\zeta)\dd\rho\dd\zeta = 1$, and supported in
\begin{equation}\label{eq:profilebox}
 [\rho_-,\rho_+]\times[\zeta_-,\zeta_+], \qquad \rho_->0, \qquad \zeta_->0.
\end{equation}
Let $R_0,H_0>0$ be the initial radial and vertical scales, and
choose $Q_*>0$. We start from
$q_0(r,z)=Q_*f_0(r/R_0,z/H_0)$ for $z>0$, with odd reflection
for $z<0$. The reference scales and normalized quantities are defined
in \eqref{eq:referencescales}--\eqref{eq:normalizedenergy} below.
No restriction on $R_0/H_0$ is imposed unless stated explicitly.
\end{setup}

Let
$A(t) = \int_0^ta(s)\dd s$ be the affine strain.
The affine field $aDx$ alone expands the radial scale by $e^A$
and compresses the vertical scale by $e^{-2A}$. Accordingly, let
\begin{equation}\label{eq:referencescales}
 \begin{gathered}
 R(t) = R_0e^{A(t)}, \qquad H(t) = H_0e^{-2A(t)},\\
 \eps=\frac RH, \qquad \lambda = \frac{Q_*R^2}{H}, \qquad \tau(t) = \int_0^t\lambda(s)\dd s.
 \end{gathered}
\end{equation}

In these coordinates, define the evolving profile by
\[
 f(t,\rho,\zeta):=\frac{1}{Q_*}q\bigl(t,R(t)\rho,H(t)\zeta\bigr),
 \qquad (\rho,\zeta)\in\Hh.
\]
Thus $f(0,\rho,\zeta)=f_0(\rho,\zeta)$.

\begin{remark}[Comparison with \cite{CMZZ25} on ansatz and profiles]\label{rem:CMZZansatz}
Using the evolving profile $f$, one may write the relative vorticity as follows:
\[
q(t,r,z)
= Q_* f\left( t, \frac{r}{R(t)}, \frac{z}{H(t)} \right).
\]
Compared to \cite{CMZZ25}, there are two main differences regarding the profile $f$. 

First, in \cite{CMZZ25}, self-deformation is perturbative, leading to closeness of the evolving profile to the prescribed reference profile (in fact, uniformly in scale index $j$). This ansatz plays a central role in approximating the PDE by a system of ODEs. However, in our setup, as self-deformation is not required to be perturbative, our evolving profile $f$ may deviate significantly from our reference profile $f_0$. Hence, such an ansatz breaks down in our setup, and cannot be used in our proofs. The definition of $f$ is introduced only for this comparison. 

Second, the ansatz used in \cite{CMZZ25} is in fact formulated for the vorticity $\omega$, not for the relative vorticity $q$. We shifted from the vorticity to the relative vorticity to utilize its transport structure, which is more suitable to allow large self-deformation.
\end{remark}

Here $\eps=R/H$ is the reference aspect ratio, the reciprocal of
the ratio used in \cite{BC26}; $\tau$ is the angular-redistribution
time, and $\lambda=\tau'$ is its rate with respect to physical time.
 The conserved quantities
$m$ and $Q$ satisfy
\begin{equation}\label{eq:normalizedmass}
 m = Q_*R^2H, \qquad F = \|f_0\|_{L^\infty(\Hh)}, \qquad Q = Q_*F.
\end{equation}

To motivate energy normalization, note that first of all, the factor $R^6H$ grows at the affine energy rate $4a$, that is $(R(t)^6 H(t))'=4a(t) R(t)^6 H(t)$. 
In addition, since $|\omega_0|=r|q_0|$ and $r\sim R_0$ on the support,
the vorticity has characteristic size $Q_*R_0$. Taking $R_0$ as
the characteristic length suggests a velocity of size $Q_*R_0^2$.
The vorticity support has volume comparable to $R_0^2H_0$, so
squaring this velocity scale and multiplying by that volume
suggests the energy scale
\[
 (Q_*R_0^2)^2 R_0^2H_0 = Q_*^2R_0^6H_0.
\]
We therefore define the normalized energy by
\begin{equation}\label{eq:normalizedenergy}
 \cE(t) = \frac{E(t)}{2\pi Q_*^2R(t)^6H(t)}.
\end{equation}

We define the radial primitive of $f_0$:
\[
 G_f(\rho,\zeta) = \int_\rho^\infty s f_0(s,\zeta)\dd s, \qquad E_f = \int_{\Hh}\rho G_f^2\dd\rho\dd\zeta.
\]
The constant $E_f$ depends only on the chosen profile. For simplicity, we use the subscript $f$ for $G_f, E_f$, but more precisely, both of them depend on the initial profile $f_0$, not on the evolving profile $f$.

The normalization controls energy without comparing the evolving
profile $f(t,\rho,\zeta)$ to $f_0(\rho,\zeta)$. Intuition behind this lemma is explained in Section~\ref{sec:constructionidea}.

\begin{lemma}[Normalized energy]\label{lem:normalizedenergy}
Under Setup~\ref{setup:profile} (an evolving, scaled pair), for every positive initial aspect ratio $\varepsilon(0)$,
\begin{equation}\label{eq:Estar}
 \cE(t) \leq E_fe^{2\mathfrak{E}(t)} \leq E_fe^2 \eqqcolon E_*.
\end{equation}
\end{lemma}

\begin{proof}
Since $(R^6H)' = 4aR^6H$, \eqref{eq:energybound} implies
\[
\cE' \leq 2\ell\cE, \quad \text{and thus }\cE (t) \leq \cE(0) e^{2\mathfrak{E}(t)}.
\]
It remains to estimate $E(0)$.

For $z\in\R$, let $G_0(r,z)=\int_r^\infty s q_0(s,z)\dd s$. Then
\[
 \curl(-G_0\mathbf e_3)
 =(\partial_rG_0)\mathbf e_\theta
 =-rq_0\mathbf e_\theta=\vectomega_0.
\]
By uniqueness of the divergence-free $L^2$ velocity with the same curl
$\vectomega_0$, it holds $u(0)=\mathbb P(-G_0\mathbf e_3)$. Hence
\[
 \frac{E(0)}{2\pi}
 =\frac{\|\mathbb P(-G_0\mathbf e_3)\|_2^2}{4\pi}
 \leq\frac{\|G_0\|_2^2}{4\pi}
 =\int_{\Hh}rG_0^2\dd r\dd z.
\]
Since $G_0(r,z) = Q_*R_0^2G_f(r/R_0,z/H_0)$,
we obtain $\cE(0) \leq E_f$.
\end{proof}

\section{Angular localization and angular flux}\label{sec:angular}

We prove the spatial angular-flux estimate explained in
Subsection~\ref{sec:constructionidea} and Figure~\ref{fig:strainmechanism}(1),
under Setup~\ref{setup:ringpair} (a fixed-time pair). The quantity to be estimated is
\begin{equation}
\label{eq:ang_flux}
 \int_{\Hh}\nabla\Phi(\mathfrak{c}_X)\cdot u(X)\dnu(X),
 \qquad \mathfrak{c}_X=z/|X|,
\end{equation}
where $\Phi$ is the smooth cutoff specified below.
We establish the angular-flux estimate through angular localization and kernel symmetrization.

\subsection{Mass in an intermediate sector}

Bounded relative vorticity controls mass near the origin but becomes
ineffective at large radii. Energy supplies a complementary bound
there. 

An intermediate sector is a region
$\mathcal S = \{X\in\Hh:c_- \leq \mathfrak{c}_X \leq c_+\}$, with fixed
$0<c_-<c_+<1$.

\begin{lemma}[Angular localization]\label{lem:polarization}
Under Setup~\ref{setup:ringpair} (a fixed-time pair), 
for every intermediate sector $\mathcal S$ and every $\varrho>0$,
\begin{equation}\label{eq:exteriorangular}
 \nu(\mathcal S\cap\{|X| \geq \varrho\}) \lesssim_{\mathcal S}E^{1/2}\varrho^{-1/2}.
\end{equation}
Consequently,
\begin{equation}\label{eq:bandmass}
 \nu(\mathcal S) \lesssim_{\mathcal S} Q^{1/7}E^{3/7}.
\end{equation}
\end{lemma}

\begin{proof}
The case $q=0$ is immediate, so assume $q\ne0$.

Choose a smooth radial cutoff $0\leq\chi_\varrho\leq1$, zero below
$\varrho/2$ and one above $\varrho$, with
$|\chi_\varrho'|\lesssim\varrho^{-1}$. Let
\[
 \varphi_\varrho(X)=\chi_\varrho(|X|)\frac{r^2z}{|X|^3}
 =\chi_\varrho(|X|)\mathfrak{c}_X(1-\mathfrak{c}_X^2).
\]
Writing $s=|X|$ and differentiating, we get
$|\nabla \varphi_\varrho|\lesssim r/s^2$ for $s\geq\varrho/2$,
since the cutoff varies only where $s\sim\varrho$. Thus
\[
 \int_{\Hh}\frac{|\nabla \varphi_\varrho|^2}{r}\dd r\dd z
 \lesssim\int_{\Hh\cap\{s\geq\varrho/2\}}\frac{r}{s^4}\dd r\dd z
 \lesssim\varrho^{-1}.
\]
Using $\dnu=(\partial_ru^z-\partial_zu^r)\dd r\dd z$
and integrating by parts, we obtain
\[
 \begin{aligned}
 \int_{\Hh}\varphi_\varrho\dnu
 &=\int_{\Hh}
 (u^r\partial_z \varphi_\varrho-u^z\partial_r \varphi_\varrho)\dd r\dd z\\
 &\leq\left(\frac{E}{2\pi}\right)^{1/2}
 \left(\int_{\Hh}\frac{|\nabla \varphi_\varrho|^2}{r}\dd r\dd z\right)^{1/2}
 \lesssim E^{1/2}\varrho^{-1/2}.
 \end{aligned}
\]
Here $E/(2\pi)=\int_{\Hh}r|u|^2\dd r\dd z$ by symmetry.
The weight vanishes on both boundaries of the quadrant.
The boundary term on $|X|=L$ tends to zero as $L\to\infty$,
since $\varphi_\varrho$ is bounded and $u=O(L^{-2})$ there. Since $\chi_\varrho=1$ on $|X|\geq\varrho$ and
$\mathfrak{c}_X(1-\mathfrak{c}_X^2)$ is bounded below by a positive constant on
$\mathcal S$, we obtain \eqref{eq:exteriorangular}.

The bound $q\leq Q$ implies
$\nu(\{|X|<\varrho\})\leq Q\int_{\{|X|<\varrho\}\cap\Hh}r\dd r\dd z
\lesssim Q\varrho^3$. Splitting the sector into its parts inside
and outside this radius, we get
\[
 \nu(\mathcal S)\lesssim_{\mathcal S}
 Q\varrho^3+E^{1/2}\varrho^{-1/2}.
\]
The two terms balance at $\varrho_*=(E/Q^2)^{1/7}$, where both
are $Q^{1/7}E^{3/7}$. This yields the bound \eqref{eq:bandmass} for the mass of the entire sector.
\end{proof}

\begin{remark}[Comparison with Choi--Jeong]
\label{rem:cjlocalization}
The localization bounds of Lemma~\ref{lem:polarization} also
follow from Choi--Jeong's positive energy-kernel estimate
\cite[Lemma~3.3, Eq.~(3.1)]{CJ25}. Our proof tests the meridional vorticity identity against the angular weight $\varphi_\varrho$.

To control angular flux, we additionally use kernel
symmetrization and the logarithmic mass estimate in
Lemma~\ref{lem:logmass} below to handle the remaining singularity;
see Proposition~\ref{prop:flux}.
\end{remark}

\begin{lemma}[Localized logarithmic mass]\label{lem:logmass}
Under Setup~\ref{setup:ringpair} (a fixed-time pair), fix an intermediate sector $\mathcal S$. For every
$0<\eta\leq1/2$ and $X\in\mathcal S$, let $\varrho=|X|$. Then
\begin{equation}\label{eq:logmass}
 \int_{\mathcal S\cap\{|X-Y| \leq \eta \varrho\}} \left(1+\log\frac \varrho{|X-Y|}\right)\dnu(Y) \lesssim_{\mathcal S} Q^{1/7}E^{3/7}.
\end{equation}
The implicit constant depends only on $\mathcal S$, independently of $\eta$.
\end{lemma}

\begin{proof}
Related logarithmic estimates against product vorticity measures
appear in \cite[Lemma~3.3]{GMT26} and
\cite[Lemma~2.1, Eq.~(2.27)]{CFQ26}. Here we combine the exterior
mass bound \eqref{eq:exteriorangular} with $q\leq Q$ to estimate
the integral uniformly in the point $X$.
Fix $X\in\mathcal S$, write $\varrho=|X|$, and take
\[
 \mathcal A_\varrho=\mathcal S\cap\{\varrho/2\leq |Y|\leq2\varrho\},
 \qquad N_\varrho=\nu(\mathcal A_\varrho).
\]
Since $\eta\leq1/2$, the integration region in \eqref{eq:logmass}
is contained in $\mathcal A_\varrho$.
Writing $Y=(s,w)$, the density $sq(s,w)$ of $\dnu(Y)$
on $\mathcal A_\varrho$ is at most $L=CQ\varrho$, and
\eqref{eq:exteriorangular} applied at radius $\varrho/2$ implies
$N_\varrho\lesssim_{\mathcal S}E^{1/2}\varrho^{-1/2}$.
For $0<\delta<2$ and $a>0$, split the integral according to
$|X-Y|\leq a$ and $|X-Y|>a$. The density bound on the first region
and the mass bound on the second imply
\[
 \begin{aligned}
 \int_{\mathcal A_\varrho\cap\{|X-Y|\leq a\}}
 |X-Y|^{-\delta}\dnu(Y)
 &\leq 2\pi L\int_0^a s^{1-\delta}\dd s
 \lesssim_\delta La^{2-\delta},\\
 \int_{\mathcal A_\varrho\cap\{|X-Y|>a\}}
 |X-Y|^{-\delta}\dnu(Y)
 &\leq a^{-\delta}N_\varrho.
 \end{aligned}
\]

If $N_\varrho=0$, the desired bound is immediate. Otherwise,
choosing $a=(N_\varrho/L)^{1/2}$ and adding the two estimates and multiplying by $\varrho^\delta$ yield
\[
 \int_{\mathcal A_\varrho}\left(\frac{\varrho}{|X-Y|}\right)^\delta\dnu(Y)
 \lesssim_\delta \varrho^\delta(La^{2-\delta}+a^{-\delta}N_\varrho)
 \lesssim_\delta L^{\delta/2}N_\varrho^{1-\delta/2}\varrho^\delta
 \lesssim_\delta Q^{\frac{\delta}{2}} E^{\frac{1}{2}\left( 1-\frac{\delta}{2}\right)}\varrho^{\frac{7\delta-2}{4}}.
\]
Taking $\delta=2/7$, we obtain
\[
 \int_{\mathcal A_\varrho}
 \left(\frac{\varrho}{|X-Y|}\right)^{2/7}\dnu(Y)
 \lesssim_{\mathcal S} Q^{1/7}E^{3/7}.
\]
Since $1+\log(\varrho/d)\lesssim(\varrho/d)^{2/7}$ for
$0<d\leq\eta\varrho\leq\varrho/2$, with an absolute implicit constant,
\eqref{eq:logmass} holds uniformly in $\eta$.
\end{proof}

\subsection{The stream kernel}

We study a stream function kernel whose symmetry yields cancellation in the angular-flux estimate. Its factorization also
helps estimate interactions between separated vortex-ring pairs
in Section~\ref{sec:interactions}.

We normalize the meridional stream function on the axis defining it as
\[
 \psi(r,z)=\int_0^r s\,u^z(s,z)\dd s.
\]
Then $\partial_r\psi=ru^z$. Integrating
$s\partial_z u^z=-\partial_s(su^r)$ from the axis, we get
$\partial_z\psi=-ru^r$. Thus $u=r^{-1}\nabla^\perp\psi$ for
$\nabla^\perp=(-\partial_z,\partial_r)$. The axisymmetric stream function formula
\cite[Section~2.2, Eq.~(2.3)]{CJ25}
(see also \cite[Eqs.~(2.8)--(2.13)]{FS15}), with our convention
$\omega=-rq$, yields the representation below after using the odd symmetry and writing
$\dnu(s,w)=sq(s,w)\dd s\dd w$:
\begin{equation}
\label{eq:reppsi}
 \psi(X)=\int_{\Hh}G(X,Y)\dnu(Y),
\end{equation}
where, for $X=(r,z)$ and $Y=(s,w)$,
\begin{equation}\label{eq:streamkernel}
 \begin{gathered}
 G(X,Y) = -\frac{rs}{4\pi}\int_0^{2\pi} \cos\theta\,(D_-^{-1}-D_+^{-1})\dd\theta,\\
 D_\pm^2 = (r-s)^2+(z\pm w)^2+2rs(1-\cos\theta).
 \end{gathered}
\end{equation}
The kernel satisfies $G(X,Y)=G(Y,X)$. In the angular-flux estimate, we use
this symmetry to interchange the integration variables $X$ and $Y$.

\begin{lemma}[Kernel factorization]\label{lem:factorization}
For the kernel $G$ defined in \eqref{eq:streamkernel} and distinct
$X=(r,z),Y=(s,w)\in\Hh$,
\begin{equation}\label{eq:kernelfactor}
 \begin{split}
 G(X,Y)& = -\frac3{2\pi}r^2s^2zw\,\mathcal H(X,Y),\\
 \mathcal H(X,Y)& = \int_0^{2\pi}\int_0^1 \frac{\sin^2\theta\dd t\dd\theta} {[(1-t)D_-^2+tD_+^2]^{5/2}}>0 .
 \end{split}
\end{equation}
In addition, under Setup~\ref{setup:ringpair} (a fixed-time pair), the strain rate at the origin $b=\partial_r u^r(0,0)$ satisfies
\[
b =-\lim_{\substack{X\to 0\\ X\in\Hh}}\frac{\psi(X)}{r^2z}
 =\frac32\int_{\Hh}\frac{s^2w}{|Y|^5}\dnu(Y).
\]
\end{lemma}

\begin{proof}
Integration by parts in $\theta$ yields
\[
 \int_0^{2\pi}\cos\theta\,D_\pm^{-1}\dd\theta = rs\int_0^{2\pi}\sin^2\theta\,D_\pm^{-3}\dd\theta.
\]
Since $D_+^2-D_-^2=4zw$, differentiating in $t$, we get
\[
 \frac{\dd}{\dd t}\bigl[(1-t)D_-^2+tD_+^2\bigr]^{-3/2}
 =-6zw\bigl[(1-t)D_-^2+tD_+^2\bigr]^{-5/2}.
\]
Integrating from $t=0$ to $t=1$ and reversing the sign, we obtain
\[
 D_-^{-3}-D_+^{-3} = 6zw\int_0^1[(1-t)D_-^2+tD_+^2]^{-5/2}\dd t.
\]
Substituting this identity into \eqref{eq:streamkernel} yields
factorization \eqref{eq:kernelfactor}. At $X=0$, both $D_-^2$ and $D_+^2$
equal $|Y|^2$, so $\mathcal H(0,Y)=\pi|Y|^{-5}$. Since $\nu$ is compactly
supported away from the origin, $D_\pm\geq |Y|-|X|\geq |Y|/2$ for
$Y\in\supp\nu$ and $X$ sufficiently close to $0$. This bound allows us
to pass to the limit by dominated convergence after substituting
\eqref{eq:kernelfactor} into \eqref{eq:reppsi} and canceling $r^2z$.
We obtain
\[
 -\lim_{\substack{X\to0\\ X\in\Hh}}\frac{\psi(X)}{r^2z}
 =\frac{3}{2\pi}\int_{\Hh}s^2w\,\mathcal H(0,Y)\dnu(Y)
 =\frac32\int_{\Hh}\frac{s^2w}{|Y|^5}\dnu(Y)
 =b,
\]
where the last equality follows from the strain-rate formula \eqref{eq:source}.
\end{proof}

\subsection{An angular flux estimate}
To estimate the angular flux \eqref{eq:ang_flux}, we will use the representation \eqref{eq:reppsi} to re-write it as a double integral with respect to two meridional points $X,Y$.  
The kernel will have singularity at $X=Y$. For points $X,Y$ that are separated, the
factors in \eqref{eq:kernelfactor} permit comparison with the strain-rate
kernel. For points $X, Y$ that are close to each other, interchanging them in the double integral through kernel symmetrization
produces a cancellation, and Lemma~\ref{lem:logmass} controls the
remaining logarithmic singularity. 

We fix a smooth nondecreasing function
$\Phi:[0,1]\to[0,1]$ such that
\begin{equation}\label{eq:angularcutoff}
 \Phi(c) = 0 \quad \text{for} \quad c \leq 1/2, \qquad \Phi(c) = 1 \quad \text{for} \quad c \geq 3/4.
\end{equation}

\begin{proposition}[Angular flux]\label{prop:flux}
Under Setup~\ref{setup:ringpair} (a fixed-time pair), for $\Phi$ as in \eqref{eq:angularcutoff},
\begin{equation}\label{eq:flux}
 \left|\int_{\Hh}\nabla\Phi(\mathfrak{c}_X)\cdot u(X)\dnu(X)\right| \lesssim_\Phi Q^{1/7}E^{3/7}b.
\end{equation}
The constant is independent of the vortex-ring pair's support scales and of all
derivatives of $q$. And here $b:=\partial_ru^r(0,0)$.
\end{proposition}

\begin{proof}
Let $\Gamma(X,Y) = r^{-1}\nabla_X^\perp G(X,Y)$. Then using the representation \eqref{eq:reppsi} of $\psi$ together with $u=r^{-1}\nabla^{\perp}\psi$, one can re-write the given integral as follows
\[
\int_{\Hh}
\nabla \Phi(\mathfrak{c}_X)\cdot u(X) \dnu(X) = \int_{\Hh}\int_{\Hh}
\nabla \Phi(\mathfrak{c}_X) \cdot \Gamma(X,Y) \dnu (X) \dnu(Y).
\]
Fix $\eta=1/16$ and split pairs into
\[
 |X-Y| \geq \eta\min(|X|,|Y|) \quad \hbox{and} \quad |X-Y|<\eta\min(|X|,|Y|).
\]

\emph{Separated pairs.} Let $L = \max(|X|,|Y|)$. For separated pairs, the bound $|X-Y|\geq\eta L/2$ follows from
the reverse triangle inequality if either $|X|$ or $|Y|$ is less than
$L/2$, and from $|X-Y|\geq\eta\min(|X|,|Y|)$ if both are at least $L/2$.
Formula \eqref{eq:kernelfactor} yields
\[
 |\mathcal H| \lesssim_\eta L^{-5}, \qquad |\nabla_X\mathcal H| \lesssim_\eta L^{-6}.
\]
Differentiating its prefactor in \eqref{eq:kernelfactor} explicitly,
\[
 \Gamma^r = \frac3{2\pi}s^2wr(\mathcal H+z\mathcal H_z), \qquad 
 \Gamma^z = -\frac3{2\pi}s^2wz(2\mathcal H+r\mathcal H_r).
\]
Consequently
\begin{equation}\label{eq:farkernel}
 \frac{|\Gamma(X,Y)|}{|X|} \lesssim_\eta s^2wL^{-5} \lesssim_\eta k(Y).
\end{equation}
The target point $X$ of $\nabla\Phi(\mathfrak{c}_X)$ lies in a fixed intermediate
sector, and $|\nabla\Phi(\mathfrak{c}_X)| \lesssim_\Phi |X|^{-1}$.
Equations \eqref{eq:bandmass} and \eqref{eq:farkernel} bound the
separated flux by $C_\Phi Q^{1/7} E^{3/7} b$.

\emph{Near pairs: estimate $\nabla_X G$.} If both $\Phi'(\mathfrak{c}_X)$ and $\Phi'(\mathfrak{c}_Y)$ vanish, the symmetrized integrand is zero. For every remaining near pair, at least one
of $\mathfrak{c}_X,\mathfrak{c}_Y$ belongs to $\operatorname{supp}\Phi'\subset[1/2,3/4]$.
Since
\[
 |\mathfrak{c}_X-\mathfrak{c}_Y|\leq\frac{2|X-Y|}{\min(|X|,|Y|)}<2\eta=1/8,
\]
both points lie in the fixed enlarged sector
$\widetilde{\mathcal S}:=\{X\in\Hh:3/8\leq \mathfrak{c}_X\leq7/8\}$.
Write $R = |X| \sim |Y|$ and
$d = |X-Y|$, so $r,s,z,w \sim R$. The local logarithmic
singularity and far-field decay of the unreflected stream kernel
are classical \cite[Lemmas~2.7--2.8]{FS15}.
Elementary angular integration,
using $D_-^2 \sim d^2+R^2\theta^2$ near $\theta = 0$, yields
\begin{equation}\label{eq:angularintegrals}
 \begin{aligned}
 \int_0^{2\pi}D_-^{-1}\dd\theta & \lesssim R^{-1}(1+\log(R/d)),\\
 \int_0^{2\pi}D_-^{-3}\dd\theta & \lesssim (Rd^2)^{-1},\\
 \int_0^{2\pi}(1-\cos\theta)D_-^{-3}\dd\theta & \lesssim R^{-3}(1+\log(R/d)).
 \end{aligned}
\end{equation}
Let $I_\pm=\int_0^{2\pi}\cos\theta\,D_\pm^{-1}\dd\theta$.
Then differentiating \eqref{eq:streamkernel}, we get
\[
 \begin{aligned}
 \partial_rG&=-\frac1{4\pi}
 \bigl[s(I_--I_+)+rs(\partial_rI_--\partial_rI_+)\bigr],\\
 \partial_zG&=-\frac{rs}{4\pi}(\partial_zI_--\partial_zI_+).
 \end{aligned}
\]
Since $D_+\geq z+w\sim R$, we have
$|I_+|\lesssim R^{-1}$ and
$|\nabla_XI_+|+|\nabla_YI_+|\lesssim R^{-2}$.
The reflected contribution to each first derivative of $G$ is
therefore $O(1)$. For the unreflected term,
\[
 \partial_rD_-^{-1}
   =-\bigl((r-s)+s(1-\cos\theta)\bigr)D_-^{-3}, \qquad
 \partial_zD_-^{-1}=-(z-w)D_-^{-3}.
\]
Using $|\cos\theta|\leq1$, $r,s\sim R$, and
$|r-s|+|z-w|\lesssim d$, we obtain
\[
 \begin{aligned}
 |\partial_rG|\lesssim{}&1
   +R\int_0^{2\pi}D_-^{-1}\dd\theta
   +R^2d\int_0^{2\pi}D_-^{-3}\dd\theta+R^3\int_0^{2\pi}(1-\cos\theta)D_-^{-3}\dd\theta,\\
 |\partial_zG|\lesssim{}&1
   +R^2d\int_0^{2\pi}D_-^{-3}\dd\theta.
 \end{aligned}
\]
The three integrals in the radial estimate arise, respectively,
from differentiating the prefactor $rs$, the squared difference $(r-s)^2$,
and the angular term $2rs(1-\cos\theta)$ in $D_-^2$.
By \eqref{eq:angularintegrals}, their contributions are bounded by
$1+\log(R/d)$, $R/d$, and $1+\log(R/d)$. Therefore, $|\nabla_X G|\lesssim R/d+ 1+ \log (R/d)$.

\emph{Near pairs: cancellation through symmetrization.} Now, let us estimate $\nabla_X G+ \nabla_Y G$.
To this end,  the difference
terms cancel:
\[
 (\partial_r+\partial_s)D_-^{-1}
 = -(r+s)(1-\cos\theta)D_-^{-3}, \qquad
 (\partial_z+\partial_w)D_-^{-1}=0.
\]
The derivative of the prefactor is
$(\partial_r+\partial_s)(rs)=r+s$, while its vertical derivative
vanishes. Hence
\[
 \begin{split}
 |\nabla_XG+\nabla_YG|\lesssim{}&1 + R\int_0^{2\pi}D_-^{-1}\dd\theta
 +R^3\int_0^{2\pi}(1-\cos\theta)D_-^{-3}\dd\theta.
 \end{split}
\]
Applying \eqref{eq:angularintegrals} again yields
\begin{equation}\label{eq:kernelderivatives}
 |\nabla_XG| \lesssim R/d+1+\log(R/d),
 \qquad 
 |\nabla_XG+\nabla_YG| \lesssim 1+\log(R/d).
\end{equation}

Let $W(X) = \nabla\Phi(\mathfrak{c}_X)/r$. On $\widetilde{\mathcal S}$,
$|W(X)|\lesssim_\Phi |X|^{-2}$ and
$|\nabla W(X)|\lesssim_\Phi |X|^{-3}$.
The segment joining a near pair stays in this enlarged sector at distance
comparable to $R$ from the origin. The mean value estimate therefore yields
\[
 |W(X)| \lesssim_\Phi R^{-2}, \qquad
 |W(X)-W(Y)| \lesssim_\Phi dR^{-3}.
\]
Interchanging source $Y$ and target $X$ to cancel a smooth-cutoff flux
singularity is also used for radial confinement in
\cite[proof of Lemma~3.2, pp.~132--134]{MM01} and
\cite[proof of Lemma~3.2, pp.~289--290]{Lim25}.
Here the cutoff is angular.
The near-pair region is symmetric, and $G(X,Y) = G(Y,X)$. Its flux is
therefore one half of the double integral of
\[
 \begin{split}
 &W(X)\cdot\nabla_X^\perp G+W(Y)\cdot\nabla_Y^\perp G\\
 & \quad = (W(X)-W(Y))\cdot\nabla_X^\perp G +W(Y)\cdot(\nabla_X^\perp G+\nabla_Y^\perp G).
 \end{split}
\]
By \eqref{eq:kernelderivatives}, the absolute value is at most
$C_\Phi R^{-2}(1+\log(R/d))$. For each fixed $X\in\widetilde{\mathcal S}$,
the contributing points $Y$ lie in $\widetilde{\mathcal S}$ and satisfy
$|X-Y|<\eta\min(|X|,|Y|)\leq\eta|X|$.
Lemma~\ref{lem:logmass}, with $\varrho=|X|$ and enlarged sector $\widetilde{\mathcal S}$,
therefore bounds the inner integral by $C_\Phi R^{-2}Q^{1/7}E^{3/7}$.
Since $k(X)$ is comparable to
$R^{-2}$ on the enlarged sector, the outer integral is at most
$C Q^{1/7}E^{3/7} b$.

All these calculations are valid with symmetric diagonal cutoffs.
The meridional singularities are locally integrable for bounded $q$;
dominated convergence removes the cutoffs.
\end{proof}

\section{Angular redistribution and the large strain estimate}\label{sec:pulse}
In this section, we derive the height-dissipation estimate and
use it to obtain angular redistribution and prove the large strain
estimate, as outlined in Figure~\ref{fig:strainmechanism}.

Consider Setup~\ref{setup:profile} (an evolving, scaled pair). 
Define three moments:
\begin{equation}\label{eq:angmoments}
 M_{\mathfrak{c}}(t)=\frac{1}{m}\int_{\Hh}\mathfrak{c}_X\dnu,\qquad
 M_\Phi(t) = \frac1{m}\int_{\Hh}\Phi(\mathfrak{c}_X)\dnu(X), \qquad M_z(t) = \frac1{m H(t)}\int_{\Hh}z\dnu.
\end{equation}
The normalized cosine angular moment $M_{\mathfrak{c}}$ is the average of $\mathfrak{c}_X=z/|X|$ with respect to the probability measure $\dnu/m$. 
The normalized cutoff angular moment $M_\Phi$ measures a smoothed fraction
of mass close to the positive $z$-axis. The normalized height moment $M_z$ measures its
mean vertical position after removing the affine compression. 

\subsection{Angular flux and the reference aspect ratio}

In this subsection, we extend the angular flux estimate (Proposition~\ref{prop:flux}) from Setup~\ref{setup:ringpair} (a fixed-time pair) to Setup~\ref{setup:profile} (an evolving scaled pair).

\begin{lemma}[Angular flux and the reference aspect ratio]
\label{lem:angularaspect}
Under Setup~\ref{setup:profile} (an evolving, scaled pair), let
$\Phi$ satisfy \eqref{eq:angularcutoff} and let $M_\Phi$ be defined
by \eqref{eq:angmoments}. Then, for $t\in I$,
\begin{equation}\label{eq:angderivative}
 |M_\Phi'(t)|\leq C_{f_0,\Phi}\eps(t)^{4/7}\bigl(b(t)+a(t)+\ell(t)\bigr),
 \qquad C_{f_0,\Phi}=C_\Phi F^{1/7}E_*^{3/7},
\end{equation}
where $C_\Phi$ depends only on the cutoff $\Phi$ and $E_*$ is defined in
\eqref{eq:Estar}. No smallness of the reference aspect ratio $\varepsilon$ or of
self-deformation is assumed.
\end{lemma}

The estimate below corresponds to Figure~\ref{fig:strainmechanism}(1). The angular-flux formula in the figure omits the contributions of the affine field $a(t)Dx$ and the error field $\mathfrak{e}(t)$, whose gradient is bounded by $\ell(t)$.

\begin{proof}
Transport of $\nu$ and conservation of $m$ imply
\begin{equation}\label{eq:angulartransport}
 M_\Phi'(t)=\frac1m\int_{\Hh}\nabla\Phi(\mathfrak{c}_X)\cdot
 (u+aDx+\mathfrak{e})\dnu.
\end{equation}
Let $\mathcal S=\{X\in\Hh:1/2\leq\mathfrak{c}_X\leq3/4\}$.
Since $\nabla\Phi(\mathfrak{c}_X)$ vanishes outside $\mathcal S$ and
$|X||\nabla\Phi(\mathfrak{c}_X)|\lesssim_\Phi1$, the bounds $|aDx|\leq2a|X|$
and $|\mathfrak{e}|\leq\ell|X|$ from the affine field and
\eqref{eq:errorcoordinates} control both external contributions.
Proposition~\ref{prop:flux} for $u$ and \eqref{eq:bandmass} for
$\nu(\mathcal S)$ therefore give
\begin{align}
 \frac1m\left|\int_{\Hh}\nabla\Phi(\mathfrak{c}_X)\cdot u\dnu\right|
 &\lesssim_\Phi \frac{Q^{1/7}E^{3/7}}mb,\label{eq:angularself}\\
 \frac1m\left|\int_{\Hh}\nabla\Phi(\mathfrak{c}_X)\cdot aDx\dnu\right|
 &\lesssim_\Phi a\frac{\nu(\mathcal S)}m\lesssim_\Phi \frac{Q^{1/7}E^{3/7}}ma,\label{eq:angularaffine}\\
 \frac1m\left|\int_{\Hh}\nabla\Phi(\mathfrak{c}_X)\cdot\mathfrak{e}\dnu\right|
 &\lesssim_\Phi\ell\frac{\nu(\mathcal S)}m\lesssim_\Phi \frac{Q^{1/7}E^{3/7}}m\ell.\label{eq:angularerror}
\end{align}
Finally, \eqref{eq:normalizedmass}--\eqref{eq:normalizedenergy}
and Lemma~\ref{lem:normalizedenergy} imply
\begin{equation}\label{eq:aspectconversion}
 \frac{Q^{1/7}E^{3/7}}m
 =(2\pi)^{3/7}F^{1/7}\cE^{3/7}\eps^{4/7}
 \leq(2\pi)^{3/7}F^{1/7}E_*^{3/7}\eps^{4/7}.
\end{equation}
Substituting \eqref{eq:angularself}--\eqref{eq:angularerror}
into \eqref{eq:angulartransport} and using
\eqref{eq:aspectconversion} proves \eqref{eq:angderivative},
with $(2\pi)^{3/7}$ absorbed into $C_\Phi$.
\end{proof}

\begin{remark}[Motivation behind $M_\Phi$]
\label{rem:mot_phi}
Our main idea in this paper is to study angular redistribution to estimate strain from below. Angular redistribution is naturally associated with the angular moment $M_{\mathfrak{c}}$, not necessarily with the cutoff angular moment $M_{\Phi}$. But we study $M_{\Phi}$ instead because the angular flux estimate \eqref{eq:angderivative} breaks down for $M_{\mathfrak{c}}$ since angular motion near the axis can also contribute at order $b$ without any small geometric factor like $\varepsilon$. 

Indeed, let us illustrate this with a simpler setup without the external field. We fix a nonnegative profile $f_0(\rho,\zeta)=f(\rho)g(\zeta)$ supported in $(1,2)^2$ and let $q_{\varepsilon,0}(r,z)=f(r/\varepsilon)g(z),z>0,$ with odd reflection. Then we can find $b(0)\sim \varepsilon^4$ whereas $|M_{\mathfrak{c}}'(0)|\sim \varepsilon^4$. And $|M_{\mathfrak{c}}'(0)|/b$ converges to a positive constant $C_{f_0}$ as $\varepsilon\to0$, which prevents an estimate of $|M_{\mathfrak{c}}'|$ in terms of $b$ with a small factor like $\varepsilon^{4/7}$. For the same data, $M_{\Phi}'(0)=0$  because $\Phi(\mathfrak{c})=1$ on the upper support when $\varepsilon$ is small.
\end{remark}

\subsection{Height dissipation and the axis kernel}\label{subsec:height-axis}

This subsection is devoted to derive the height-dissipation estimate explained in Section~\ref{sec:constructionidea} and Figure~\ref{fig:strainmechanism}. The first ingredient is a height identity (Lemma~\ref{lem:height}), which is a differential equation of a height moment involving an integral of the square of an axial velocity. The second ingredient is a relation of this axis integral with a height moment and the normalized cosine angular moment $M_{\mathfrak{c}}$ (Lemma~\ref{lem:axis}). Using these two, we derive the height-dissipation estimate in \eqref{eq:heightdissipation}.

\begin{lemma}[Height identity]\label{lem:height}
Under Setup~\ref{setup:evolution} (an evolving pair), let
$M_z^{\mathrm{ph}}=\int_{\Hh}z\dnu$. Then
\begin{equation}\label{eq:heightphysical}
 (M_z^{\mathrm{ph}})'+2aM_z^{\mathrm{ph}} = -\frac12\int_0^\infty|u^z(0,z)|^2\dd z -\int_{\Hh}\frac{|u^r|^2}{r}\dd r\dd z
 +\int_{\Hh}\mathfrak{e}^z\dnu.
\end{equation}
\end{lemma}

\begin{proof}
Transport of $\nu$ implies
$(M_z^{\mathrm{ph}})' = \int(u^z-2az+\mathfrak{e}^z)\dnu$.
The spatial integration-by-parts identity below is the one
underlying \cite[Eq.~(4.4)]{EY25}; their scalar vorticity
equals $-\omega$ in our sign convention.
Since $\omega = \partial_zu^r-\partial_ru^z$,
\[
 \begin{split}
 \int_{\Hh}u^z\dnu & = -\frac12\int_0^\infty|u^z(0,z)|^2\dd z +\int_{\Hh}(\partial_zu^z)u^r\dd r\dd z\\
 & = -\frac12\int_0^\infty|u^z(0,z)|^2\dd z -\int_{\Hh}\frac{|u^r|^2}{r}\dd r\dd z.
 \end{split}
\]
The first equality retains the axis term from $u^z\partial_ru^z$;
the boundary product $u^zu^r$ vanishes at $z=0$ by oddness of $u^z$.
The second uses $\partial_zu^z=-\partial_ru^r-u^r/r$ and
$u^r(0,z)=0$. The boundary terms at infinity vanish since
the Biot--Savart law and compactly supported vorticity imply $u(x)=O(|x|^{-2})$.
\end{proof}

\begin{lemma}[Axis mass and first moment]\label{lem:axis}
Under Setup~\ref{setup:ringpair} (a fixed-time pair), suppose that the pair is nonzero. Let
$v_{\rm ax}(z)=-u^z(0,z)$ and $M_z^{\mathrm{ph}}=\int_{\Hh}z\dnu$. Then
\begin{equation}\label{eq:axiskernel}
 v_{\rm ax}(z)=\int_{\Hh}\{P_s(z-w)-P_s(z+w)\}\dnu(s,w),
 \qquad P_h(s)=\frac{h^2}{2(s^2+h^2)^{3/2}}.
\end{equation}
In particular, $v_{\rm ax}\geq0$ on $z>0$, and
\begin{equation}\label{eq:axismoments}
 \int_0^\infty v_{\rm ax}\dd z=mM_{\mathfrak{c}},\qquad
 \int_0^\infty z v_{\rm ax}\dd z=M_z^{\mathrm{ph}}.
\end{equation}
Consequently,
\begin{equation}\label{eq:axisL2}
 \int_0^\infty v_{\rm ax}^2\dd z\geq
 \frac{4m^3M_{\mathfrak{c}}^3}{9M_z^{\mathrm{ph}}}.
\end{equation}
\end{lemma}

\begin{proof}
Evaluating \eqref{eq:bs} for $u^z$ on the axis and pairing the
odd reflection establishes \eqref{eq:axiskernel}, the physical
axis formula in \cite[Eq.~(4.6)]{EY25}.
Its positivity follows from $|z-w|\leq z+w$.
The first identity in \eqref{eq:axismoments} is \cite[Eq.~(4.7)]{EY25}, and the second is derived in \cite[proof of Proposition~5.1, Case~1]{EY25}.
For completeness, we derive both from \eqref{eq:axiskernel}.

For $h,\eta>0$, elementary integration yields
\[
 \begin{split}
 \int_0^\infty\{P_h(z-\eta)-P_h(z+\eta)\}\dd z
 =\frac{\eta}{\sqrt{\eta^2+h^2}},
 \qquad
 \int_0^\infty z\{P_h(z-\eta)-P_h(z+\eta)\}\dd z
 =\eta.
 \end{split}
\]
Thus \eqref{eq:axismoments} follows by Fubini's theorem.
For any nonnegative $w$ on $(0,\infty)$ with mass $M>0$, first moment $V>0$,
and $L>V/M$,
\[
 M-\frac VL\leq\int_0^Lw(z)(1-z/L)\dd z
 \leq\|w\|_{L^2(0,\infty)}(L/3)^{1/2}.
\]
The first inequality follows by writing
$M-V/L=\int_0^\infty w(z)(1-z/L)\dd z$ and discarding the
nonpositive integral over $(L,\infty)$.
Choosing $L=3V/M$ yields $\int_0^\infty w^2\dd z\geq4M^3/(9V)$;
apply this with $w=v_{\rm ax}$, $M=mM_{\mathfrak{c}}$, and $V=M_z^{\mathrm{ph}}$
to obtain \eqref{eq:axisL2}.
\end{proof}

\begin{remark}[Comparison with earlier height identities and axis estimates]
\label{rem:heightcomparison}
For evolution under the pair's own velocity ($a=\mathfrak{e}=0$),
Choi--Jeong proved that $M_z^{\mathrm{ph}}$ decreases
\cite[Lemma~3.4, Eq.~(3.3)]{CJ25}.  And in the case of $a=\mathfrak{e}=0$, Lemma~\ref{lem:height} reduces to the height identity
of Egamberganov--Yao \cite[Eq.~(4.4)]{EY25}.
Lemma~\ref{lem:axis} is a reformulation of their axis identities
and Cauchy--Schwarz argument
\cite[Eqs.~(4.6)--(4.7) and the proof of Proposition~5.1,
Case~1]{EY25}.

Here we extend the height evolution to include the prescribed
fields $aDx$ and $\mathfrak{e}$. Our main contribution in this section
is to combine the resulting normalized height-dissipation
estimate with the angular-flux and energy estimates.
This yields large strain at the origin through angular
redistribution, while allowing large self-deformation;
see Theorem~\ref{thm:pulse}.
\end{remark}

In Setup~\ref{setup:profile} (an evolving, scaled pair), we have
$M_z^{\mathrm{ph}}=mH M_z$, $H'=-2aH$, and $\lambda=m/H^2$.
Using $\mathfrak{e}^z\leq\ell z$, Lemmas~\ref{lem:height} and~\ref{lem:axis}
therefore imply

\[
 mH M_z'\leq-\frac{2m^3M_{\mathfrak{c}}^3}{9mH M_z}+\ell mH M_z.
\]
Hence we obtain the height-dissipation estimate
\begin{equation}\label{eq:heightdissipation}
 (M_z^2)'\leq-\frac49\lambda M_{\mathfrak{c}}^3+2\ell M_z^2
 \leq-\frac1{18}\lambda M_\Phi^3+2\ell M_z^2,
\end{equation}
where the last inequality uses $M_\Phi\leq2M_{\mathfrak{c}}$, since
$\Phi(c)\leq2c$. The negative term expresses height dissipation;
in the presence of the contribution of the error field, the estimate does not
assert that $M_z$ decreases.

\subsection{A large strain estimate}

The next theorem is devoted to the strain estimate in
Figure~\ref{fig:strainmechanism}. If the angular-redistribution
time $\tau$ reaches the prescribed threshold, which will be denoted by $\tau_*$, while the reference aspect $\varepsilon(t)$
and integrated error bound $\mathfrak{E}(t)$ remain controlled, the pair generates large
strain at the origin. 

Write
\[
 M_{z,0}:=M_z(0)=\int_{\Hh}\zeta\rho f_0\dd\rho\dd\zeta
\]
for the initial normalized height moment, and define
\begin{equation}\label{eq:pulseconstants}
 \tau_*=144e^2M_{z,0}^2,
\end{equation}
and 
\begin{equation}\label{eq:selfstrainconstant}
 c_*=\frac1{4C_{f_0,\Phi}},
\end{equation}
where $C_{f_0,\Phi}$ is the constant in Lemma~\ref{lem:angularaspect}.

Fix $\eps_*>0$, depending only on $f_0$ and $\Phi$, such that
\begin{equation}\label{eq:activationaspect}
 \eps_*\leq
 \min\left\{1,\frac{\sqrt7\,\zeta_-}{3\rho_+},
             \left(\frac{12c_*}{19}\right)^{7/4}\right\}.
\end{equation}
Thus $R_0/H_0\leq\eps_*$ implies
$z/(r^2+z^2)^{1/2}\geq3/4$ on the upper initial support and $M_\Phi(0)=1$,
while $(19/12)\eps_*^{4/7}\leq c_*$. If the angular-redistribution time $\tau$ reaches the threshold $\tau_*$, let $t_*$
denote the unique physical time satisfying $\tau(t_*)=\tau_*$.

\begin{theorem}[Large strain estimate]\label{thm:pulse}
Under Setup~\ref{setup:profile} (an evolving, scaled pair), take the cutoff $\Phi$
\eqref{eq:angularcutoff} and the constants $\tau_*,c_*,
\eps_*$ fixed by \eqref{eq:pulseconstants},
\eqref{eq:selfstrainconstant}, and \eqref{eq:activationaspect}.
Suppose $\tau(t_*)=\tau_*$ and
\[
 \mathfrak{E}(t) \leq 1, \qquad \eps(t) \leq \bar\eps \leq \eps_* \qquad (0 \leq t \leq t_*).
\]
Then
\begin{equation}\label{eq:pulse}
 \int_0^{t_*}b(t)\dd t \geq c_*\bar\eps^{-4/7}.
\end{equation}
The constants $\tau_*, c_*, \varepsilon_*$ depend only on $f_0$ and $\Phi$.
\end{theorem}

\begin{proof}
Write $B=\int_0^{t_*}b(t)\dd t$. By
Lemma~\ref{lem:angularaspect}, $b\geq0$, and
$\eps\leq\bar\eps$,
\[
 B\geq
 \frac{\bar\eps^{-4/7}}{C_{f_0,\Phi}}
 \int_0^{t_*}|M_\Phi'(t)|\dd t
 -
 \bar\eps^{-4/7}
 \int_0^{t_*}\eps(t)^{4/7}(a(t)+\ell(t))\dd t.
\]
Since $\eps'=3a\eps$ and $\mathfrak{E}(t_*)\leq1$,
\[
 \begin{aligned}
 \int_0^{t_*}\eps^{4/7}(a+\ell)\dd t
 &\leq
 \frac7{12}\bigl(\eps(t_*)^{4/7}-\eps(0)^{4/7}\bigr)
 +\bar\eps^{4/7}\mathfrak{E}(t_*)\leq\frac{19}{12}\bar\eps^{4/7}.
 \end{aligned}
\]
Thus the angular-flux and energy estimates imply
\[
 B\geq
 \frac{\bar\eps^{-4/7}}{C_{f_0,\Phi}}
 \int_0^{t_*}|M_\Phi'(t)|\dd t-\frac{19}{12}.
\]

It remains to use height dissipation and obtain angular redistribution.
The initial aspect condition ensures $M_\Phi(0)=1$.
If $M_\Phi>1/2$ throughout $[0,t_*]$, multiplying
\eqref{eq:heightdissipation} by $e^{-2\mathfrak{E}}$ and integrating
would imply
\[
 \begin{aligned}
 0\leq e^{-2\mathfrak{E}(t_*)}M_z(t_*)^2
 &\leq M_{z,0}^2
 -\frac1{18}\int_0^{t_*}
 e^{-2\mathfrak{E}(t)}\lambda(t)M_\Phi(t)^3\dd t\\
 &<M_{z,0}^2-\frac{\tau_*}{144e^2}=0,
 \end{aligned}
\]
a contradiction. Here we used $\mathfrak{E}\leq1$ and
$\tau(t_*)=\tau_*=144e^2M_{z,0}^2$.
By continuity, $M_\Phi(t_1)=1/2$ for some $t_1\leq t_*$,
and hence
\[
 \int_0^{t_*}|M_\Phi'(t)|\dd t
 \geq |M_\Phi(t_1)-M_\Phi(0)|=\frac12.
\]
Substitution into the strain lower bound yields
\[
 B\geq\frac1{2C_{f_0,\Phi}}\bar\eps^{-4/7}-\frac{19}{12}
 =2c_*\bar\eps^{-4/7}-\frac{19}{12}
 \geq c_*\bar\eps^{-4/7},
\]
where the last inequality uses
$\bar\eps\leq\eps_*$ and
$(19/12)\eps_*^{4/7}\leq c_*$.
\end{proof}

\section{Strain-rate capacity and conditional upper bounds}
\label{sec:strainbounds}

As explained in Section~\ref{sec:constructionidea} around \eqref{eq:capacity_bd}, we also need to control our solution from above in order to preserve its desired structures. In this Section~\ref{sec:strainbounds}, we derive a key ingredient for that goal in a single-pair setup (Setup~\ref{setup:ringpair} or Setup~\ref{setup:evolution}): a strain-rate capacity bound. Using this bound, we then estimate $A(t)=\int_0^t a(s) ds$ and $\varepsilon$ from above.

\subsection{The strain-rate capacity}

\begin{lemma}[Strain-rate capacity]\label{lem:capacity}
Under Setup~\ref{setup:ringpair} (a fixed-time pair),
\begin{equation}\label{eq:capacity}
 0\leq b\leq Q^{2/3}m^{1/3}.
\end{equation}
\end{lemma}

\begin{proof}
Use the strain-rate formula \eqref{eq:source}.
For a cutoff distance $\varrho>0$ from the origin, write
$b=b_{\rm in}+b_{\rm out}$, where
\[
 b_{\rm in}:=\int_{\Hh\cap\{r^2+z^2<\varrho^2\}}k(r,z)\dnu,
 \qquad
 b_{\rm out}:=\int_{\Hh\cap\{r^2+z^2\geq\varrho^2\}}k(r,z)\dnu.
\]
These contributions satisfy
\[
 b_{\rm in}\leq\frac38Q\varrho,
 \qquad
 b_{\rm out}\leq\frac{m}{\sqrt3\,\varrho^2}.
\]
The first inequality follows by polar integration and the second
from $k(r,z)(r^2+z^2)\leq1/\sqrt3$.
For $Qm>0$, take $\varrho=(m/Q)^{1/3}$
and use $3/8+1/\sqrt3<1$. The zero case of $Qm=0$ is immediate.
\end{proof}

For positive initial scales $Q_*,R_0,H_0$, define the capacity
\begin{equation}\label{eq:singlecapacity}
 P:=Q_*(R_0^2H_0)^{1/3}.
\end{equation}
Under Setup~\ref{setup:profile} (an evolving, scaled pair),
\eqref{eq:normalizedmass} and Lemma~\ref{lem:capacity} imply
\[
 Q^{2/3}m^{1/3}=F^{2/3}P,
 \qquad
 0\leq b(t)\leq F^{2/3}P.
\]
This upper bound depends only on the initial data.

\subsection{Conditional bounds on strain and reference aspect}

\begin{lemma}[Conditional upper bounds]\label{lem:singleclockbounds}
Let $I\subset[0,\infty)$ be a time interval containing $0$.
Take $Q_*,R_0,H_0>0$ and let $a\in C(I)$ satisfy
$0\leq a(t)\leq\bar a$ on $I$, where $\bar a\geq0$ is a fixed constant.
Set $A(t)=\int_0^t a(s)\dd s$ and define
$R,H,\eps,\lambda,\tau$ by \eqref{eq:referencescales}
and $P$ by \eqref{eq:singlecapacity}.
Write $\lambda_0=\lambda(0)$ and $\eps_0=\eps(0)$.
Then, for every $t\in I$,
\begin{equation}\label{eq:singleclockcomparison}
 \lambda(t)\leq\lambda_0+4\bar a\,\tau(t).
\end{equation}
Consequently, if $\tau(t)\leq\tau_*$, then
\begin{equation}\label{eq:singleclockbounds}
 \begin{gathered}
 A(t)\leq\frac14\log\left(1+
               \frac{4\bar a\,\tau_*}{\lambda_0}\right), 
               \qquad
 \eps(t)\leq
       \left(\eps_0^{4/3}+\frac{4\bar a\,\tau_*}{P}\right)^{3/4}.
 \end{gathered}
\end{equation}
\end{lemma}

\begin{proof}
The definitions of the reference scales imply
\begin{equation}\label{eq:singlereferenceidentities}
 R^2H=R_0^2H_0,
 \qquad
 \lambda=\lambda_0e^{4A}=P\eps^{4/3},
 \qquad
 \lambda'=4a\lambda.
\end{equation}
Therefore,
\[
 \lambda(t)-\lambda_0
 =4\int_0^t a(s)\lambda(s)\dd s
 \leq4\bar a\,\tau(t).
\]
This proves \eqref{eq:singleclockcomparison}.
If $\tau(t)\leq\tau_*$, dividing by $\lambda_0$ or $P$
and using \eqref{eq:singlereferenceidentities} yields
\eqref{eq:singleclockbounds}.
\end{proof}

Currently, this upper bound is conditional on the bound $a(t)\leq \bar a$.  But, in the multiscale application later on, the strain-rate capacity bound \eqref{eq:capacity} will supply this bound; see Lemma~\ref{lem:clockbounds}.

\section{Support bounds}\label{sec:geometry}
For our multiscale solution, ring pairs of different scales will remain well separated over time, as the suitable setup explained in Section~\ref{sec:constructionidea}.

In this Section~\ref{sec:geometry}, we now establish support estimates under a single pair setup (Setup~\ref{setup:evolution} or Setup~\ref{setup:profile}), which will be needed for scale separation in Section~\ref{sec:cascade}.  

\subsection{Velocity interpolation and an upper bound}

Without assuming small self-deformation, the reference lengths
$R,H$ do not by themselves control the actual support. We estimate its maximal radial and vertical extents directly. 

To this end, we first estimate the velocity through interpolation in terms of a vorticity bound and energy, which in turn can be estimated without requiring $f(t,\rho,\zeta)$ remain close to $f_0$.

\begin{lemma}[Velocity interpolation]\label{lem:velocityinterpolation}
For a decaying divergence-free velocity on $\R^3$ with
finite energy and bounded vorticity,
\begin{equation}\label{eq:velocityinterpolation}
 \|u\|_{L^\infty(\R^3)} \lesssim \|\curl u\|_{L^\infty(\R^3)}^{3/5}\|u\|_{L^2(\R^3)}^{2/5}.
\end{equation}
\end{lemma}

\begin{proof}
Let $K$ be the Biot--Savart kernel, so that
$|K(y)|\lesssim |y|^{-2}$ and $|\nabla K(y)|\lesssim |y|^{-3}$.
Choose a smooth cutoff $\chi$ equal to one on the unit ball and
zero outside the ball of radius two, and take
$\chi_\varrho(y)=\chi(y/\varrho)$. Split the kernel as
$K=\chi_\varrho K+(1-\chi_\varrho)K$.
The near contribution satisfies
\[
 \| (\chi_\varrho K)*\curl u\|_{L^\infty(\R^3)}
 \lesssim \|\curl u\|_{L^\infty(\R^3)}
 \int_{|y|\leq2\varrho}|y|^{-2}\dd y
 \lesssim \varrho\|\curl u\|_{L^\infty(\R^3)}.
\]
For the far contribution, integration by parts transfers the curl
from $u$ to $(1-\chi_\varrho)K$. The differentiated kernel obeys
\[
 \bigl|\nabla\bigl((1-\chi_\varrho)K\bigr)(y)\bigr|
 \lesssim |y|^{-3}\mathbf 1_{\{|y|\geq\varrho\}},
 \qquad
 \|\nabla((1-\chi_\varrho)K)\|_{L^2(\R^3)}
 \lesssim\varrho^{-3/2}.
\]
Here the cutoff derivative is supported on
$\varrho\leq|y|\leq2\varrho$ and satisfies the same bound.
The integration by parts is justified by truncation at infinity;
Cauchy--Schwarz therefore bounds the far contribution by
$C\varrho^{-3/2}\|u\|_{L^2(\R^3)}$. Consequently,
\[
 \|u\|_{L^\infty(\R^3)} \lesssim \varrho\|\curl u\|_{L^\infty(\R^3)}+\varrho^{-3/2}\|u\|_{L^2(\R^3)}.
\]
When both norms on the right are positive, choose
$\varrho=(\|u\|_{L^2(\R^3)}/\|\curl u\|_{L^\infty(\R^3)})^{2/5}$
to balance the two terms and obtain \eqref{eq:velocityinterpolation}.
The zero cases follow directly from the same bound.
\end{proof}

Now, to bound the support from above, we use Lemma~\ref{lem:velocityinterpolation} to estimate the self-induced velocity and combine it with the contribution of the prescribed external field $aDx+\mathfrak{e}$.

\begin{proposition}[Upper bound of support]\label{prop:confinement}
Under Setup~\ref{setup:profile} (an evolving, scaled pair), with arbitrary positive initial scales $Q_*,R_0,H_0$, let $V_0 := Q_*R_0^{6/5}H_0^{1/5}$.
Then
\begin{equation}\label{eq:confinement}
 {\rho^+}(t):= \max_{\supp\vectomega(t)}|x| \lesssim \{H_0e^{-2A}+R_0e^A+V_0^{5/2}t^{5/2}e^{2A}\}.
\end{equation}
\end{proposition}

\begin{proof}
Let $\mathcal R(t)$ and $\mathcal Z(t)$ be the maximal radial
and absolute vertical coordinates of the support of $\vectomega(t)$.
The conserved scalar bound implies
$\|\vectomega\|_{L^\infty(\R^3)}\leq Q_*F\mathcal R$.
The normalized-energy bound \eqref{eq:Estar} and Lemma~\ref{lem:velocityinterpolation} yield
\begin{equation}\label{eq:Vforce}
 \|u(t)\|_{L^\infty(\R^3)} \lesssim V_0e^{4A(t)/5}\mathcal R(t)^{3/5}.
\end{equation}
On compact classical intervals the support extrema are locally
Lipschitz. The trajectory equations and the envelope argument for
extrema over fixed compact initial labels (see \cite[Lemma~4.4]{BC26})
imply, with the upper right Dini derivative
$D^+f(t)=\limsup_{h\to0^+}\frac{f(t+h)-f(t)}h$,
\[
 \begin{aligned}
 D^+\mathcal R&\leq(a+\ell)\mathcal R+\|u\|_{L^\infty(\R^3)},\\
 D^+\mathcal Z&\leq(-2a+\ell)\mathcal Z+\|u\|_{L^\infty(\R^3)}.
 \end{aligned}
\]
These inequalities hold almost everywhere with ordinary derivatives.
Let
\[
 \mathcal Z_{\rm hom}(t)=\mathcal Z(0)\exp\bigl(-2A(t)+\mathfrak{E}(t)\bigr),\qquad
 \mathcal M(t)=\max\{\mathcal R(t),\mathcal Z(t)-\mathcal Z_{\rm hom}(t)\}.
\]

The function $\mathcal M$ is positive and absolutely continuous,
with $\mathcal M(0)=\mathcal R(0)$. Since
$\mathcal Z_{\rm hom}'=(-2a+\ell)\mathcal Z_{\rm hom}$,
whether $\mathcal M=\mathcal R$ or $\mathcal M=\mathcal Z-\mathcal Z_{\rm hom}$, the preceding differential inequalities, together with \eqref{eq:Vforce} and $a\geq0$, imply, almost everywhere,
\begin{equation}\label{eq:DE_upper}
 \mathcal M'\leq(a+\ell)\mathcal M+CV_0e^{4A/5}\mathcal M^{3/5},
\end{equation}
which we can rewrite as
\[
 (\mathcal M^{2/5})'\leq\frac25(a+\ell)\mathcal M^{2/5}+CV_0e^{4A/5},
\]

Multiplying the inequality by
$\exp\bigl(-\frac25(A+\mathfrak{E})\bigr)$, integrating, and using $0\leq\mathfrak{E}\leq1$, we obtain
\[
 \mathcal M(t)^{2/5}\lesssim e^{2A(t)/5}
 \left(R_0^{2/5}+V_0\int_0^t e^{2A(s)/5}\dd s\right).
\]
Since $A$ is nondecreasing,
\[
 \mathcal M(t)\lesssim R_0e^{A(t)}+V_0^{5/2}t^{5/2}e^{2A(t)}.
\]
Finally, $\rho^+\leq\mathcal R+\mathcal Z\leq2\mathcal M+\mathcal Z_{\rm hom}$,
$\mathcal Z(0)\lesssim H_0$, and $\mathfrak{E}\leq1$ establish
\eqref{eq:confinement}.
\end{proof}

\subsection{Linear approximation and a lower bound}

An upper bound on the support distances alone cannot separate adjacent
scales; a positive lower bound is also needed. We retain the linear part $bDx$ of the pair's induced velocity at the origin and estimate the quadratic remainder
$u-bDx$. 

\begin{lemma}[Linear approximation]\label{lem:origin}
Under Setup~\ref{setup:ringpair} (a fixed-time pair), with self-induced
velocity $u$, strain rate $b$ at the origin, and $\|q\|_{L^\infty(\R^3)}\leq Q$,
\begin{equation}\label{eq:originremainder}
 |u(x)-b Dx| \lesssim Q|x|^2, \qquad \forall x\in\R^3.
\end{equation}
\end{lemma}

\begin{proof}
We fix a time and suppress it in the notation. Recall that $\vectomega(x)$ is even. Indeed, in Cartesian coordinates,
\[
 \vectomega(x)=q(x)(x_2,-x_1,0),\qquad q(-x)=-q(x),
\]
where the second identity follows from axisymmetry and oddness in
$z$. Both factors change sign, so
$\vectomega(-x)=\vectomega(x)$. Also,
$|\vectomega(y)|\leq Qr_y\leq Q|y|$, where
$r_y=(y_1^2+y_2^2)^{1/2}$.

Write the Biot--Savart law using the vector kernel
\[
 \mathcal K(z)=\frac{z}{4\pi|z|^3},\qquad
 u(x)=\int_{\R^3}\vectomega(y)\times\mathcal K(x-y)\dd y.
\]
The kernel is odd, and its Cartesian derivatives satisfy
$|\nabla^k\mathcal K(z)|\lesssim_k|z|^{-2-k}$.
Below, $\partial_i\mathcal K$ denotes differentiation in the
$i$-th coordinate of the kernel argument; $D$ continues to denote
only the matrix $\diag(1,1,-2)$.
Differentiating at the origin and using Lemma~\ref{lem:originstrain},
we obtain the linear term
\[
 bDx = \sum_{i=1}^3x_i\partial_i u(0) = \sum_{i=1}^3x_i\int_{\R^3}
 \vectomega(y)\times\partial_i\mathcal K(-y)\dd y.
\]
These integrals converge absolutely: near the origin the integrands
are bounded by $CQ|y|^{-2}$, which is locally integrable in
three dimensions, and the vorticity is compactly supported.

The assertion at $x=0$ follows from $u(0)=0$, since an even
vorticity and an odd kernel produce an odd velocity. Now fix $x\ne0$
and take $\varrho_x=|x|$. We split both the velocity and its linear term
into $|y|\leq2\varrho_x$ and $|y|>2\varrho_x$.
On the first region, $|\vectomega(y)|\leq2Q\varrho_x$ and
$|x-y|\leq3\varrho_x$, so
\[
 \left|\int_{|y|\leq2\varrho_x}\vectomega(y)\times\mathcal K(x-y)\dd y\right|
 \lesssim Q\varrho_x\int_{|z|\leq3\varrho_x}|z|^{-2}\dd z
 \lesssim Q\varrho_x^2.
\]
Similarly, the contribution to the linear term is bounded by
\[
 \left|\sum_{i=1}^3x_i\int_{|y|\leq2\varrho_x}
 \vectomega(y)\times\partial_i\mathcal K(-y)\dd y\right|
 \lesssim Q\varrho_x\int_{|y|\leq2\varrho_x}|y|^{-2}\dd y \lesssim Q\varrho_x^2.
\]

On $|y|>2\varrho_x$, thanks to Taylor's formula,
\[
 \mathcal K(x-y)=\mathcal K(-y)
 +\sum_{i=1}^3x_i\partial_i\mathcal K(-y)
 +\frac12\sum_{i,j=1}^3x_ix_j\partial_i\partial_j\mathcal K(-y)
 +\mathcal K_{\rm rem}(x,y),
\]
where
\[
 |\mathcal K_{\rm rem}(x,y)|\lesssim |x|^3|y|^{-5}.
\]
Indeed, the segment $-y+s x$, $0\leq s\leq1$,
stays at distance at least $|y|/2$ from the origin, so the bound
follows from the third derivatives of the kernel.
Both $\mathcal K$ and $\partial_i\partial_j\mathcal K$ are odd.
Their cross products with the even function $\vectomega(y)$
therefore integrate to zero on the symmetric region $|y|>2\varrho_x$.
Thus the constant and quadratic Taylor terms vanish after integration,
and the linear term is exactly the far-field part of $bDx$.
The remaining contribution is bounded by
\[
 \begin{aligned}
 \left|\int_{|y|>2\varrho_x}\vectomega(y)\times\mathcal K_{\rm rem}(x,y)\dd y\right|
 &\lesssim Q\varrho_x^3\int_{|y|>2\varrho_x}|y|^{-4}\dd y\\
 &\lesssim Q\varrho_x^3\int_{2\varrho_x}^\infty s^{-2}\dd s
 \lesssim Q\varrho_x^2.
 \end{aligned}
\]
Combining the near- and far-field bounds, we obtain
\eqref{eq:originremainder}. 
\end{proof}

Outer-region strain-rate approximations that are relevant to the lemma above also appear in
\cite[Lemma~2.1]{KJ22} and \cite[Proposition~2.1]{BC26}; their
hypotheses and error bounds differ from the estimate above.

\begin{proposition}[Lower bound of support]\label{prop:gap}
Consider Setup~\ref{setup:evolution} (an evolving pair).
Let ${\rho^-}(t) := \min_{\supp\vectomega(t)}|x|$, with ${\rho^-}(0)>0$, and 
\[
 A_{\rm tot}(t) := \int_0^t(a(s)+b(s))\dd s=A(t)+B(t).
\]
Then
\begin{equation}\label{eq:gap}
 {\rho^-}(t) \geq \frac{{\rho^-}(0)\exp\bigl(-2A_{\rm tot}(t)-\mathfrak{E}(t)\bigr)} {1+CQ{\rho^-}(0)t}.
\end{equation}
\end{proposition}

\begin{proof}
Fix an initial point $x_0\in\supp\vectomega(0)$. Its trajectory $x(t)$ is defined by the full transporting velocity:
\begin{equation}\label{eq:trajectory_eq}
 \frac{\dd x(t)}{\dd t}
 =u(t,x(t))+a(t)Dx(t)+\mathfrak{e}(t,x(t)),\qquad x(0)=x_0.
\end{equation}
As $x_0$ ranges over the initial support, the points $x(t)$ fill $\supp\vectomega(t)$.
By Lemma~\ref{lem:origin}, the self-induced velocity equals $bDx$
up to an error bounded by $CQ|x|^2$. Moreover,
\[
x\cdot Dx=x_1^2+x_2^2-2x_3^2\geq-2|x|^2,
\]
and
$|\mathfrak{e}(t,x)|\leq\ell(t)|x|$, since $\mathfrak{e}(t,0)=0$.
Taking the scalar product of the trajectory equation~\eqref{eq:trajectory_eq} with
$x(t)/|x(t)|$ therefore yields
\begin{equation}
\label{eq:DE_lower}
 \frac{\dd}{\dd t}|x(t)|
 \geq-[2a+2b+\ell]|x(t)|-CQ|x(t)|^2.
\end{equation}
On a compact classical interval the trajectory stays away from the
origin. From the this differential inequality, we can obtain one of its reciprocal $1/|x(t)|$:
\[
\frac{d}{dt} \left( \frac{1}{|x(t)|} \right) \leq \frac{2a+2b+\ell}{|x(t)|}+CQ.
\]

Multiplying it by the integrating
factor $\exp\bigl(-2A_{\rm tot}-\mathfrak{E}\bigr)$, we obtain
\[
 \frac{\dd}{\dd t}
 \left(\frac{\exp\bigl(-2A_{\rm tot}(t)-\mathfrak{E}(t)\bigr)}{|x(t)|}\right)
 \leq CQ\exp\bigl(-2A_{\rm tot}(t)-\mathfrak{E}(t)\bigr)\leq CQ.
\]
Here $a,b,\ell\geq0$, so the integrating factor is at most one.
Integrating from $0$ to $t$ we get
\[
 |x(t)|\geq
 \frac{|x_0|\exp\bigl(-2A_{\rm tot}(t)-\mathfrak{E}(t)\bigr)}{1+CQ|x_0|t}.
\]
For fixed $t\geq0$, the function $s\mapsto s/(1+CQt s)$
is increasing on $s>0$. Since $|x_0|\geq\rho^-(0)$, we can take
the minimum over all initial points obtaining \eqref{eq:gap}.
Since we integrate along individual trajectories before taking this minimum,
no differentiability of $\rho^-$ is needed. The strain $B$ generated by the vortex-ring pair is included in $A_{\rm tot}=A+B$.
\end{proof}

We clarify an important structural difference between the upper \eqref{eq:confinement} and lower bounds \eqref{eq:gap}  of the support. The upper bound \eqref{eq:confinement} does not depend on the solution $q(t)$ whereas the lower bound does depend on the solution because it involves $B(t)$. The main reason of the difference is that the self contribution enters the dynamics at a leading order only for the lower bound~\eqref{eq:gap}. Compare the leading order with the order of the self-induced terms in \eqref{eq:DE_upper} and \eqref{eq:DE_lower}.

We can in fact estimate $B$ further and remove the dependence on the solution. However, we retain the current lower bound, which will suffice for its application to our multiscale solution. See \eqref{eq:comparisoncancellation} and the paragraph after this equation.

\endgroup

\begingroup

\section{Interaction estimates}
\label{sec:interactions}
The local interaction error refers to the error of the approximation underlying the main decomposition~\eqref{eq:main_dec}. 
Here we estimate  the local interaction error based on scale separation. We will mainly consider Setup~\ref{setup:ringpair} and estimate the local interaction error in terms of relative separation between the source and observation regions.

In addition, we will extend the local error fields globally to construct the error field $\mathfrak{e}$ required by Setup~\ref{setup:evolution}. To this end, we will use a special kind of stream functions and cutoffs. 
Let us first collect their basic properties below.

\subsection{Stream function preliminaries}\label{subsec:streamfunction}

For $X=(r,z) \in \mathbb{R}^2$, define $h$ by
\begin{equation}\label{eq:representation}
h(X):=\frac3{2\pi}\int_{\Hh}s^2w\mathcal H(X,Y)\dnu(Y), \quad Y=(s,w).
\end{equation}
Here, we allow signed $r$. We use the letter $h$ due to its relation with $\mathcal{H}$.
Now using \eqref{eq:reppsi} and \eqref{eq:kernelfactor}, it holds
\[
\psi(X)= -r^2 z h(X), \qquad X=(r,z)\in \mathbb{H}.
\]
We write $\nabla_X=(\partial_r,\partial_z)$ for meridional
derivatives of scalar functions and kernels, including the radial
cutoffs below. Gradients of velocity fields are taken in Cartesian coordinates.

Let $\Sigma=\supp (q|_{\mathbb{H}}) \subset (0,\infty)^2$, and
\[
\Sigma_{\mathrm{sym}}
=\{
(\epsilon_1 r, \epsilon_2 z):
(r,z)\in \Sigma, \epsilon_1,\epsilon_2\in \{-1,1\}
\}.
\]

\begin{lemma}[Regularity and parity of $h$]\label{lem:hregularity}
 On $\mathbb{R}^2\setminus \Sigma_{\mathrm{sym}}$, the function $h$ is $C^2$ and even in each coordinate.
At the origin, $h(0)=b$ and $\nabla_Xh(0)=0$.
\end{lemma}

\begin{proof}
Use the integral in \eqref{eq:representation} for signed $r,z$ together with \eqref{eq:kernelfactor}.
On every compact set disjoint from $\Sigma_{\mathrm{sym}}$,
both $D_-^2$ and $D_+^2$ have a uniform positive lower bound for source
points $Y$ in the support $\Sigma$ and all $\theta$. Thus $h$ is $C^2$ on $\mathbb{R}^2\setminus \Sigma_{\mathrm{sym}}$.

Replacing $r$ by $-r$ amounts to shifting $\theta$
by $\pi$, and replacing $z$ by $-z$ interchanges $t$ and $1-t$ in
\eqref{eq:kernelfactor}. Hence $h$ is even. 
Thanks to $\mathcal H(0,Y)=\pi|Y|^{-5}$ and \eqref{eq:source} we have $h(0)=b$,
and due to evenness we also have $\nabla_Xh(0)=0$.
\end{proof}

If
$|X| \leq |Y|/2$ or $|Y| \leq |X|/2$, then, for $k = 0,1,2$,
\begin{equation}\label{eq:separatedkernelderivatives}
 (1-t)D_-^2+tD_+^2 \geq (|X|-|Y|)^2,
 \qquad 
 |\nabla_X^k\mathcal H(X,Y)| \lesssim_k\max(|X|,|Y|)^{-5-k}.
\end{equation}
The derivative bounds follow by differentiating the integral in
\eqref{eq:kernelfactor}.

For any $h\in C^2(\R^2)$ even in $r$ and $z$, denote by
$\mathcal V[h]$ the velocity generated by the stream function $-r^2zh$.
Its components are defined by
\[
 \mathcal V[h]^r := r^{-1}\partial_z (r^2zh)= r(h+zh_z), 
 \qquad \mathcal V[h]^z :=
 -r^{-1}\partial_r (r^2 z h) =
 -z(2h+rh_r).
\]
By evenness in $r$, this field is $C^1$ in Cartesian coordinates,
including at the axis. It is divergence-free and has the prescribed
plane parity. Differentiating, we obtain
\begin{equation}\label{eq:streamvelocitybounds}
 \begin{split}
 |\mathcal V[h]|& \lesssim |X|\bigl(|h|+|X||\nabla_X h|\bigr),\\
 |\nabla \mathcal V[h]|& \lesssim \bigl(|h|+|X||\nabla_X h|+|X|^2|\nabla_X^2h|\bigr).
 \end{split}
\end{equation}

As described above at the beginning of this Section~\ref{sec:interactions}, we will cut off stream functions. To this end,
fix once and for all smooth cutoff functions
$\chi^{\mathrm{out}},\chi^{\mathrm{in}}:[0,\infty)\to[0,1]$ satisfying the following conditions:
\begin{equation}\label{eq:radialcutoffs}
 \begin{gathered}
 \chi^{\mathrm{out}}(s)=1\quad(s\leq3/2),\qquad
 \chi^{\mathrm{out}}(s)=0\quad(s\geq2),\\
 \chi^{\mathrm{in}}(s)=0\quad(s\leq1/4),\qquad
 \chi^{\mathrm{in}}(s)=1\quad(s\geq1/2),\\
 \chi^{\mathrm{out}}_{\overline\rho}(X)
   =\chi^{\mathrm{out}}(|X|/\overline\rho),\qquad
 \chi^{\mathrm{in}}_{\underline\rho}(X)
   =\chi^{\mathrm{in}}(|X|/\underline\rho).
 \end{gathered}
\end{equation}

\subsection{An outer vortex-ring pair}

In the next two lemmas, the positive radii $\overline\rho$ and
$\underline\rho$ are, respectively, an upper bound on the size of the
inner region and a lower bound on the distance of the outer region
from the origin. In their application, these radii refer to different
vortex-ring pairs.
We denote the error fields by $\mathfrak{e}^{\mathrm{out}}$ and
$\mathfrak{e}^{\mathrm{in}}$. The superscripts describe the source vortex-ring
pair's position relative to the observation region: the outer field $\mathfrak{e}^{\mathrm{out}}$
extends $u-bDx$, whereas the inner field $\mathfrak{e}^{\mathrm{in}}$ extends $u$ on the distant
observation region.

\begin{lemma}[An outer vortex-ring pair]\label{lem:outer}
Under Setup~\ref{setup:ringpair} (a fixed-time pair), regard the pair $q$ as a source and suppose
$\supp\vectomega\cap B_{\underline\rho}=\varnothing$ and
$4\overline\rho\leq\underline\rho$.
\begin{enumerate}
\item Its induced velocity $u$ satisfies
\begin{equation}\label{eq:outertaylor}
 \|\nabla(u-b Dx)\|_{L^\infty(B_{2\overline\rho})} \lesssim b(\overline\rho/\underline\rho)^2.
\end{equation}
\item There is a global divergence-free axisymmetric velocity
$\mathfrak{e}^{\mathrm{out}}\in C^1(\R^3)$, without swirl and with the prescribed parity, such that
\[
 \mathfrak{e}^{\mathrm{out}} = u-b Dx\text{ on }B_{3\overline\rho/2},
 \qquad
 \|\nabla \mathfrak{e}^{\mathrm{out}}\|_{L^\infty(\R^3)}
 \lesssim b(\overline\rho/\underline\rho)^2.
\]
\end{enumerate}
\end{lemma}

\begin{proof}
For $|X| \leq \underline\rho/2$, the kernel derivative bound
\eqref{eq:separatedkernelderivatives} implies
\[
 |\nabla_X^2h(X)| \lesssim \int_{\Hh}s^2w|Y|^{-7}\dnu(Y) \lesssim b(\underline\rho)^{-2}.
\]
Since $h(0) = b$ and $\nabla_X h(0) = 0$,
Taylor's theorem yields
\[
 |h(X)-b|+|X||\nabla_X h(X)| +|X|^2|\nabla_X^2h(X)| \lesssim b|X|^2(\underline\rho)^{-2}.
\]
Applying the gradient bound in \eqref{eq:streamvelocitybounds}
to $h-b$, and using $u-b Dx = \mathcal V[h-b]$, we obtain
\eqref{eq:outertaylor}.

Using the cutoff \eqref{eq:radialcutoffs}, define
$\widetilde h = \chi^{\mathrm{out}}_{\overline\rho}(h-b)$, extended by
zero outside $B_{2\overline\rho}$, and let $\mathfrak{e}^{\mathrm{out}} = \mathcal V[\widetilde h]$.
The preceding bounds, the product rule, and
\eqref{eq:streamvelocitybounds} imply
\[
 \|\nabla \mathfrak{e}^{\mathrm{out}}\|_{L^\infty(\R^3)} \lesssim b(\overline\rho/\underline\rho)^2.
\]
This $C^1$ divergence-free field has the prescribed parity and
equals $u-bDx$ on $B_{3\overline\rho/2}$.
\end{proof}

\subsection{An inner vortex-ring pair}

For an inner source, its conserved mass and support size control
the whole distant field. 

\begin{lemma}[An inner vortex-ring pair]\label{lem:inner}
Under Setup~\ref{setup:ringpair} (a fixed-time pair), let the pair $q$ be a source supported
in $B_{\overline\rho}$, with mass $m$.
\begin{enumerate}
\item Its induced velocity $u$ satisfies
\begin{equation}\label{eq:innerfar}
 |\nabla_x^k u(x)|\lesssim_k m \overline\rho^3|x|^{-4-k}, \qquad \text{for} \quad |x| \geq 2\overline\rho,\ k = 0,1.
\end{equation}
\item If $\underline\rho \geq 16\overline\rho$, there is a global divergence-free
axisymmetric velocity $\mathfrak{e}^{\mathrm{in}}\in C^1(\R^3)$, without swirl and with the
prescribed parity, equal to $u$ outside $B_{\underline\rho/2}$, zero on
$B_{\underline\rho/4}$, and satisfying
\begin{equation}\label{eq:innerextension}
 (\underline\rho)^{-1}\|\mathfrak{e}^{\mathrm{in}}\|_{L^\infty(\R^3)}+\|\nabla \mathfrak{e}^{\mathrm{in}}\|_{L^\infty(\R^3)} \lesssim m \overline\rho^3(\underline\rho)^{-5}.
\end{equation}
\end{enumerate}
\end{lemma}

\begin{proof}
Since $s^2w \leq \overline\rho^3$ on the support, the kernel bounds
\eqref{eq:separatedkernelderivatives} imply
\[
 |\nabla_X^k h(X)| \lesssim_k m \overline\rho^3|X|^{-5-k}, \qquad \text{for} \quad |X| \geq 2\overline\rho,\ k = 0,1,2.
\]
The velocity and gradient bounds \eqref{eq:streamvelocitybounds}
imply \eqref{eq:innerfar}.

Using the cutoff \eqref{eq:radialcutoffs}, define
$\widetilde h = \chi^{\mathrm{in}}_{\underline\rho}h$, continued by zero
through the inner ball. On the cutoff annulus $|X|$ is
comparable to $\underline\rho$; the product rule, the preceding bounds,
and \eqref{eq:streamvelocitybounds} yield \eqref{eq:innerextension} for
$\mathfrak{e}^{\mathrm{in}} = \mathcal V[\widetilde h]$.
This $C^1$ divergence-free field has the prescribed parity and
equals $u$ outside $B_{\underline\rho/2}$.
\end{proof}

The two main estimates of this section are closely related to \cite[Lemma 2.1]{KJ22} and \cite[Proposition 2.1]{BC26}. One of the main differences is that they treated self interaction as error whereas we don't.

\subsection{The combined error field}

\begin{setup}[A separated family]\label{setup:family}
On a time interval $I$, consider vortex-ring pairs
$(q_j)_{j=1}^N$, where $N\in\mathbb N\cup\{\infty\}$, each satisfying
Setup~\ref{setup:ringpair} (a fixed-time pair) at every time. They are transported by one
Euler velocity:
\[
 \partial_tq_j+u\cdot\nabla q_j=0,\qquad
 u_j=\BS[\vectomega_j],\qquad u=\sum_j u_j.
\]
Assume the velocity series converges in $C^1(\mathbb{R}^3)$ at each fixed time $t\in I$. Here $u$ denotes the
total velocity, while $u_j$ denotes the self-induced velocity of
vortex-ring pair $j$. Write $m_j=\int_{\Hh}rq_j\dd r\dd z$
for its conserved mass, and $\rho_j^+,\rho_j^-$ for its maximal and
minimal support distances from the origin. Let
$0<\underline\rho_j(t)\leq\overline\rho_j(t)<\infty$ satisfy, for every $t\in I$,
\[
 \rho_j^+\leq\overline\rho_j,\qquad
 \rho_j^-\geq\underline\rho_j,\qquad
 16\overline\rho_j\leq\underline\rho_i\quad(i<j).
\]
Let $a_0=0$ and $a_j=\sum_{i\leq j}b_i$. For $N=\infty$, assume that
$\sum_{i>j}m_i\overline\rho_i^3$ converges locally uniformly in time
for each fixed $j$.
\end{setup}

\begin{lemma}[Combined error field]\label{lem:combinederror}
Under Setup~\ref{setup:family} (a separated family), the following hold for each $j$.
\begin{enumerate}
\item For every $t\in I$, there is a global divergence-free axisymmetric field
$\mathfrak{e}_j(t,\cdot)\in C^1(\R^3)$, without swirl and with the prescribed parity,
agreeing with $\sum_{i\ne j}u_i(t,x)-a_{j-1}(t)Dx$ on a neighborhood of $\supp q_j(t)$, such that
\begin{equation}\label{eq:combinederror}
 \begin{gathered}
 \|\nabla \mathfrak{e}_j(t,\cdot)\|_{L^\infty(\R^3)}\leq\ell_j(t),\\
 \ell_j(t)=C\left[\overline\rho_j(t)^2\sum_{i=1}^{j-1}
 \frac{b_i(t)}{\underline\rho_i(t)^2}
 +\underline\rho_j(t)^{-5}\sum_{i>j}m_i\overline\rho_i(t)^3\right].
 \end{gathered}
\end{equation}
The constant $C$ depends only on the fixed cutoffs $\chi^{\mathrm{out}}, \chi^{\mathrm{in}}$ and is independent
of time, the indices, and the support radii.
\item If $u_i\in C(I;C^1(\R^3))$ for every $i$ and
$\underline\rho_j,\overline\rho_j\in C(I)$, then
$\mathfrak{e}_j\in C(I;C^1(\R^3))$.
\end{enumerate}
\end{lemma}

\begin{proof}
\textup{(i)} Fix $j$ and $t\in I$, and suppress the time variable.
For each source pair $i$, let $h_i$ be defined by
\eqref{eq:representation} with $q=q_i$. Define
\[
 \mathfrak{e}^{\mathrm{out}}_{i\to j}
 :=\mathcal V\bigl[\chi^{\mathrm{out}}_{\overline\rho_j}(h_i-b_i)\bigr],
 \quad (i<j),
 \qquad\mathfrak{e}^{\mathrm{in}}_{i\to j}
 :=\mathcal V\bigl[\chi^{\mathrm{in}}_{\underline\rho_j}h_i\bigr],
 \quad (i>j).
\]
These are the extensions from source pair $i$ to receiving pair $j$
constructed in Lemmas~\ref{lem:outer} and~\ref{lem:inner}, with the
same zero continuations. Set
\[
 \mathfrak{e}_j
 :=\sum_{i<j}\mathfrak{e}^{\mathrm{out}}_{i\to j}
   +\sum_{i>j}\mathfrak{e}^{\mathrm{in}}_{i\to j}.
\]
The outer sum is finite. By \eqref{eq:innerextension},
\[
 \|\mathfrak{e}^{\mathrm{in}}_{i\to j}\|_{C^1(\R^3)}
 \lesssim\bigl(\underline\rho_j^{-4}+\underline\rho_j^{-5}\bigr)
 m_i\overline\rho_i^3,
 \qquad \sum_{i>j}m_i\overline\rho_i^3<\infty.
\]
Thus the inner sum converges in $C^1(\R^3)$, and $\mathfrak{e}_j$
is divergence-free, axisymmetric, without swirl, and has the prescribed parity.
On the neighborhood
$\{\underline\rho_j/2<|x|<3\overline\rho_j/2\}$ of $\supp q_j$, we have
\[
 \mathfrak{e}_j(x)
 =\sum_{i<j}\bigl(u_i(x)-b_iDx\bigr)+\sum_{i>j}u_i(x)
 =\sum_{i\ne j}u_i(x)-a_{j-1}Dx.
\]
Summing the gradient bounds from the two lemmas yields
\[
 \|\nabla\mathfrak{e}_j\|_{L^\infty(\R^3)}
 \leq C\left[\overline\rho_j^2\sum_{i<j}\frac{b_i}{\underline\rho_i^2}
 +\underline\rho_j^{-5}\sum_{i>j}m_i\overline\rho_i^3\right]=\ell_j,
\]
with $C$ the maximum of the two extension constants.

\textup{(ii)} Fix a compact subinterval $I_0\subseteq I$. By positivity
and continuity of the radii,
\[
 0<\inf_{t\in I_0}\underline\rho_j(t)
 \leq\sup_{t\in I_0}\overline\rho_j(t)<\infty,
 \qquad\text{and}\qquad
 \chi^{\mathrm{out}}_{\overline\rho_j(\cdot)},\, 
 \chi^{\mathrm{in}}_{\underline\rho_j(\cdot)}
 \in C(I_0;C^2(\R^2)).
\]

Let $v=\mathcal V[h]$. Its axis-normalized stream function satisfies
\[
 \psi_v(r,z)=-r^2zh(r,z)=\int_0^r s\,v^z(s,z)\dd s.
\]
For a radial cutoff $\chi$, we have
$\mathcal V[\chi h]=r^{-1}\nabla^\perp(\chi\psi_v)$, so the product rule yields
\[
 \mathcal V[\chi h]
 =\chi v+r\left(\int_0^1 s\,v^z(sr,z)\dd s\right)
 (-\chi_z\mathbf e_r+\chi_r\mathbf e_3).
\]
In Cartesian coordinates, the integral factor is
\[
 I_v(x):=\int_0^1 s\,v_3(sx_1,sx_2,x_3)\dd s,
 \qquad x=(x_1,x_2,x_3),
\]
where $v_3$ is the axial component of $v$.
Differentiation under the integral yields, for every $R>0$,
\[
 \|I_v-I_w\|_{C^1(B_R)}\leq C\|v-w\|_{C^1(B_R)}.
\]
Hence $v\in C(I_0;C^1(B_R))$ implies
$I_v\in C(I_0;C^1(B_R))$.
For $r>0$, the correction term equals
\[
 rI_v(-\chi_z\mathbf e_r+\chi_r\mathbf e_3)
 =I_v\begin{pmatrix}
 -x_1\partial_3\chi\\
 -x_2\partial_3\chi\\
 x_1\partial_1\chi+x_2\partial_2\chi
 \end{pmatrix}.
\]
Since $I_v$ is $C^1$ and the radial cutoff $\chi$ is $C^2$ in
Cartesian coordinates, the right-hand side defines a $C^1$ field
also on the axis.
Consequently, we can conclude that on every fixed ball $B_R$,
if
 $\mathcal V[h]\in C(I_0;C^1(B_R))$ and $
 \chi\in C(I_0;C^2(B_R))$
then  $\mathcal V[\chi h]\in C(I_0;C^1(B_R))$.

Apply this conclusion with $h=h_i$ and
$\chi=\chi^{\mathrm{in}}_{\underline\rho_j}$ for the inner extensions,
and with $h=h_i-b_i$ and
$\chi=\chi^{\mathrm{out}}_{\overline\rho_j}$ for the outer extensions.
Here $\mathcal V[h_i]=u_i$, $\mathcal V[h_i-b_i]=u_i-b_iDx$, and
$b_i=\partial_r u_i^r(\cdot,0,0)\in C(I_0)$.
Therefore,
\[
 \mathfrak{e}^{\mathrm{out}}_{i\to j}\in C(I_0;C^1(B_R))\quad(i<j),
 \qquad
 \mathfrak{e}^{\mathrm{in}}_{i\to j}\in C(I_0;C^1(B_R)),\quad (i>j)
\]
for every fixed $R>0$. 

To obtain continuity in $C^1(\R^3)$, choose
$R>2\sup_{t\in I_0}\overline\rho_j(t)$, independently of time.
The cutoff definitions imply, for every $t\in I_0$ and $|x|\geq R$,
\[
 \mathfrak{e}^{\mathrm{out}}_{i\to j}(t,x)=0\quad(i<j),
 \qquad
 \mathfrak{e}^{\mathrm{in}}_{i\to j}(t,x)=u_i(t,x)\quad(i>j).
\]
Since $u_i\in C(I_0;C^1(\R^3))$ by hypothesis, these identities
and the continuity on $B_R$ yield
\[
 \mathfrak{e}^{\mathrm{out}}_{i\to j}\in C(I_0;C^1(\R^3))\quad(i<j),
 \qquad
 \mathfrak{e}^{\mathrm{in}}_{i\to j}\in C(I_0;C^1(\R^3)),\quad (i>j).
\]
For $N=\infty$ and $M\geq j$, the inner-extension bound and local
uniform summability imply
\[
 \sup_{t\in I_0}
 \left\|\sum_{i>M}\mathfrak{e}^{\mathrm{in}}_{i\to j}(t)\right\|_{C^1(\R^3)}
 \leq C_{I_0,j}\sup_{t\in I_0}\sum_{i>M}m_i\overline\rho_i(t)^3
 \to0\qquad\text{as }M\to\infty.
\]
Adding the finite outer sum, we conclude
$\mathfrak{e}_j\in C(I_0;C^1(\R^3))$. Since $I_0$ is arbitrary,
$\mathfrak{e}_j\in C(I;C^1(\R^3))$.
\end{proof}
\endgroup

\section{{The multiscale construction}}
\label{sec:cascade}

In this section we construct an initial family of vortex-ring pairs
whose Euler evolution cannot remain classical on
$[0,T_{\mathrm{cmp}}]$, where
$T_{\mathrm{cmp}}=\tau_*/\lambda_{1,0}$. We transport each initial
component by the common Euler flow and apply the single-pair
estimates of Sections~\ref{sec:ringpairs}--\ref{sec:geometry} using
the error fields from Section~\ref{sec:interactions}.

We first derive conditions on the scales that keep the supports separated and the interaction errors small.
Later we choose scales that satisfy all the required conditions. On a classical
interval, the support estimates rule out a first failure of the
support bounds. If a pair were the first to reach the
angular-redistribution threshold $\tau_*$, then by Theorem~\ref{thm:pulse}, that pair would have to accumulate more strain than the upper bound allows.
The joint bootstrap thus keeps every angular-redistribution time below $\tau_*$.
Then by classical continuation to $T_{\mathrm{cmp}}$ we would have
$\tau_1(T_{\mathrm{cmp}})=\tau_*$, a contradiction.

\Needspace{18\baselineskip}
The following table records the constants used repeatedly in this
section and the estimates or choices that introduce them.

\smallskip
\begingroup
\small
\setlength{\tabcolsep}{4pt}
\renewcommand{\arraystretch}{1.18}
\noindent
\begin{tabularx}{\linewidth}{@{}
 >{\raggedright\arraybackslash}p{.23\linewidth}
 >{\raggedright\arraybackslash}X
 >{\raggedright\arraybackslash}p{.31\linewidth}@{}}
\toprule
Constants & Role & Reference \\
\midrule
$c_f,C_0,C_a$
 & Capacity, strain, and aspect bounds
 & Setup~\ref{setup:parameters} \\
\addlinespace[3pt]
$C_{\mathrm{upp}},C_{\mathrm{low}}$
 & Upper and lower support bounds
 & \eqref{eq:comparisonupperbudget} and \eqref{eq:comparisonlowerbudget} \\
\addlinespace[3pt]
$C_{\mathrm{int}},\widetilde C_{\mathrm{int}}$
 & Interaction and integrated error bounds
 & \eqref{eq:combinederror} and \eqref{eq:integratederror} \\
\addlinespace[3pt]
$K_f,c_d$
 & Comparison prefactors
 & \eqref{eq:comparisonprefactors}; proof of Lemma~\ref{lem:comparisonchoice} \\
\addlinespace[3pt]
$C_{\mathrm{sep}}$
 & Minimum scale growth
 & \conditionref{cond:geometricgrowth} \\
\addlinespace[3pt]
$C_{\mathrm{geom}}$
 & Lower threshold on $\widehat L_1$ for the comparison construction
 & Lemma~\ref{lem:comparisonchoice} \\
\bottomrule
\end{tabularx}
\par
\endgroup
\medskip

\subsection{Initial data}

Fix the profile $f_0$ and its support box from
Setup~\ref{setup:profile} {(an evolving, scaled pair)}. Let $\alpha_j$ be an increasing sequence converging to $1/3$, with
$1/6\leq\alpha_j<1/3$, and define the regularity gap by
\begin{equation}\label{eq:regularityparameters}
 \beta_j:=1-3\alpha_j,\qquad 0<\beta_j\leq\frac12.
\end{equation}
We define the radial scale and vorticity amplitude as
\begin{equation}\label{eq:radialparameters}
 L_j:=\log\frac1{R_{j,0}}>0,\qquad
 R_{j,0}=e^{-L_j},\qquad Q_{*,j}=R_{j,0}^{\alpha_j-1}.
\end{equation}
The letter $L$ is chosen to represent Logarithmic scales. In addition, the minus sign in front of $L_j$ in the definition of $R_{j,0}$ is chosen so that increasing $L_j$ is generally equivalent to increasing relative separation between ring pairs of different scales.
With this amplitude choice, by \eqref{eq:holder} in Lemma~\ref{lem:initialregularity} we will have
\[
 \|\vectomega_{j,0}\|_{C^\alpha(\R^3)}
 \lesssim_{f_0,\alpha}e^{-(\alpha_j-\alpha)L_j}.
\]
For each fixed $0<\alpha<1/3$, we have
$\alpha_j-\alpha\geq\frac12(1/3-\alpha)>0$ for all sufficiently large $j$.

For pair $j$, the capacity defined in
\eqref{eq:singlecapacity} satisfies
\begin{equation}\label{eq:capacityparameter}
 P_j:=Q_{*,j}(R_{j,0}^2H_{j,0})^{1/3},
 \qquad
 P_j^3=\frac{H_{j,0}}{R_{j,0}^{\beta_j}}.
\end{equation}
To make the heights shrink and the capacities grow, we need
$R_{j,0}^{\beta_j}\ll H_{j,0}\ll1$. We therefore choose
\begin{equation}\label{eq:heightparameters}
 \widehat L_j:=\beta_jL_j,\qquad
 H_{j,0}=R_{j,0}^{\kappa\beta_j}=e^{-\kappa\widehat L_j},
 \qquad 0<\kappa<1.
\end{equation}
Then $H_{j,0}\to0$ and $P_j=e^{(1-\kappa)\widehat L_j/3}\to\infty$
as $\widehat L_j\to\infty$. Note that it is essential that $\alpha_j <1/3$.
To close the support estimates, we will also impose
$\kappa>1/3$; see \eqref{eq:uppervelocity}.

\begin{remark}[Motivation behind the relative scale between $H_{j,0}$ and $R_{j,0}$]
\label{rem:relative_scale}
If instead of \eqref{eq:heightparameters}, we choose $H_{j,0}$ to have the same scale as $R_{j,0}$, that is, $H_{j,0}\sim R_{j,0}$, as in \cite{KJ22,BC26}, then $P_j\sim R_{j,0}^{\alpha_j}$. Thus in a natural scale choice for $R_{j,0}$, the capacity $P_j$ decays exponentially fast, leading to global regularity of the solution. Therefore, it is necessary to choose $H_{j,0}$ to have a bigger scale than $R_{j,0}$ with $\kappa \beta_j <1$ as in \eqref{eq:heightparameters}. 

For this reason, we can connect this choice of relative scale to a large cone slope, the key idea of \cite{BC26}. In that paper, we need to choose a large fixed cone slope $1/\varepsilon_{j,0}=H_{j,0}/R_{j,0}$ to achieve a sharp instantaneous blow up whereas in this paper, we need to choose cone slopes tending to infinity as we approach the origin in order to establish finite-time blow up.
\end{remark}

We distinguish the requirements used in the evolution argument,
numbered \ref{cond:initialsummability}--\ref{cond:errorcondition}, from
sufficient restrictions on the initial parameters, labeled
\ref{cond:nextscaleupper}--\ref{cond:initialthreshold}, which ensure the numerically numbered conditions. See Lemma~\ref{lem:parameterconditions}.

\begin{condition}[Initial decay and summability]\label{cond:initialsummability}
We require
\[
 \widehat L_j\to\infty,\qquad
 \sum_{j\geq1}e^{-\eta L_j}<\infty\quad\text{for every }\eta>0.
\]
\end{condition}
These conditions ensure that the initial supports shrink to the origin
and that the vorticity and velocity series have the required regularity.

\begin{setup}[{An initial family}]\label{setup:parameters}
Fix $f_0$, the cutoffs $\Phi, \chi^{\mathrm{in}}, \chi^{\mathrm{out}}$ of \eqref{eq:angularcutoff} and
\eqref{eq:radialcutoffs}, and $0<\kappa<1$. With the parameters in
\eqref{eq:regularityparameters}--\eqref{eq:heightparameters}, assume
\begin{itemize}
\item \conditionref{cond:initialsummability}.
\end{itemize}
On $z>0$ define
\begin{equation}\label{eq:initialdata}
 \begin{gathered}
 q_{j,0}(r,z)=Q_{*,j}f_0(r/R_{j,0},z/H_{j,0}),
 \qquad 
 \vectomega_{j,0}(x)=q_{j,0}(r,z)(x_2,-x_1,0),
 \end{gathered}
\end{equation}
with odd reflection of $q_{j,0}$ in $z$. Let
$\vectomega_0=\sum_j\vectomega_{j,0}$ and
$u_0=\sum_j\BS[\vectomega_{j,0}]$.
Since all pairs use the same initial profile $f_0$,
$M_{z,j}(0)=M_{z,0}$ for every $j$.
{Retain the constants $F,E_*,M_{z,0},C_{f_0,\Phi},\tau_*,c_*,\eps_*$ of
Sections~\ref{sec:ringpairs} and~\ref{sec:pulse}, and take}
\[
 c_f=F^{2/3},\qquad C_0=8c_f\tau_*,\qquad C_a=1+C_0.
\]
\end{setup}

Since $\kappa\beta_j<1$, we have $R_{j,0}\leq H_{j,0}<1$.
By Lemma~\ref{lem:initialregularity} and \conditionref{cond:initialsummability},
the vorticity and velocity series converge absolutely in $C^\alpha(\R^3)$
and $C^{1,\alpha}(\R^3)\cap L^2(\R^3)$, respectively, for every
$0<\alpha<1/3$. Classical local existence then applies.
Constants in the conditional estimates depend only on the fixed
profile $f_0$, cutoffs $\Phi, \chi^{\mathrm{in}}, \chi^{\mathrm{out}}$, and $\kappa$, unless otherwise indicated.

\subsection{Transport of the vortex-ring pairs}

\begin{setup}[Classical multiscale evolution]\label{setup:cascade}
{Under} Setup~\ref{setup:parameters} {(an initial family)}, let $u$ be a
finite-energy classical Euler solution with initial data $u_0$ on
$[0,t_0]$, with
$\sup_{0\leq t\leq t_0}\|u(t)\|_{C^{1,\alpha}}<\infty$
for some $0<\alpha<1/3$. Let $q_j$ solve
\[
 (\partial_t+u\cdot\nabla)q_j=0,\qquad q_j|_{t=0}=q_{j,0},
\]
and write
\[
 \begin{gathered}
 \omega_j=-rq_j,\qquad
 \vectomega_j=q_j(x_2,-x_1,0),\qquad u_j=\BS[\vectomega_j],\\
 A_{\mathrm{Lip}}(t_0):=\int_0^{t_0}
 \|\nabla u(t)\|_{L^\infty(\R^3)}\dd t<\infty.
 \end{gathered}
\]
\end{setup}

For each fixed $j$, the set
$\bigcup_{0\leq t\leq t_0}\supp q_j(t)$ is compact and disjoint from
$\{r=0\}\cup\{z=0\}$. By classical transport and elliptic regularity we have
\[
 \begin{gathered}
 q_j(t)\in C_c^{1,\alpha}(\R^3),\qquad
 u_j(t)\in C^{2,\alpha}(\R^3), \qquad \forall t\in[0,t_0],\\
 u_j\in C([0,t_0];L^2(\R^3)\cap C^1(\R^3)).
 \end{gathered}
\]
The single-pair identities follow by integration by parts and smooth
approximation.

\begin{lemma}[Decomposition of the velocity]\label{lem:materialdecomposition}
Under Setup~\ref{setup:cascade} {(classical multiscale evolution)},
$u=\sum_{j\geq1}u_j$ in $C([0,t_0];L^2(\R^3)\cap C^1(\R^3))$.
Moreover,
\begin{gather}
 \sup_{0\leq t\leq t_0}\|u_j(t)\|_{L^2(\R^3)}
 \leq e^{A_{\mathrm{Lip}}(t_0)}\|u_j(0)\|_{L^2(\R^3)},\label{eq:taggedenergy}\\
 \sup_{0\leq t\leq t_0}\|\vectomega_j(t)\|_{C^\alpha(\R^3)}
 \lesssim_{u,t_0,\alpha}\|\vectomega_{j,0}\|_{C^\alpha(\R^3)}.
 \label{eq:taggedpersistence}
\end{gather}
\end{lemma}
\begin{proof}
Each component satisfies
\[
 (\partial_t+u\cdot\nabla)\vectomega_j = \vectomega_j\cdot\nabla u,\qquad
 \vectomega_j|_{t=0}=\vectomega_{j,0}.
\]
The H\"older estimate for this linear equation implies
\eqref{eq:taggedpersistence}. Since
$\partial_tu_j=\mathbb P(u\times\vectomega_j)$,
\[
 \frac12\frac{\dd}{\dd t}\|u_j\|_{L^2(\R^3)}^2
 \leq\|\nabla u\|_{L^\infty(\R^3)}\|u_j\|_{L^2(\R^3)}^2.
\]
Gronwall's inequality yields \eqref{eq:taggedenergy}.
By \eqref{eq:schauder} and Lemma~\ref{lem:initialregularity},
\[
 \begin{aligned}
 &\sum_{j\geq1}\sup_{0\leq t\leq t_0}
 \left(\|u_j(t)\|_{L^2}+\|u_j(t)\|_{C^1}\right)\\
 &\qquad\lesssim_{u,t_0,\alpha}
 \sum_{j\geq1}\left(\|u_j(0)\|_{L^2}
 +\|\vectomega_{j,0}\|_{C^\alpha}\right)<\infty.
 \end{aligned}
\]
Thus the series converges uniformly on $[0,t_0]$ in $L^2(\R^3)$ and $C^1(\R^3)$.
By uniqueness for the linear vorticity equation, we have
$\sum \vectomega_j=\curl u$.
Since both $u$ and $\sum u_j$ are divergence free and belong
to $L^2(\R^3)$, they coincide.
\end{proof}

Take $a_0=A_0=0$ and
\begin{equation}\label{eq:cumulative}
 \begin{aligned}
 b_j(t)&:=\partial_ru_j^r(t,0,0),
 & B_j(t)&:=\int_0^tb_j(s)\dd s,\\
 a_j(t)&:=\sum_{i=1}^jb_i(t),
 & A_j(t)&:=\sum_{i=1}^jB_i(t)=\int_0^ta_j(s)\dd s.
 \end{aligned}
\end{equation}
{The strain rate at the origin generated by pair $j$ is $b_j$,
whereas the outer pairs contribute $a_{j-1}$. Their strains are $B_j$
and $A_{j-1}$, respectively.} By \eqref{eq:mass}, the following quantities are conserved:
\[
 \begin{gathered}
 Q_j=\|q_j(t)\|_{L^\infty}=FQ_{*,j},\\
 m_j=\int_{\Hh}rq_j(t,r,z)\dd r\dd z
 =Q_{*,j}R_{j,0}^2H_{j,0}
 =e^{-(1+\alpha_j)L_j-\kappa\widehat L_j},\\
 Q_j^{2/3}m_j^{1/3}=c_fP_j.
 \end{gathered}
\]

\begin{lemma}[Strain and capacity bounds]\label{lem:cascadestrain}
Under Setup~\ref{setup:cascade} {(classical multiscale evolution)}, for $0\leq t\leq t_0$,
\begin{equation}\label{eq:totalsourcebound}
 0\leq\sum_{j\geq1}b_j(t)=\partial_ru^r(t,0,0)
 \leq\|\nabla u(t)\|_{L^\infty(\R^3)},\qquad
 0\leq A_j(t)\leq A_{\mathrm{Lip}}(t_0).
\end{equation}
For every $j\geq1$,
\begin{equation}\label{eq:tagcapacity}
 0\leq b_j(t)\leq c_fP_j.
\end{equation}
\end{lemma}
\begin{proof}
Lemma~\ref{lem:materialdecomposition} implies
$\sum_{j\geq1}b_j(t)=\partial_ru^r(t,0,0)$.
By Lemma~\ref{lem:capacity},
$0\leq b_j(t)\leq Q_j^{2/3}m_j^{1/3}=c_fP_j$.
Integrating $0\leq a_j(t)\leq\|\nabla u(t)\|_{L^\infty(\R^3)}$, we
obtain $0\leq A_j(t)\leq A_{\mathrm{Lip}}(t_0)$.
\end{proof}

\subsection{{Bounds on angular-redistribution time}}

For each $j\geq1$, define 
\begin{equation}\label{eq:cascadeclocks}
 \begin{gathered}
 R_j:=R_{j,0}e^{A_{j-1}},\qquad H_j:=H_{j,0}e^{-2A_{j-1}},\qquad
 \eps_j:=\frac{R_j}{H_j} = \eps_{j,0}e^{3A_{j-1}},\\
 \lambda_j:=\frac{m_j}{H_j^2} = \lambda_{j,0}e^{4A_{j-1}}=P_j\eps_j^{4/3},\qquad
 \tau_j(t):=\int_0^t\lambda_j(s)\dd s,\\
 \lambda_{j,0}=e^{-(1+\alpha_j)L_j+\kappa\widehat L_j},\qquad
 \eps_{j,0}=e^{-L_j+\kappa\widehat L_j}.
 \end{gathered}
\end{equation}
Here, the time variable is suppressed for simplicity.
And $R_j,H_j$ are the reference scales for the affine field
$a_{j-1}Dx$. {The variable $\tau_j$ is the angular-redistribution time
of pair $j$, and $\lambda_j=\tau_j'$ is its rate with respect to physical
time.}

\begin{condition}[Capacity growth]\label{cond:capacitygap}
For every $j\geq2$, we require 
\[
 P_j\geq2P_{j-1},
 \qquad j\geq2,
\]
which is equivalent to $(1-\kappa)(\widehat L_j-\widehat L_{j-1})\geq3\log 2$ for $j\geq 2$.
\end{condition}

\begin{lemma}[{Conditional bounds on $\tau_j$, strain, and aspect}]\label{lem:clockbounds}
Consider Setup~\ref{setup:cascade} {(classical multiscale evolution)} and assume
\begin{itemize}
    \item \conditionref{cond:capacitygap}.
\end{itemize} Then the
following estimates hold for every $t\in[0,t_0]$:
\begin{equation}\label{eq:clockcomparison}
 \lambda_j(t)\leq\lambda_{j,0}+8c_fP_{j-1}\tau_j(t), \qquad j\geq 2.
\end{equation}
Consequently, for $j\geq1$, if $\tau_{j+1}(t)\leq\tau_*$, then
\begin{equation}\label{eq:clockstraincap}
 A_j(t)\leq\frac14\log\left(1+
 C_0\frac{P_j}{\lambda_{j+1,0}}\right)=:\overline A_j.
\end{equation}
For $j\geq2$, if $\tau_j(t)\leq\tau_*$, then
\begin{equation}\label{eq:aspectuniform}
 \eps_j(t)\leq\left(\eps_{j,0}^{4/3}
 +C_0\frac{P_{j-1}}{P_j}\right)^{3/4}=:\bar\eps_j.
\end{equation}
For $j=1$,
\begin{equation}\label{eq:firstclock}
 \eps_1(t)=\eps_{1,0}=: \bar\eps_1,\qquad
 \tau_1(t)=\lambda_{1,0}t.
\end{equation}
\end{lemma}
We also take $\overline A_0=0$.
\begin{proof}
By \eqref{eq:tagcapacity} and \conditionref{cond:capacitygap}, for $j\geq2$,
\[
    a_{j-1}(t)=\sum_{i<j}b_i(t)
 \leq c_f\sum_{i<j}P_i\leq2c_fP_{j-1}.
\]
Apply Lemma~\ref{lem:singleclockbounds} to the reference
scales of pair $j$, with $a=a_{j-1}$, $\bar a=2c_fP_{j-1}$,
$P=P_j$, and $\lambda_0=\lambda_{j,0}$.
Since $C_0=8c_f\tau_*$, this proves
\eqref{eq:clockcomparison} and \eqref{eq:aspectuniform} for $j\geq2$.
Applying the same lemma to pair $j+1$ yields
\eqref{eq:clockstraincap} for every $j\geq1$.

For $j=1$, $a_0=A_0=0$ implies
$\tau_1(t)=\lambda_{1,0}t$ and $\eps_1(t)=\eps_{1,0}=\bar\eps_1$, hence
\eqref{eq:firstclock} holds.
\end{proof}

\subsection{{Application of the large strain estimate}}

{
\begin{condition}[Small reference aspect]\label{cond:smallreferenceaspect}
For every $j\geq1$, we require
\[
 \bar\eps_j\leq\eps_*.
\]
\end{condition}
 
This is the reference-aspect condition in Theorem~\ref{thm:pulse}.
Since $\eps_{j,0}\leq\bar\eps_j$, it also ensures
$M_{\Phi,j}(0)=1$ by \eqref{eq:activationaspect}.

\begin{condition}[Strain comparison]\label{cond:straincomparison}
For every $j\geq1$, we require
\[
 \overline A_j<c_*\bar\eps_j^{-4/7}.
\]
\end{condition}

This condition makes the strain produced when a pair reaches
the threshold exceed its allowed contribution to $A_j$.

\begin{lemma}[Large strain at the prescribed threshold]\label{lem:cascadepulse}
{Under} Setup~\ref{setup:cascade} {(classical multiscale evolution)}, 
{impose}
\begin{itemize}
\item
\conditionref{cond:capacitygap}, 
\item \conditionref{cond:smallreferenceaspect}, and 
\item 
\conditionref{cond:straincomparison}.
\end{itemize} Fix $j\geq1$. 
Suppose there exists a time
$t_{*,j}\in(0,t_0]$ such that $\tau_j(t_{*,j})=\tau_*$, and
assume that pair $j$ admits an error field as in
Setup~\ref{setup:evolution} {(an evolving pair)} on $[0,t_{*,j}]$, with outer strain
rate $a_{j-1}$ and $\mathfrak{E}_j\leq1$. Then
\begin{equation}\label{eq:cascadepulse}
 B_j(t_{*,j})\geq c_*\bar\eps_j^{-4/7}>\overline A_j.
\end{equation}
\end{lemma}

\begin{proof}
Since $\tau_j$ is increasing, $\tau_j(t)\leq\tau_*$ on
$[0,t_{*,j}]$. Lemma~\ref{lem:clockbounds} therefore implies
$\eps_j(t)\leq\bar\eps_j\leq\eps_*$ throughout this interval.
Apply Theorem~\ref{thm:pulse} with the initial profile $f_0$,
scales $Q_{*,j},R_{j,0},H_{j,0}$, and affine strain rate
$a_{j-1}$. The bound on $\mathfrak{E}_j$ verifies the remaining hypothesis,
so \eqref{eq:pulse}, \conditionref{cond:smallreferenceaspect}, and \conditionref{cond:straincomparison} yield
\eqref{eq:cascadepulse}.
\end{proof}
}

\begin{parametercondition}[Upper bound on scale growth]\label{cond:nextscaleupper}
For every $j\geq1$, we require
\[
 L_{j+1}\leq\widehat L_j^2=(\beta_jL_j)^2,\qquad j\geq1.
\]
\end{parametercondition}

\begin{lemma}[Sufficient conditions for the aspect and strain bounds]
\label{lem:aspectstrainconditions}
With the parameters in
\eqref{eq:regularityparameters}--\eqref{eq:heightparameters}, fix
$0<\kappa<1$. Assume 
\begin{itemize}
    \item \conditionref{cond:nextscaleupper} and
    \item
$\widehat L_j\geq2\widehat L_{j-1}$ for $j\geq2.$
\end{itemize}
There exists a constant $S_{\mathrm{asp}}>0$, depending only on the profile $f_0$,
cutoffs $\Phi, \chi^{\mathrm{in}}, \chi^{\mathrm{out}}$, and $\kappa$, such that $\widehat L_1\geq S_{\mathrm{asp}}$
ensures 
\begin{itemize}
    \item 
\conditionref{cond:smallreferenceaspect} and
    \item
\conditionref{cond:straincomparison} for every $j\geq1$.
\end{itemize}
\end{lemma}

\begin{proof}
By the definitions in
\eqref{eq:clockstraincap}--\eqref{eq:aspectuniform},
\begin{align}
 \overline A_j&=\frac14\log\left[1+C_0\exp\left(
 \frac{1-\kappa}{3}\widehat L_j
 +(1+\alpha_{j+1})L_{j+1}-\kappa\widehat L_{j+1}\right)\right],
 \label{eq:explicitstrain}\\
 \bar\eps_j^{4/3}&=e^{-\frac43(L_j-\kappa\widehat L_j)}
 +C_0e^{-\frac{1-\kappa}{3}(\widehat L_j-\widehat L_{j-1})},
 \qquad j\geq2.\label{eq:explicitaspect}
\end{align}
Since $L_j\geq2\widehat L_j$, \eqref{eq:explicitaspect} implies, for $j\geq 2$,
\begin{equation}\label{eq:coarseaspect}
 \bar\eps_j^{4/3}
 \leq C_a e^{-(1-\kappa)\widehat L_j/6},\quad
 \text{and hence}\quad
 \bar\eps_j^{4/7}
 \leq C_a^{3/7}e^{-(1-\kappa)\widehat L_j/14}.
\end{equation}
The same estimate holds for $j=1$, since
$\bar\eps_1=e^{-L_1+\kappa\widehat L_1}$ and $L_1\geq2\widehat L_1$.
Moreover, \eqref{eq:explicitstrain}, the identity
$(1+\alpha_{j+1})L_{j+1}=(4L_{j+1}-\widehat L_{j+1})/3$,
and $\widehat L_j\leq\widehat L_{j+1}$ yield
\[
 \overline A_j\leq\frac{L_{j+1}}3+\frac14\log C_a.
\]
{
Using \conditionref{cond:nextscaleupper}, \eqref{eq:coarseaspect}, and the
preceding strain bound, we obtain
\[
 \begin{aligned}
 \bar\eps_j&\leq C_a^{3/4}
 e^{-(1-\kappa)\widehat L_j/8},\\
 \overline A_j\bar\eps_j^{4/7}
 &\leq C_a^{3/7}\left(\frac{\widehat L_j^2}{3}
                  +\frac14\log C_a\right)
 e^{-(1-\kappa)\widehat L_j/14}.
 \end{aligned}
\]
Both right-hand sides tend to zero as $\widehat L_j\to\infty$.
Choose $S_{\mathrm{asp}}$ so large that they are at most $\eps_*$ and
$c_*/2$, respectively, whenever $\widehat L_j\geq S_{\mathrm{asp}}$.
Since $\widehat L_j\geq\widehat L_1$, this establishes
\conditionref{cond:smallreferenceaspect} and
\conditionref{cond:straincomparison} for every $j\geq1$.
}
\end{proof}

{

\subsection{Interaction and support estimates}

Let
\[
 \rho_j^+(t):=\max_{\supp\vectomega_j(t)}|x|,\qquad
 \rho_j^-(t):=\min_{\supp\vectomega_j(t)}|x|.
\]
Fix $T\geq\max\{1,t_0\}$ and take
$V_{0,j}:=Q_{*,j}R_{j,0}^{6/5}H_{j,0}^{1/5}$.
If the $j$-th pair admits an error field as in
Setup~\ref{setup:evolution} {(an evolving pair)}, with outer strain rate $a_{j-1}$
and $\mathfrak{E}_j\leq1$, then, thanks to Propositions~\ref{prop:confinement}
and~\ref{prop:gap}, we have
\begin{align}
 e^{2A_{j-1}(t)}\rho_j^+(t)
 &\leq C_{\mathrm{upp}}\left[
 H_{j,0}+R_{j,0}e^{3A_{j-1}(t)}
 +T^{5/2}V_{0,j}^{5/2}e^{4A_{j-1}(t)}\right],
 \label{eq:comparisonupperbudget}\\
 e^{2A_j(t)}\rho_j^-(t)
 &\geq\frac{e^{-1}\rho_j^-(0)}
 {1+C_{\mathrm{low}}Q_j\rho_j^-(0)T}.
 \label{eq:comparisonlowerbudget}
\end{align}
Fix $C_{\mathrm{upp}},C_{\mathrm{low}}\geq1$, increasing
$C_{\mathrm{upp}}$ so that
$\rho_j^+(0)\leq C_{\mathrm{upp}}H_{j,0}$ for every $j$.
For positive constants $\overline d_j,\underline d_j$, depending
only on $T$ and the initial parameters, define
\begin{equation}\label{eq:comparisons}
 \overline\rho_j(t)=\overline d_j e^{-2A_{j-1}(t)},\qquad
 \underline\rho_j(t)=\underline d_j e^{-2A_j(t)}.
\end{equation}
Then
\begin{equation}\label{eq:comparisoncancellation}
 \frac{\overline\rho_{j+1}(t)}{\underline\rho_j(t)}
 =\frac{\overline d_{j+1}}{\underline d_j}.
\end{equation}
The exact cancellation of $e^{-2A_j(t)}$ in this ratio is the reason we retained $A_{\mathrm{tot}}$ in the lower bound \eqref{eq:gap} instead of estimating $B$ out of $A_{\mathrm{tot}}$ as explained right after the proof of \eqref{eq:gap}.

We call $\bar \rho_j, \underline \rho_j$ comparison radii and $\bar d_j, \underline d_j$ prefactor sequences. 

\begin{condition}[Support margins]\label{cond:supportmargins}
For every $j\geq1$, we require
\[
\begin{aligned}
 \overline d_j
 &\geq2C_{\mathrm{upp}}\left[
 H_{j,0}+R_{j,0}e^{3\overline A_{j-1}}
 +T^{5/2}V_{0,j}^{5/2}e^{4\overline A_{j-1}}\right],
 \\
 \underline d_j
 &\leq\frac{e^{-1}\rho_j^-(0)}
 {4\bigl(1+C_{\mathrm{low}}Q_j\rho_j^-(0)T\bigr)}.
\end{aligned}
\]
\end{condition}

Under \conditionref{cond:capacitygap} and the bounds
$\tau_i\leq\tau_*$ for all $i\geq1$, Lemma~\ref{lem:clockbounds}
implies $A_i\leq\overline A_i$. The support margins therefore yield
\[
\rho_j^+\leq\overline\rho_j/2,
\qquad
\rho_j^-\geq4\underline\rho_j
\]
 from the conditional estimates
\eqref{eq:comparisonupperbudget}--\eqref{eq:comparisonlowerbudget}.

\begin{condition}[Separation of comparison radii]\label{cond:prefactorseparation}
For every $j\geq1$, we require
\[
 0<\underline d_j\leq\overline d_j,\qquad
 \overline d_{j+1}\leq\frac1{32}\underline d_j,
 \qquad j\geq1.
\]
\end{condition}
By \eqref{eq:comparisoncancellation}, these conditions separate the
comparison radii throughout the evolution.

\begin{condition}[Summability of inner contributions]\label{cond:prefactorsummability}
We require
\[
 \sum_{i\geq1}m_i\overline d_i^3<\infty.
\]
\end{condition}
Let $C_{\mathrm{int}}$ be the constant in
\eqref{eq:combinederror}. Define
\begin{equation}\label{eq:interactionbudget}
 \mathcal I_j:=C_{\mathrm{int}}\left[
 \left(\frac{\overline d_j}{\underline d_{j-1}}\right)^2
 \overline A_{j-1}
 +\frac{\tau_*}{\lambda_{j+1,0}\underline d_j^5}
 \sum_{i>j}m_i\overline d_i^3\right],
\end{equation}
with the first term omitted for $j=1$.

\begin{condition}[Small interaction errors]\label{cond:errorcondition}
For every $j\geq1$, we require
\[
 \mathcal I_j\leq\frac12,\qquad j\geq1.
\]
\end{condition}
This condition controls the integrated interaction errors in
Lemma~\ref{lem:cascadeerrors} below. Define
\[
 \ell_j:=C_{\mathrm{int}}\left[
 \overline\rho_j^2\sum_{i<j}\frac{b_i}{\underline\rho_i^2}
 +\frac1{\underline\rho_j^5}\sum_{i>j}m_i\overline\rho_i^3
 \right],\qquad
 \mathfrak{E}_j(t):=\int_0^t\ell_j(s)\dd s.
\]

\Needspace{8\baselineskip}
\begin{lemma}[Conditional interaction and support improvement]
\label{lem:cascadeerrors}\label{lem:supportimprovement}
{Under} Setup~\ref{setup:cascade} {(classical multiscale evolution)}, 
{impose}
\begin{itemize}
\item 
\conditionref{cond:capacitygap} {holds and} $T\geq\max\{1,t_0\}$.
\end{itemize}
Let positive constants
$\overline d_j,\underline d_j$ satisfy
\begin{itemize}
\item
\conditionref{cond:supportmargins}, 
\item \conditionref{cond:prefactorseparation}, 
\item
\conditionref{cond:prefactorsummability},
and 
\item 
\conditionref{cond:errorcondition}. 
\end{itemize}
Suppose $0\leq t_1\leq t_0$ and
\begin{equation}\label{eq:weakclock}
 \tau_j(t)\leq\tau_*,
 \qquad j\geq1,\quad 0\leq t\leq t_1.
\end{equation}
Then the following hold on $[0,t_1]$.
\smallskip

\noindent\textup{(i)} \emph{Comparison estimates.}
Both sequences of comparison radii decrease with $j$, the series
defining each $\ell_j$ converges uniformly on $[0,t_1]$, and
\begin{equation}\label{eq:separation}
 \begin{aligned}
 \frac{\overline\rho_j(t)}{\underline\rho_i(t)}
 &\leq\frac{\overline d_j}{\underline d_{j-1}}\leq\frac1{32},
 &&1\leq i<j,\\
 \mathfrak{E}_j(t)&\leq\mathcal I_j\leq\frac12,
 &&j\geq1.
 \end{aligned}
\end{equation}
\smallskip

\noindent\textup{(ii)} \emph{Support improvement.}
If
\begin{equation}\label{eq:weaksupport}
 \rho_j^+(t)\leq\overline\rho_j(t),\qquad
 \rho_j^-(t)\geq\underline\rho_j(t),
 \qquad j\geq1,\quad 0\leq t\leq t_1,
\end{equation}
then the error fields of Lemma~\ref{lem:combinederror} exist
and satisfy
\[
 \mathfrak{e}_j\in C([0,t_1];C^1(\R^3)),\qquad
 \|\nabla \mathfrak{e}_j(t)\|_{L^\infty(\R^3)}\leq\ell_j(t).
\]
They agree with $\sum_{i\ne j}u_i-a_{j-1}Dx$ near
$\supp\vectomega_j(t)$. Moreover,
\begin{equation}\label{eq:supportimprovement}
 \rho_j^+(t)\leq\tfrac12\overline\rho_j(t),\qquad
 \rho_j^-(t)\geq4\underline\rho_j(t).
\end{equation}
The inequalities in \eqref{eq:weaksupport} are strict at time zero.
\end{lemma}

\begin{proof}
\textup{(i)} By \conditionref{cond:prefactorseparation}, both prefactor
sequences $\bar d_j, \underline d_j$ decrease. Since $A_j\geq A_{j-1}$, the comparison radii $\bar \rho_j, \underline \rho_j$
also decrease with $j$. For $i<j$, $A_{j-1}\geq A_i$, and hence
\[
 \frac{\overline\rho_j}{\underline\rho_i}
 =\frac{\overline d_j}{\underline d_i}e^{-2(A_{j-1}-A_i)}
 \leq\frac{\overline d_j}{\underline d_{j-1}}\leq\frac1{32}.
\]
For each fixed $j$,
\[
 \sum_{i>j}\sup_{0\leq t\leq t_1}
 \frac{m_i\overline\rho_i(t)^3}{\underline\rho_j(t)^5}
 \leq
 \frac{\sum_{i>j}m_i\overline d_i^3}
 {\min_{0\leq t\leq t_1}\underline\rho_j(t)^5}<\infty.
\]
Thus the series defining $\ell_j$ converges uniformly on $[0,t_1]$.

The outer contribution satisfies, for $j\geq2$,
\[
 \overline\rho_j^2\sum_{i<j}\frac{b_i}{\underline\rho_i^2}
 \leq\left(\frac{\overline d_j}{\underline d_{j-1}}\right)^2a_{j-1}.
\]
For the inner contribution, $A_{i-1}\geq A_j$ when $i>j$, so
\[
 \frac1{\underline\rho_j^5}\sum_{i>j}m_i\overline\rho_i^3
 =\frac1{\underline d_j^5}\sum_{i>j}
 m_i\overline d_i^3e^{10A_j-6A_{i-1}}
 \leq\frac{e^{4A_j}}{\underline d_j^5}
 \sum_{i>j}m_i\overline d_i^3.
\]
By \eqref{eq:cascadeclocks},
\begin{equation}\label{eq:clockintegral}
 \int_0^t e^{4A_j(s)}\dd s
 =\frac{\tau_{j+1}(t)}{\lambda_{j+1,0}}
 \leq\frac{\tau_*}{\lambda_{j+1,0}}.
\end{equation}
Integrating the preceding estimates and using
$A_{j-1}(t)\leq\overline A_{j-1}$ we obtain
$\mathfrak{E}_j(t)\leq\mathcal I_j\leq1/2$. Hence, \textup{(i)} holds.

\textup{(ii)} By \eqref{eq:weaksupport}, part~\textup{(i)}, and
Lemma~\ref{lem:materialdecomposition}, 
Setup~\ref{setup:family} {(a separated family)} {applies} on $[0,t_1]$.
Since each $u_i$ belongs to $C([0,t_1];C^1(\R^3))$ and the
comparison radii are positive and continuous,
Lemma~\ref{lem:combinederror} gives the asserted error fields.
Now $\mathfrak{E}_j\leq1/2$, so
\eqref{eq:comparisonupperbudget}--\eqref{eq:comparisonlowerbudget}
and \conditionref{cond:supportmargins}
yield
\[
 \begin{aligned}
 \rho_j^+(t)&\leq\tfrac12\overline d_j e^{-2A_{j-1}(t)}
 =\tfrac12\overline\rho_j(t),\\
 \rho_j^-(t)&\geq4\underline d_j e^{-2A_j(t)}
 =4\underline\rho_j(t).
 \end{aligned}
\]
Finally, \conditionref{cond:supportmargins} implies
\[
 \rho_j^+(0)\leq\tfrac12\overline d_j<\overline d_j,\qquad
 \rho_j^-(0)\geq4e\,\underline d_j>\underline d_j.
\]
Thus the initial inequalities in \eqref{eq:supportimprovement} are strict.
\end{proof}

}

{

\subsection{Choice of comparison radii}

In this subsection, we choose the constants $\overline d_j,\underline d_j$ in
\eqref{eq:comparisons}. For $j\geq2$,
the definition in \eqref{eq:clockstraincap}, $C_a=1+C_0$, and
$\lambda_{j,0}<1<P_{j-1}$ imply
\[
 e^{4\overline A_{j-1}}
 =1+C_0\frac{P_{j-1}}{\lambda_{j,0}}
 \leq C_a\frac{P_{j-1}}{\lambda_{j,0}}.
\]
Substituting the initial parameters, we obtain
\begin{align}
 R_{j,0}e^{3\overline A_{j-1}}
 &\leq C_a^{3/4}
 \exp\left(-(1+3\kappa)\widehat L_j/4+(1-\kappa)\widehat L_{j-1}/4\right),
 \label{eq:upperradial}\\
 V_{0,j}^{5/2}e^{4\overline A_{j-1}}
 &\leq C_a
 \exp\left(-(3\kappa-1)\widehat L_j/2+(1-\kappa)\widehat L_{j-1}/3\right).
 \label{eq:uppervelocity}
\end{align}
\begin{parametercondition}[Height exponent]\label{cond:heightmargin}
We require
\[
 \frac13<\kappa<1,\qquad
 \sigma:=\frac{3\kappa-1}{4}>0.
\]
\end{parametercondition}
The restriction $\kappa>1/3$ permits decay in \eqref{eq:uppervelocity}.
We use the positive constant $\sigma$ in choosing the comparison radii.

For the lower support bound, the initial support and
$(1-\alpha_j)-\kappa\beta_j\geq2\alpha_j\geq1/3$ yield
\[
 Q_j\rho_j^-(0)\geq F\zeta_-
 e^{(1-\alpha_j)L_j-\kappa\widehat L_j}
 \geq F\zeta_-e^{L_j/3}.
\]
If $Q_j\rho_j^-(0)\geq1$ and $T\geq1$, then
\[
 \frac{e^{-1}\rho_j^-(0)}
 {4\bigl(1+C_{\mathrm{low}}Q_j\rho_j^-(0)T\bigr)}
 \geq\frac{e^{-1}}{4(1+C_{\mathrm{low}})TQ_j}.
\]
We define
\begin{equation}\label{eq:comparisonprefactors}
 \overline d_j=K_fT^{5/2}e^{-\sigma\widehat L_j},\qquad
 \underline d_j=\frac{c_d}{4TQ_{*,j}}
 =\frac{c_d}{4T}e^{-(1-\alpha_j)L_j},
\end{equation}
where $K_f\geq1$ and $c_d>0$ will be fixed independently of the
initial scales. Then
\begin{equation}\label{eq:adjacentseparation}
 \frac{\overline d_{j+1}}{\underline d_j}
 =\frac{4K_f}{c_d}T^{7/2}
 e^{-\sigma\widehat L_{j+1}+(1-\alpha_j)L_j}.
\end{equation}
\begin{parametercondition}[Lower bound on scale growth]\label{cond:geometricgrowth}
Fix $C_{\mathrm{sep}}\geq 12/\sigma.$ For every $j\geq1$, we require
\[
 \widehat L_{j+1}\geq C_{\mathrm{sep}}L_j,\qquad
\]
\end{parametercondition}
We use this restriction to verify separation of comparison radii and
estimate the interaction errors below. For this verification, take
\begin{equation}\label{eq:comparison-time}
 T=T_{\mathrm{cmp}}:=\frac{\tau_*}{\lambda_{1,0}}
 =\tau_*e^{(1+\alpha_1)L_1-\kappa\widehat L_1}.
\end{equation}
This choice fixes the time bound in our estimates. We consider evolution on classical intervals $[0,t_0]$ with $t_0\leq T$.

\begin{lemma}[Verification of the comparison radii]\label{lem:comparisonchoice}
Under Setup~\ref{setup:parameters} {(an initial family)}, assume
\begin{itemize}
\item
\conditionref{cond:heightmargin}, \item \conditionref{cond:geometricgrowth}, 
\end{itemize}
and let $T=T_{\mathrm{cmp}}\geq1$. There exist constants
$c_d>0$, $K_f\geq1$, $\widetilde C_{\mathrm{int}}>0$, and $C_{\mathrm{geom}}>0$,
depending only on the profile $f_0$, cutoffs $\Phi, \chi^{\mathrm{out}}, \chi^{\mathrm{in}}$, and $\kappa$, such that
\[
 \widehat L_1\geq C_{\mathrm{geom}}
\]
implies that the constants in \eqref{eq:comparisonprefactors} satisfy
\begin{itemize}
\item \conditionref{cond:supportmargins}, 
\item
\conditionref{cond:prefactorseparation}, and 
\item \conditionref{cond:prefactorsummability}.
\end{itemize}
Moreover, for every $j\geq1$,
\begin{equation}\label{eq:integratederror}
 \mathcal I_j\leq \widetilde C_{\mathrm{int}}\left[
 e^{-\sigma\widehat L_j}\overline A_{j-1}
 +e^{-(2\kappa+3\sigma)\widehat L_{j+1}/2}\right].
\end{equation}
In particular, \conditionref{cond:errorcondition} follows from
\begin{equation}\label{eq:expliciterrorcondition}
 \widetilde C_{\mathrm{int}}\left[
 e^{-\sigma\widehat L_j}\overline A_{j-1}
 +e^{-(2\kappa+3\sigma)\widehat L_{j+1}/2}\right]\leq\frac12,
 \qquad j\geq1.
\end{equation}
\end{lemma}

\begin{proof}
\emph{Support margins.}
For $j\geq2$, \conditionref{cond:geometricgrowth} and
$\beta_{j-1}\leq1/2$ imply
\[
 \widehat L_{j-1}\leq\frac{\widehat L_j}{2C_{\mathrm{sep}}}.
\]
Since
\[
 \frac{1-\kappa}{8C_{\mathrm{sep}}}\leq\frac12,\qquad
 \frac{1-\kappa}{6C_{\mathrm{sep}}}\leq\sigma,
\]
the bounds \eqref{eq:upperradial}--\eqref{eq:uppervelocity} yield
\[
 R_{j,0}e^{3\overline A_{j-1}}
 \leq C_a^{3/4}e^{-\sigma\widehat L_j},\qquad
 V_{0,j}^{5/2}e^{4\overline A_{j-1}}
 \leq C_a e^{-\sigma\widehat L_j}.
\]
For $j=1$, the same bounds follow from $\overline A_0=0$,
$L_1\geq2\widehat L_1$, and
\[
 V_{0,1}^{5/2}
 =\exp\left(-(1+5\alpha_1)L_1/2-\kappa\widehat L_1/2\right)
 \leq \exp\left( -\sigma\widehat L_1\right).
\]
Since $\kappa>\sigma$ and $T\geq1$, for every $j\geq1$,
\[
 H_{j,0}+R_{j,0}e^{3\overline A_{j-1}}
 +T^{5/2}V_{0,j}^{5/2}e^{4\overline A_{j-1}}
 \leq(1+C_a^{3/4}+C_a)T^{5/2}e^{-\sigma\widehat L_j}.
\]
Choose
\[
 K_f\geq2C_{\mathrm{upp}}(1+C_a^{3/4}+C_a),\qquad
 0<c_d\leq\frac{e^{-1}}{(1+C_{\mathrm{low}})F}.
\]
The first choice implies the first inequality in \conditionref{cond:supportmargins}.
Increase $C_{\mathrm{geom}}$ so that
$F\zeta_-e^{2C_{\mathrm{geom}}/3}\geq1$.
Since $L_j\geq2\widehat L_j\geq2\widehat L_1$, we then have
$Q_j\rho_j^-(0)\geq1$ for every $j$.
The lower bound above and $Q_j=FQ_{*,j}$ now yields
the second inequality in \conditionref{cond:supportmargins}.

\smallskip
\noindent\emph{Separation.}
By \eqref{eq:comparison-time}, $T\leq\tau_*e^{4L_1/3}$.
Thus \eqref{eq:adjacentseparation} implies, for $j\geq2$,
\begin{align*}
 \frac{\overline d_j}{\underline d_{j-1}}
 &\leq\frac{4K_f}{c_d}\tau_*^{7/2}
 \exp\left(-\sigma\widehat L_j+(1-\alpha_{j-1})L_{j-1}+(14/3)L_1\right)\\
 &\leq\frac{4K_f}{c_d}\tau_*^{7/2}
 \exp\left(-[\sigma-11/(2C_{\mathrm{sep}})]\widehat L_j\right)
 \lesssim \exp\left(-\sigma\widehat L_j/2\right)\lesssim \exp\left(-\sigma C_{\mathrm geom}/2\right).
\end{align*}
Here we used
\[
 1-\alpha_{j-1}\leq\frac56,\qquad
 L_1\leq L_{j-1}\leq\frac{\widehat L_j}{C_{\mathrm{sep}}},\qquad
 \frac{11}{2C_{\mathrm{sep}}}\leq\frac{\sigma}{2}
\]
and the hypothesis $\hat L_1 \geq C_{\mathrm{geom}}$. 
Enlarging $C_{\mathrm{geom}}$ makes the ratio at most $1/32$
for every $j\geq2$. The verifiesupport margins also give
\[
 \underline d_j\leq\frac{\rho_j^-(0)}{4e}
 \leq\rho_j^+(0)\leq\frac12\overline d_j,
\]
which completes \conditionref{cond:prefactorseparation}.

\smallskip
\noindent\emph{Inner mass and interaction estimates.}
Let
\begin{equation}\label{eq:M}
 N_i:=m_i\overline d_i^3
 =K_f^3T^{15/2}
 \exp\left(-(1+\alpha_i)L_i-(\kappa+3\sigma)\widehat L_i\right).
\end{equation}
Since $L_{i+1}\geq2C_{\mathrm{sep}}L_i$ and
$\widehat L_{i+1}\geq\widehat L_i$,
\[
 \frac{N_{i+1}}{N_i}
 \leq \exp\left(-7L_{i+1}/6+4L_i/3\right)
 \leq \exp\left(-(7C_{\mathrm{sep}}-4)L_i/3 \right)
 \leq\frac12
\]
after a further enlargement of $C_{\mathrm{geom}}$. Hence
\begin{equation}\label{eq:innertail}
 \sum_{i\geq1}N_i<\infty,\qquad
 \sum_{i>j}N_i\leq2N_{j+1}.
\end{equation}
This establishes \conditionref{cond:prefactorsummability}. The separation estimate,
\eqref{eq:innertail}, and \eqref{eq:interactionbudget} imply
\begin{equation}\label{eq:adjacenterror}
 \mathcal I_j\lesssim e^{-\sigma\widehat L_j}\overline A_{j-1}
 +\frac{\tau_*N_{j+1}}{\lambda_{j+1,0}\underline d_j^5}.
\end{equation}
To bound the second term, note that
\[
 \frac{N_{j+1}}{\lambda_{j+1,0}}
 =K_f^3T^{15/2}\exp\left(-(2\kappa+3\sigma)\widehat L_{j+1}\right).
\]
Using \eqref{eq:comparisonprefactors} and
$T\leq\tau_*e^{4L_1/3}$, we obtain
\begin{align*}
 \frac{\tau_*N_{j+1}}{\lambda_{j+1,0}\underline d_j^5}
 &\lesssim\tau_*T^{25/2}
\exp\left(-(2\kappa+3\sigma)\widehat L_{j+1}+5(1-\alpha_j)L_j\right)\\
 &\lesssim\tau_*^{27/2}
\exp\left(-(2\kappa+3\sigma)\widehat L_{j+1}
       +5(1-\alpha_j)L_j+(50/3)L_1\right)\\
 &\lesssim \exp\left(-(2\kappa+3\sigma)\widehat L_{j+1}/2\right).
\end{align*}
Indeed, $C_{\mathrm{sep}}\geq12/\sigma
\geq125/[3(2\kappa+3\sigma)]$ implies
\[
 5(1-\alpha_j)L_j+\frac{50}{3}L_1
 \leq\frac{125}{6}L_j
 \leq\frac{125}{6C_{\mathrm{sep}}}\widehat L_{j+1}
 \leq\frac{2\kappa+3\sigma}{2}\widehat L_{j+1}.
\]
Combining this bound with \eqref{eq:adjacenterror} we obtain
\eqref{eq:integratederror}.
\end{proof}

}

{

\subsection{Choice of scales}
\label{sec:Choice_of_scales}

Combining \conditionref{cond:geometricgrowth} and
\conditionref{cond:nextscaleupper}, we require, for every $j\geq1$,
\[
 \frac{C_{\mathrm{sep}}L_j}{\beta_{j+1}}
 \leq L_{j+1}\leq\widehat L_j^2,
 \qquad j\geq1.
\]
The next lemma verifies the earlier requirements once the initial
scale $\widehat L_1$ is sufficiently large.

\begin{lemma}[Sufficient conditions on the logarithmic scales]
\label{lem:parameterconditions}
Fix the profile $f_0$, cutoffs $\Phi, \chi^{\mathrm{in}},\chi^{\mathrm{out}}$, $1/3<\kappa<1$, and
$C_{\mathrm{sep}}\geq12/\sigma$.
There exists $S_*>0$, depending only on the profile, cutoffs, and
$\kappa$, such that every choice of parameter sequences $(\alpha_j)_{j\geq 1}$ and $(L_j)_{j\geq 1}$ satisfying
\begin{itemize}
\item
\eqref{eq:regularityparameters}, 
\item $\widehat L_1\geq S_*$, 
\item \conditionref{cond:geometricgrowth}, and
\item
\conditionref{cond:nextscaleupper}
\end{itemize}
has the following properties:
\begin{enumerate}
\item \conditionref{cond:initialsummability} holds.
\item 
\begin{itemize}
\item
\conditionref{cond:capacitygap}, \item \conditionref{cond:smallreferenceaspect}, and \item 
\conditionref{cond:straincomparison} hold.
\end{itemize}
\item With $T=T_{\mathrm{cmp}}$, one has $T\geq1$ and
$\widehat L_1\geq C_{\mathrm{geom}}$. The constants
$\bar d_j,\underline d_j$ in \eqref{eq:comparisonprefactors}
satisfy 
\begin{itemize}
    \item \conditionref{cond:supportmargins},
    \item
\conditionref{cond:prefactorseparation},
    \item
\conditionref{cond:prefactorsummability},
and \item\conditionref{cond:errorcondition}.
\end{itemize}
\end{enumerate}
Moreover,
\begin{equation}\label{eq:coarsestrain}
 \overline A_j<\frac{L_{j+1}}3
 \leq\frac{\widehat L_j^2}{3},
 \qquad j\geq1.
\end{equation}
\end{lemma}

\begin{proof}
\noindent\emph{Part (i)}
Since $0<\beta_j\leq1/2$, \conditionref{cond:geometricgrowth} implies
\[
 L_{j+1}\geq2C_{\mathrm{sep}}L_j,\qquad
 \widehat L_{j+1}\geq C_{\mathrm{sep}}L_j
 \geq2C_{\mathrm{sep}}\widehat L_j.
\]
Hence
\[
 \widehat L_j\geq2^{j-1}\widehat L_1\to \infty, \qquad
 \sum_{j\geq1}e^{-\eta L_j}
 \leq\sum_{j\geq1}e^{-\eta2^{j-1}L_1}<\infty, \quad \eta>0,
\]
and hence (i) is satisfied. 

\smallskip
\noindent\emph{Part (ii)}
Taking $S_*\geq S_{\mathrm{asp}}$,
Lemma~\ref{lem:aspectstrainconditions} applies and establishes
\conditionref{cond:smallreferenceaspect} and
\conditionref{cond:straincomparison}.
Take also $S_*\geq C_{\mathrm{geom}}$ so large that
\[
 (1-\kappa)S_*\geq3\log2,\qquad
 \kappa S_*>\tfrac34\log C_a,\qquad
 \tau_*e^{4S_*/3}\geq1.
\]
Then
\[
 \begin{gathered}
 (1-\kappa)(\widehat L_j-\widehat L_{j-1})
 \geq(1-\kappa)S_*\geq3\log2, \qquad j\geq2,\\
 T_{\mathrm{cmp}} \geq\tau_*e^{2L_1/3}
 \geq\tau_*e^{4\widehat L_1/3}\geq 1.
 \end{gathered}
\]
Thus \conditionref{cond:capacitygap} holds, completing (ii).

\smallskip
\noindent\emph{Part (iii)}
Lemma~\ref{lem:comparisonchoice} establishes
\conditionref{cond:supportmargins}, \conditionref{cond:prefactorseparation}, \conditionref{cond:prefactorsummability},
and \eqref{eq:integratederror}.

Since $\widehat L_j\leq\widehat L_{j+1}$,
\eqref{eq:explicitstrain} yields
\[
 \overline A_j
 \leq\frac{L_{j+1}}3-\frac{\kappa\widehat L_{j+1}}3
      +\frac14\log C_a
 <\frac{L_{j+1}}3
 \leq\frac{\widehat L_j^2}{3},
\]
and thus \eqref{eq:coarsestrain} holds. In particular,
$\overline A_{j-1}\leq\widehat L_j^2/3$ for every $j\geq1$,
with $\overline A_0=0$.
Using \eqref{eq:integratederror} and
$\widehat L_{j+1}\geq\widehat L_j$, we obtain
\[
 \mathcal I_j\leq C(1+\widehat L_j^2)e^{-\sigma\widehat L_j},
 \qquad j\geq1,
\]
where $C$ depends only on the profile, cutoffs, and $\kappa$.
Enlarge $S_*$ so that
\begin{equation}\label{eq:Ystrain}
 C\sup_{s\geq S_*}(1+s^2)e^{-\sigma s}\leq\frac12.
\end{equation}
Since $\widehat L_j\geq S_*$, \conditionref{cond:errorcondition}
follows. This completes (iii).
\end{proof}

\begin{parametercondition}[Initial scale threshold]\label{cond:initialthreshold}
We require $\widehat L_1\geq S_*$, with $S_*$ as in
Lemma~\ref{lem:parameterconditions}.
\end{parametercondition}

Thus \conditionref{cond:nextscaleupper}, \conditionref{cond:heightmargin},
\conditionref{cond:geometricgrowth}, and \conditionref{cond:initialthreshold}
ensure Conditions~\ref{cond:initialsummability}--\ref{cond:errorcondition},
with the prescribed comparison radii and $T=T_{\mathrm{cmp}}$.

We now choose the explicit family to satisfy this threshold,
\conditionref{cond:geometricgrowth}, and \conditionref{cond:nextscaleupper}.
For the family in Theorem~\ref{thm:constructive}, take
\begin{equation}\label{eq:finalparameters}
 \alpha_j=\frac{j}{3(j+1)},\qquad
 \beta_j=\frac1{j+1},\qquad
 \kappa=\frac56,\qquad C_{\mathrm{sep}}=100.
\end{equation}
Then $\sigma=3/8$ and $12/\sigma=32<100$.
For $Y>0$, define
\begin{equation}\label{eq:explicitlogscales}
 L_j=Y100^{j-1}(j+1)!,\qquad
 \widehat L_j=Y100^{j-1}j!.
\end{equation}
Since $(j+1)(j+2)\leq6j!$ for $j\geq1$,
\begin{equation}\label{eq:scales}
 \begin{aligned}
 L_{j+1}
 &=\frac{C_{\mathrm{sep}}L_j}{\beta_{j+1}},\\
 \frac{L_{j+1}}{\widehat L_j^2}
 &=\frac{100^{2-j}(j+1)(j+2)}{Yj!}
 \leq\frac{600}{Y}.
 \end{aligned}
\end{equation}
Choose
\begin{equation}\label{eq:finalthreshold}
 Y_0:=\max\{S_*,600\}.
\end{equation}
For $Y\geq Y_0$, we have $\widehat L_1=Y\geq S_*$, and
\conditionref{cond:geometricgrowth} and \conditionref{cond:nextscaleupper} hold. Lemma~\ref{lem:parameterconditions}
therefore applies, and \eqref{eq:initialdata} gives the datum
of Theorem~\ref{thm:constructive}.
With the cutoffs fixed, $Y_0$ depends only on $f_0$.

}

\subsection{Classical evolution and finite-time blow-up}
\label{sec:comparisonargument}

\begin{figure}[!ht]
\centering
\begin{tikzpicture}[
 x=1cm,y=1cm,
 font=\small,text=black,
 box/.style={draw=black!75,fill=black!3,
   line width=.6pt,rounded corners=2pt,align=flush center,
   text width=12.2cm,inner sep=7pt},
 flow/.style={-{Stealth[length=2mm,width=1.5mm]},
   draw=black,line width=.6pt,shorten <=1pt,shorten >=1pt}
]
\node[box] (first) at (0,0)
 {\textbf{Suppose the first equality occurs at $t_{\rm e}$}\\[4pt]
 $\tau_j\leq\tau_*,\qquad\rho_j^+\leq\overline\rho_j,
 \qquad\rho_j^-\geq\underline\rho_j$\\[3pt]
 These bounds hold on $[0,t_{\rm e}]$ for every $j$.};
\node[box,below=7mm of first] (inputs)
 {\textbf{Strain, aspect, and interaction bounds}\\[4pt]
 $A_j\leq\overline A_j,\qquad\eps_j\leq\bar\eps_j,
 \qquad\mathfrak{E}_j\leq\tfrac12$\\[3pt]
 Lemmas~\ref{lem:clockbounds} and~\ref{lem:cascadeerrors}};
\draw[flow] (first.south)--(inputs.north);
\node[box,below=7mm of inputs] (supports)
 {\textbf{Strict support improvement}\\[4pt]
 $\rho_j^+\leq\tfrac12\overline\rho_j,\qquad
 \rho_j^-\geq4\underline\rho_j$\\[3pt]
 Lemma~\ref{lem:cascadeerrors}};
\draw[flow] (inputs.south)--(supports.north);
\node[box,below=7mm of supports] (threshold)
 {\textbf{Support equality is excluded}\\[4pt]
 Hence $\tau_k(t_{\rm e})=\tau_*$ for some $k$.};
\draw[flow] (supports.south)--(threshold.north);
\node[box,line width=.9pt,below=7mm of threshold] (contradiction)
 {\textbf{The large strain estimate contradicts the strain bound}\\[5pt]
 $\overline A_k<B_k(t_{\rm e})\leq A_k(t_{\rm e})\leq\overline A_k$\\[4pt]
 Theorem~\ref{thm:pulse} and the parameter conditions};
\draw[flow] (threshold.south)--(contradiction.north);
\end{tikzpicture}
\caption{The joint bootstrap in Proposition~\ref{prop:classicalcascade}.
Classical regularity reduces possible first equalities to finitely many
indices. At a hypothetical first equality time, the support estimates
exclude geometric failure, and the large strain estimate excludes a
pair reaching the threshold $\tau_*$. The argument takes place on an
existing classical interval and assumes no continuation to
$T_{\mathrm{cmp}}$.}
\label{fig:jointbootstrap}
\end{figure}
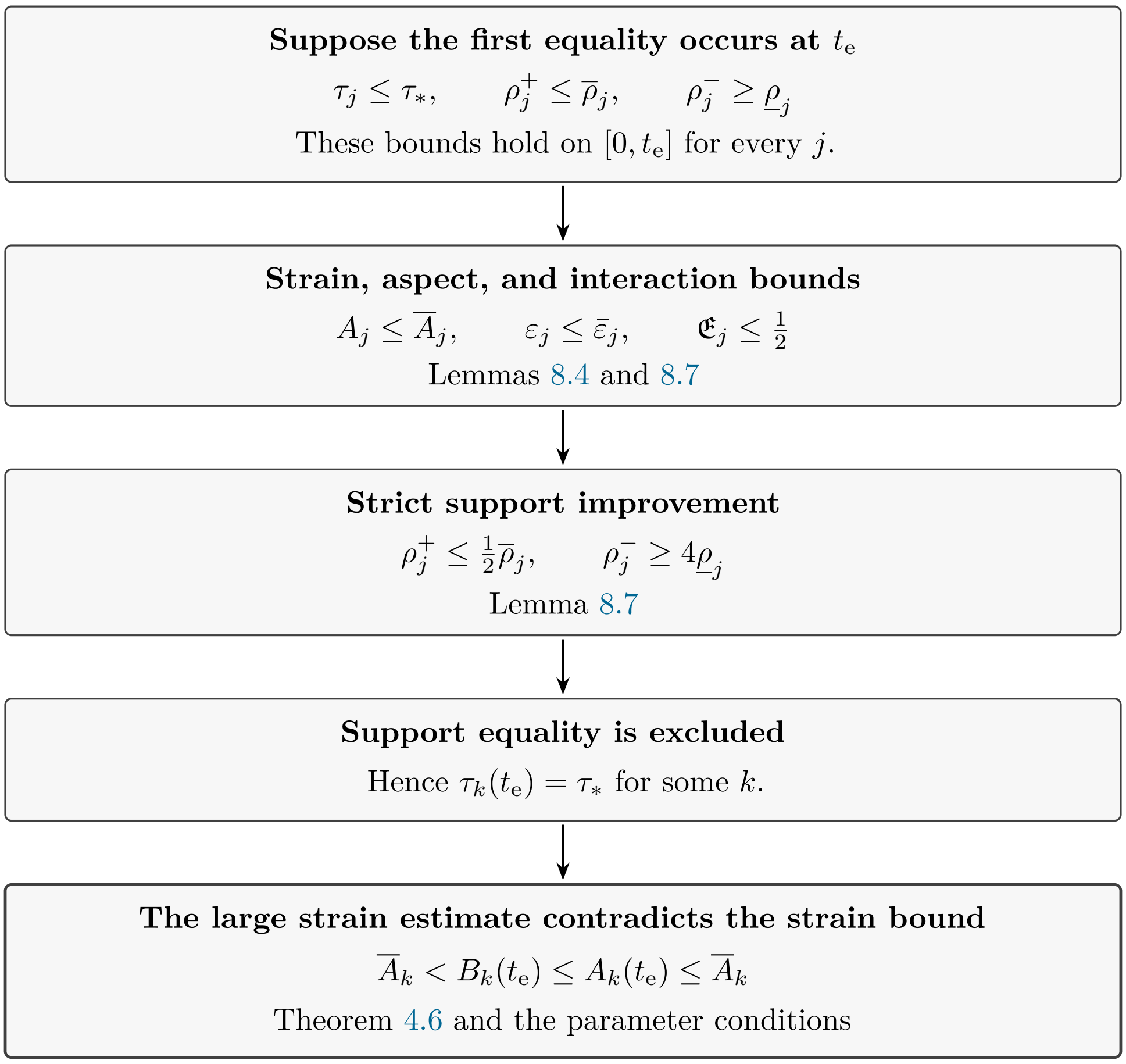

We propagate the comparison bounds on every classical interval
and then exclude classical existence up to $T_{\mathrm{cmp}}$, as summarized in
Figure~\ref{fig:jointbootstrap}.

\begin{proposition}[Propagation of the comparison bounds]
\label{prop:classicalcascade}
Under Setup~\ref{setup:cascade} (classical multiscale evolution), impose
\begin{itemize}
    \item 
\conditionref{cond:capacitygap},
    \item
\conditionref{cond:smallreferenceaspect}, and 
    \item \conditionref{cond:straincomparison}.
\end{itemize}
Take $T\geq\max\{1,t_0\}$ and
positive constants $\overline d_j,\underline d_j$ satisfying
\begin{itemize}
\item \conditionref{cond:supportmargins}, 
\item \conditionref{cond:prefactorseparation}, 
\item
\conditionref{cond:prefactorsummability}, and
\item
 \conditionref{cond:errorcondition}, 
\end{itemize}
and define the comparison radii $\bar \rho_j, \underline \rho_j$ by
\eqref{eq:comparisons}.
Then the error fields $\mathfrak{e}_j$ exist on $[0,t_0]$ and, for every
$j\geq1$ and $0\leq t\leq t_0$,
\begin{equation}\label{eq:classicalcascadebounds}
 \begin{gathered}
 \tau_j(t)<\tau_*,\qquad
 A_j(t)\leq\overline A_j,\qquad
 \eps_j(t)\leq\bar\eps_j,\\
 \rho_j^+(t)\leq\tfrac12\overline\rho_j(t),\qquad
 \rho_j^-(t)\geq4\underline\rho_j(t),\qquad
 \mathfrak{E}_j(t)\leq\tfrac12.
 \end{gathered}
\end{equation}
\end{proposition}

\begin{proof}
We prove that for every $t\in[0,t_0],$
\begin{equation}\label{eq:geometricbootstrap}
 \tau_j<\tau_*, \qquad \rho_j^+<\overline\rho_j,\qquad
 \rho_j^->\underline\rho_j, \qquad j\geq 1.
\end{equation}
These inequalities hold at $t=0$ by $\tau_j(0)=0$ and
Lemma~\ref{lem:cascadeerrors}.

Let $X_t$ be the flow of $u$. Without assuming
\eqref{eq:geometricbootstrap}, Lemma~\ref{lem:cascadestrain}
and $u(t,0)=0$ yield, for $0\leq t\leq t_0$,
\[
 0\leq A_j(t)\leq A_{\mathrm{Lip}}(t_0),\qquad
 e^{-A_{\mathrm{Lip}}(t_0)}|a|
 \leq|X_t(a)|\leq e^{A_{\mathrm{Lip}}(t_0)}|a|.
\]
For $j\geq2$,
\[
 \frac{V_{0,j}^{5/2}}{\lambda_{j,0}H_{j,0}}=P_j^{3/2},
 \qquad
 e^{4\overline A_{j-1}}
 =1+C_0\frac{P_{j-1}}{\lambda_{j,0}}.
\]
Consequently, \conditionref{cond:supportmargins} implies
\begin{align*}
 \frac{\overline d_j}{H_{j,0}}
 &\geq2C_{\mathrm{upp}}T^{5/2}C_0P_{j-1}P_j^{3/2}
  \to \infty,\\
 \frac{\underline d_j}{\rho_j^-(0)}
 &\leq
 \frac{e^{-1}}{4\bigl(1+C_{\mathrm{low}}Q_j\rho_j^-(0)T\bigr)}
  \to 0,
\end{align*}
since
\[
 P_j\geq2^{j-1}P_1 \to \infty,\qquad
 Q_j\rho_j^-(0)\geq F\zeta_-e^{L_j/3} \to \infty.
\]
Using $1+\alpha_j-\kappa\beta_j\geq4\alpha_j\geq2/3$ and
$\rho_j^+(0)\lesssim H_{j,0}$, we obtain, as $j\to\infty$,
\begin{align}
 \sup_{0\leq t\leq t_0}\tau_j(t)
 &\leq t_0e^{4A_{\mathrm{Lip}}(t_0)}\lambda_{j,0}
 \leq t_0e^{4A_{\mathrm{Lip}}(t_0)-2L_j/3}
  \to 0,\label{eq:tailclock}\\
 \sup_{0\leq t\leq t_0}
 \frac{\rho_j^+(t)}{\overline\rho_j(t)}
 &\lesssim e^{3A_{\mathrm{Lip}}(t_0)}
 \frac{H_{j,0}}{\overline d_j}
  \to 0,\label{eq:tailupper}\\
 \sup_{0\leq t\leq t_0}
 \frac{\underline\rho_j(t)}{\rho_j^-(t)}
 &\leq e^{A_{\mathrm{Lip}}(t_0)}
 \frac{\underline d_j}{\rho_j^-(0)}
  \to 0.\label{eq:taillower}
\end{align}
Thus \eqref{eq:geometricbootstrap} holds strictly on $[0,t_0]$
for all $j>J$, for some $J$. All quantities are continuous;
for $\rho_j^\pm$, this follows from the continuous flow of the
compact initial supports. If \eqref{eq:geometricbootstrap}
fails, there is therefore a first equality time
$t_{\rm e}\in(0,t_0]$, attained at some $j\leq J$.

On $[0,t_{\rm e}]$, \eqref{eq:geometricbootstrap} holds
non-strictly. Lemmas~\ref{lem:clockbounds}
and~\ref{lem:cascadeerrors} give the error fields and
\[
 \begin{gathered}
 A_j\leq\overline A_j,\qquad
 \eps_j\leq\bar\eps_j,\qquad \mathfrak{E}_j\leq\tfrac12,\\
 \rho_j^+\leq\tfrac12\overline\rho_j<\overline\rho_j,\qquad
 \rho_j^-\geq4\underline\rho_j>\underline\rho_j.
 \end{gathered}
\]
Hence $\tau_k(t_{\rm e})=\tau_*$ for some $1\leq k\leq J$.
Lemma~\ref{lem:cascadepulse} yields
\[
 \overline A_k<B_k(t_{\rm e})
 \leq A_k(t_{\rm e})\leq\overline A_k,
\]
a contradiction. Thus \eqref{eq:geometricbootstrap} holds on
$[0,t_0]$. Applying Lemmas~\ref{lem:clockbounds}
and~\ref{lem:cascadeerrors} on this time interval, we obtain
\eqref{eq:classicalcascadebounds} and the existence of the
error fields.
\end{proof}

\begin{theorem}[Finite comparison-time obstruction]
\label{thm:cascadeobstruction}
Let the initial parameters and comparison prefactors $\overline d_j, \underline d_j$ satisfy
the parameter hypotheses of Proposition~\ref{prop:classicalcascade}
with $T=T_{\mathrm{cmp}}\geq1$.
The solution from \eqref{eq:initialdata} cannot remain classical on
$[0,T_{\mathrm{cmp}}]$. In particular, this holds for every
$Y\geq Y_0$ in the explicit family
\eqref{eq:finalparameters}--\eqref{eq:finalthreshold}
and \eqref{eq:explicitlogscales}.
\end{theorem}

\begin{proof}
Suppose that $u$ is classical on $[0,T_{\mathrm{cmp}}]$.
By Proposition~\ref{prop:classicalcascade}, applied with
$t_0=T_{\mathrm{cmp}}$, and
\eqref{eq:firstclock}, \eqref{eq:comparison-time}, we have
\[
 \tau_*=\lambda_{1,0}T_{\mathrm{cmp}}
 =\tau_1(T_{\mathrm{cmp}})<\tau_*,
\]
a contradiction. For the explicit family, the parameter
hypotheses follow from Lemma~\ref{lem:parameterconditions}.
\end{proof}

{

\section{Initial regularity and proof of the main theorem}
\label{sec:limit}

We first verify the regularity and summability of the datum
\eqref{eq:initialdata}. Write
\begin{equation}\label{eq:infinitedatum}
 \vectomega_0=\sum_{j\geq1}\vectomega_{j,0},
 \qquad u_0=\BS[\vectomega_0].
\end{equation}

\begin{lemma}[Initial regularity]\label{lem:initialregularity}
Under Setup~\ref{setup:parameters}, the vorticity $\vectomega_0$
is compactly supported, divergence-free, and axisymmetric, 
and $u_0$ is divergence-free and axisymmetric without swirl.
Moreover,
\[
 \begin{gathered}
 \vectomega_0\in
 \bigcap_{0<\alpha<1/3}C^\alpha(\R^3)
 \cap C^\infty(\R^3\setminus\{0\}),\\
 u_0\in L^2(\R^3)\cap
 \bigcap_{0<\alpha<1/3}C^{1,\alpha}(\R^3)
 \cap C^\infty(\R^3\setminus\{0\}),
 \end{gathered}
\]
and, for every $0<\alpha<1/3$,
\[
 \sum_{j\geq1}\left(
 \|\vectomega_{j,0}\|_{C^\alpha} +\|u_j(0)\|_{L^2}
 +\|u_j(0)\|_{C^{1,\alpha}}\right)<\infty.
\]
\end{lemma}

\begin{proof}
Fix $0<\alpha<1/3$. Since
$Q_{*,j}=R_{j,0}^{\alpha_j-1}$ and $0<R_{j,0}\leq H_{j,0}<1$,
\[
 \|\vectomega_{j,0}\|_{L^\infty}
 \lesssim_{f_0}R_{j,0}^{\alpha_j}, \qquad
 \|\nabla\vectomega_{j,0}\|_{L^\infty}
 \lesssim_{f_0}R_{j,0}^{\alpha_j-1}.
\]
By interpolation, we have
\begin{equation}\label{eq:holder}
 \|\vectomega_{j,0}\|_{C^\alpha}
 \lesssim_{f_0,\alpha}R_{j,0}^{\alpha_j-\alpha}.
\end{equation}
Choose $J_\alpha$ so that
$\alpha_j-\alpha\geq\frac12(1/3-\alpha)$ for $j\geq J_\alpha$.
Then \conditionref{cond:initialsummability} yields
\[
 \sum_{j\geq J_\alpha}\|\vectomega_{j,0}\|_{C^\alpha}
 \lesssim_{f_0,\alpha} \sum_{j\geq J_\alpha}
 e^{-\frac12(1/3-\alpha)L_j} < \infty.
\]
Also,
\[
 \supp\vectomega_{j,0}\subset \overline{B_{C_{f_0}H_{j,0}}(0)},\qquad
 H_{j,0}=e^{-\kappa\widehat L_j}\to 0.
\]
Thus the sum is compactly supported, and every compact subset of
$\R^3\setminus\{0\}$ intersects only finitely many initial supports.
Thus $\vectomega_0\in C^\infty(\R^3\setminus\{0\})$.
The divergence and symmetry properties follow from the azimuthal
form of the summands and the Biot--Savart law.

By Lemma~\ref{lem:normalizedenergy} at $t=0$,
\[
 \|u_j(0)\|_{L^2}
 \lesssim_{f_0}
 Q_{*,j}R_{j,0}^3H_{j,0}^{1/2}
 =R_{j,0}^{\alpha_j+2}H_{j,0}^{1/2}
 \leq e^{-13L_j/6},
\]
where $\alpha_j\geq1/6$ and $H_{j,0}<1$.
Hence \conditionref{cond:initialsummability} implies $\sum_j\|u_j(0)\|_{L^2}<\infty$.
The inhomogeneous Biot--Savart estimate (see \cite[Theorem~2.4.2]{Che95}) is
\begin{equation}\label{eq:schauder}
 \|v\|_{C^{1,\alpha}} \lesssim_\alpha
 \|v\|_{L^2}+\|\curl v\|_{C^\alpha},
 \qquad \Div v=0.
\end{equation}
Applying it to each $u_j(0)$, we obtain
\[
 \sum_{j\geq1}
 \left(\|u_j(0)\|_{L^2} +\|u_j(0)\|_{C^{1,\alpha}}\right)
 \lesssim_\alpha \sum_{j\geq1}
 \left(\|u_j(0)\|_{L^2}
       +\|\vectomega_{j,0}\|_{C^\alpha}\right)<\infty.
\]
Finally, $-\Delta u_0=\curl\vectomega_0$ and local elliptic
regularity imply $u_0\in C^\infty(\R^3\setminus\{0\})$.
\end{proof}

\begin{proof}[Proof of Theorems~\ref{thm:main}
and~\ref{thm:constructive}]
Fix the datum \eqref{eq:constructive-data}--\eqref{eq:constructive-scales}
with $Y\geq Y_0(f_0)$, as in
\eqref{eq:finalparameters}--\eqref{eq:finalthreshold}.
Lemma~\ref{lem:parameterconditions} verifies the parameter
hypotheses of Theorem~\ref{thm:cascadeobstruction}, and
Lemma~\ref{lem:initialregularity} gives the stated initial
regularity and absolute convergence. In addition,
\[
 \|\nabla u_0\|_{L^2}
 =\|\vectomega_0\|_{L^2}<\infty.
\]
By classical local existence in $C^{1,1/6}(\R^3)\cap L^2(\R^3)$, there exists
a maximal solution $u$ with lifespan $T_*>0$.
Axisymmetry, absence of swirl, and parity are preserved.

If $T_*>T_{\mathrm{cmp}}$, then $u$ is classical on
$[0,T_{\mathrm{cmp}}]$, contradicting
Theorem~\ref{thm:cascadeobstruction}. Thus
\[
 0<T_*\leq T_{\mathrm{cmp}}<\infty.
\]
By the Beale--Kato--Majda continuation criterion in this
H\"older class \cite[Theorem~4.2.3]{Che95}, we have
\[
 \int_0^{T_*}\|\vectomega(t)\|_{L^\infty(\R^3)}\dd t=\infty.
\]
\end{proof}

\subsection{Smoothness away from the origin and localization of blow-up}
\label{subsec:blowuplocation-proof}

\begin{proof}[Proof of Remark~\ref{rem:blowuplocation}]
Fix $0<t_0<T_*$ and $0<\alpha<1/6$. Since $u$ is classical, $u \in L^\infty(0,t_0;C^{1,\alpha}(\R^3))$.
Let $X(t,a)=X_t(a)$ be the flow map generated by $u$.
By symmetry, $u(t,0)=0$, and hence $X_t(0)=0$.
Continuity and invertibility imply, for every compact
$\mathfrak{K}\subset\R^3\setminus\{0\}$, the set
\[
 \bigcup_{0\leq t\leq t_0} \bigl(X_t(\mathfrak{K})\cup X_t^{-1}(\mathfrak{K})\bigr)
\]
is compact and does not contain the origin.
Each such compact set intersects only finitely many initial supports of $q_{j,0}$.
Thus $q_0=\sum_jq_{j,0}\in C^\infty(\R^3\setminus\{0\})$, and
\[
 q(t)=q_0\circ X_t^{-1},\qquad \vectomega(t,x)=q(t,x)(x_2,-x_1,0).
\]
Let $V$ be a bounded open set and $\mathfrak{K}$ be compact such that 
$\mathfrak{K} \subset V \subset \overline V \subset \R^3\setminus\{0\}$.
For $k\geq1$, the interior Schauder estimate for $-\Delta u=\curl\vectomega$ yields
\[
 \|u(t)\|_{C^{k+1,\alpha}(\mathfrak{K})}
 \leq C_{\mathfrak{K},V,k,\alpha}\bigl(
 \|u(t)\|_{L^\infty(V)}
 +\|q(t)\|_{C^{k,\alpha}(\overline V)}\bigr).
\]
Since $u$ is classical, 
$X_t,X_t^{-1}\in C^{1,\alpha}_{\rm loc}$, so the transport
formula yields $q(t)\in C^{1,\alpha}_{\rm loc}$,
uniformly on compact subsets away from the origin.
Using this estimate with
the flow equation and the transport formula $q(t,x)=q_0(X_t^{-1}(x))$, we obtain
\[
 q\in C^{k,\alpha}_{\rm loc}
  \to  u\in C^{k+1,\alpha}_{\rm loc}
  \to  X_t,X_t^{-1}\in C^{k+1,\alpha}_{\rm loc}
  \to  q\in C^{k+1,\alpha}_{\rm loc},
\]
with all H\"older norms bounded uniformly for $0\leq t\leq t_0$
on every compact subset of $\R^3\setminus\{0\}$.
Starting with $k=1$, we get a bound on every
spatial derivative and obtain smoothness in space. By Lagrangian analyticity
\cite[Theorem~1.1]{CVW15}, every pure time derivative of $X$ is also
locally bounded. For a cutoff function
$\chi\in C_c^\infty((0,t_0)\times(\R^3\setminus\{0\}))$, the Fourier
estimate
\[
 \| \chi X \|_{H^m(\R^4)} \lesssim_m\| \chi X \|_2+\|\partial_t^m (\chi X)\|_2 + \sum_{i=1}^3\|\partial_{a_i}^m (\chi X)\|_2<\infty, \qquad m \in \mathbb{N},
\]
implies $X\in C^\infty((0,t_0)\times(\R^3\setminus\{0\}))$.
Then using the inverse function theorem and
$u(t,x)=\partial_tX(t,X_t^{-1}(x))$, we obtain the space-time smoothness of $u$.

By Proposition~\ref{prop:classicalcascade},
\eqref{eq:comparisons}, and \conditionref{cond:prefactorseparation},
\[
 \supp\vectomega_j(t)\subset\overline{B_{\bar d_j/2}(0)},
 \qquad \lim _{j\to \infty} \bar d_j = 0,\qquad 0\leq t<T_*.
\]
Choose $J$ such that $\bar d_j<2\delta$ for $j>J$.
Since $\|q_j(t)\|_\infty=\|q_{j,0}\|_\infty$ and
$|\vectomega_j|=r|q_j|$,
\[
 \sup_{0\leq t<T_*}
 \|\vectomega(t)\|_{L^\infty(\R^3\setminus B_\delta(0))}
 \leq \frac12\sum_{j=1}^{J}\|q_{j,0}\|_\infty \bar d_j < \infty.
\]
Hence, Theorem~\ref{thm:main} implies that the integral of $\|\vectomega(t)\|_{L^\infty(B_\delta(0))}$ diverges.
\end{proof}

}

\subsection*{AI Acknowledgment}
This project started two years ago. ChatGPT was used to optimize arguments and speed up the preparation of the first draft. The final draft was prepared and checked by the authors.

\end{document}